\documentclass[12pt]{article}%
\usepackage{amsmath}
\usepackage{amsfonts}
\usepackage{amssymb}
\usepackage{graphicx}
\usepackage{color}%
\providecommand{\U}[1]{\protect\rule{.1in}{.1in}}
\newtheorem{theorem}{Theorem}[section]
\newtheorem{acknowledgement}[theorem]{Acknowledgement}

\newtheorem{corollary}[theorem]{Corollary}

\newtheorem{definition}[theorem]{Definition}

\newtheorem{lemma}[theorem]{Lemma}
\newtheorem{notation}[theorem]{Notation}

\newtheorem{proposition}[theorem]{Proposition}
\newtheorem{remark}[theorem]{Remark}

\newenvironment{proof}[1][Proof]{\textbf{#1.} }{\ \rule{0.5em}{0.5em}}

\makeatletter\@addtoreset{equation}{section}\makeatother
\evensidemargin=\oddsidemargin
\newdimen\dummy
\dummy=\oddsidemargin
\begin{document}

\title{The inf-sup constant for $hp$-Crouzeix-Raviart triangular elements}
\author{S. Sauter\thanks{(stas@math.uzh.ch), Institut f\"{u}r Mathematik,
Universit\"{a}t Z\"{u}rich, Winterthurerstr 190, CH-8057 Z\"{u}rich,
Switzerland}}
\maketitle

\begin{abstract}
In this paper, we consider the discretization of the two-dimensional
stationary Stokes equation by Crouzeix-Raviart elements for the velocity of
polynomial order $k\geq1$ on conforming triangulations and discontinuous
pressure approximations of order $k-1$. We will bound the inf-sup constant
from below independent of the mesh size and show that it depends only
logarithmically on $k$. Our assumptions on the mesh are very mild: for odd $k$
we require that the triangulations contain at least one inner vertex while for
even $k$ we assume that the triangulations consist of more than a single triangle.

\end{abstract}



\section{Introduction}

In this paper we consider the numerical discretization of the two-dimensional
stationary Stokes problem by Crouzeix-Raviart elements. They were introduced
in the seminal paper \cite{CrouzeixRaviart} in 1973 by Crouzeix and Raviart
with the goal to obtain a stable and economic discretization of the Stokes
equation. They can be considered as an non-conforming enrichment of conforming
finite elements of polynomial degree $k$ for the velocity and discontinuous
pressures of degree $k-1$. It is well known that the conforming $\left(
k,k-1\right)  $ pair of finite elements can be unstable; for two-dimensions
the proof of the inf-sup stability of Crouzeix-Raviart discretizations of
general order $k$ has been evolved over the last 50 years, the inf-sup
stability for $k=1$ has been proved in \cite{CrouzeixRaviart} and only
recently the last open case $k=3$, has been proved in \cite{CCSS_CR_2}. We
mention the papers \cite{Fortin_Soulie}, \cite{ScottVogelius},
\cite{Crouzeix_Falk}, \cite{Baran_Stoyan}, \cite{GuzmanScott2019},
\cite{CCSS_CR_1} which contain essential milestones in this development. There
is a vast of literature on various further aspects of Crouzeix-Raviart
elements; we omit to present a comprehensive review here but refer to the
overview article \cite{Brenner_Crouzeix} instead.

Since higher order methods are becoming increasingly popular a natural
question arises how the inf-sup constant depends on the polynomial degree $k$.
It is the goal of this paper to investigate this dependency.

The paper is organized as follows. In Section \ref{NuDiscrete} we introduce
the Stokes problem and the Crouzeix-Raviart discretization of polynomial order
$k$. We state our main theorem that the discrete inf-sup constant can be
estimated from below by $c\left(  \log\left(  k+1\right)  \right)  ^{-\alpha}%
$, where the positive constant $c$ depends only on the shape-regularity of the
mesh and on the maximal outer angle of the domain $\Omega$. The explicit value
of $\alpha$ depends on the mesh topology. The simplest case is that each
triangle in triangulation contains at least one inner vertex and then
$\alpha=1/2$ holds. We will give the value of $\alpha$ also for more general triangulations.

The proof is given in Section \ref{SecProof}. The key ingredient is to show
that for any discrete pressure $q$, there exists a velocity field
$\mathbf{v}_{q}$ from the Crouzeix-Raviart space such that $\tilde
{q}:=q-\operatorname{div}_{\mathcal{T}}\mathbf{v}_{q}$ belongs to the
\textit{Scott-Vogelius pressure space} \cite[R.1, R.2]{vogelius1983right},
\cite[R.1, R.2]{ScottVogelius} (see (\ref{defMSV})) and $\mathbf{v}_{q}$
depends continuously on $q$. These properties allow us to construct a
conforming velocity field $\mathbf{\tilde{v}}_{\tilde{q}}$ of order $k$ with
$\operatorname*{div}\mathbf{\tilde{v}}_{\tilde{q}}=\tilde{q}$. For this step,
the construction in \cite[Thm. 5.1]{ScottVogelius}, \cite[Thm. 1]%
{GuzmanScott2019} is modified such that the norm of the right-inverse does not
deteriorate if a triangle vertex is a \textquotedblleft
nearly-critical\textquotedblright\ point -- a notion which will be introduced
in Definition \ref{Defetacriticalpoint}. This key result is proved for odd
polynomial degree in Section \ref{PTpge5} and for even polynomial degree in
Section \ref{SecEvenk}.

In the conclusions (Sec. \ref{SecConcl}) we summarize our main findings and
compare our results with existing results in the literature on some other
pairs of finite elements for the Stokes equation.

In the appendices, we prove a technical result for a Gram matrix related to
the bilinear form $\left(  \operatorname*{div}_{\mathcal{T}}\cdot
,\cdot\right)  _{L^{2}\left(  \Omega\right)  }$ applied to the
Crouzeix-Raviart element and discontinuous pressure space (see \S \ref{AppT_L}%
), an estimate of traces of non-conforming Crouzeix-Raviart basis functions
(see \S \ref{AppH1/2Est}), give some explicit formulae for integrals related
to orthogonal polynomials (see \S \ref{LegPoly}), and prove a discrete
Friedrichs inequality for Crouzeix-Raviart spaces (see \S \ref{SecNormEq}).

\subsection{The Stokes problem and its numerical
discretization\label{NuDiscrete}}

Let $\Omega\subset\mathbb{R}^{2}$ denote a bounded polygonal Lipschitz domain
with boundary $\partial\Omega$. For a vertex $\mathbf{z}$ in $\partial\Omega$,
we denote by $\alpha_{\mathbf{z}}$ the exterior angle between the two segments
in $\partial\Omega$ with joint $\mathbf{z}$. The minimal outer angle at the
boundary vertices is given by%
\begin{equation}
\alpha_{\Omega}:=\min_{\mathbf{z}\text{ is a vertex in }\partial\Omega}%
\alpha_{\mathbf{z}} \label{defalphatau}%
\end{equation}
and satisfies $0<\alpha_{\Omega}<2\pi$ since $\Omega$ is Lipschitz. We
consider the Stokes equation%
\[%
\begin{array}
[c]{llll}%
-\Delta\mathbf{u} & -\nabla p & =\mathbf{f} & \text{in }\Omega,\\
\operatorname*{div}\mathbf{u} &  & =0 & \text{in }\Omega
\end{array}
\]
with Dirichlet boundary conditions for the velocity and a usual normalization
condition for the pressure%
\[
\mathbf{u}=\mathbf{0}\quad\text{on }\partial\Omega\quad\text{and\quad}%
\int_{\Omega}p=0.
\]
To formulate this equation in a variational form we first introduce the
relevant function spaces. Throughout the paper we restrict to vector spaces
over the field of real numbers.

For $s\geq0$, $1\leq p\leq\infty$, $W^{s,p}\left(  \Omega\right)  $ denote the
classical Sobolev spaces of functions with norm $\left\Vert \cdot\right\Vert
_{W^{s,p}\left(  \Omega\right)  }$. As usual we write $L^{p}\left(
\Omega\right)  $ instead of $W^{0,p}\left(  \Omega\right)  $ and $H^{s}\left(
\Omega\right)  $ for $W^{s,2}\left(  \Omega\right)  $. For $s\geq0$, we denote
by $H_{0}^{s}\left(  \Omega\right)  $ the closure of the space of infinitely
smooth functions with compact support in $\Omega$ with respect to the
$H^{s}\left(  \Omega\right)  $ norm. Its dual space is denoted by
$H^{-s}\left(  \Omega\right)  $. For the pressure $p$, the space $L_{0}%
^{2}\left(  \Omega\right)  :=\left\{  u\in L^{2}\left(  \Omega\right)
:\int_{\Omega}u=0\right\}  $ will be relevant.

The scalar product and norm in $L^{2}\left(  \Omega\right)  $ are denoted
respectively by%
\[%
\begin{array}
[c]{llll}%
\left(  u,v\right)  _{L^{2}\left(  \Omega\right)  }:=\int_{\Omega}uv &
\text{and} & \left\Vert u\right\Vert _{L^{2}\left(  \Omega\right)  }:=\left(
u,u\right)  _{L^{2}\left(  \Omega\right)  }^{1/2} & \text{in }L^{2}\left(
\Omega\right)  .
\end{array}
\]
Vector-valued and $2\times2$ tensor-valued analogues of the function spaces
are denoted by bold and blackboard bold letters, e.g., $\mathbf{H}^{s}\left(
\Omega\right)  =\left(  H^{s}\left(  \Omega\right)  \right)  ^{2}$ and
$\mathbb{H}^{s}=\left(  H^{s}\left(  \Omega\right)  \right)  ^{2\times2}$.

The $\mathbf{L}^{2}\left(  \Omega\right)  $ scalar product and norm for vector
valued functions are given by%
\[
\left(  \mathbf{u},\mathbf{v}\right)  _{\mathbf{L}^{2}\left(  \Omega\right)
}:=\int_{\Omega}\left\langle \mathbf{u},\mathbf{v}\right\rangle \quad
\text{and\quad}\left\Vert \mathbf{u}\right\Vert _{\mathbf{L}^{2}\left(
\Omega\right)  }:=\left(  \mathbf{u},\mathbf{u}\right)  _{\mathbf{L}%
^{2}\left(  \Omega\right)  }^{1/2},
\]
where $\left\langle \mathbf{u},\mathbf{v}\right\rangle $ denotes the Euclidean
scalar product in $\mathbb{R}^{3}$. In a similar fashion, we define for
$\mathbf{G},\mathbf{H}\in\mathbb{L}^{2}\left(  \Omega\right)  $ the scalar
product and norm by%
\[
\left(  \mathbf{G},\mathbf{H}\right)  _{\mathbb{L}^{2}\left(  \Omega\right)
}:=\int_{\Omega}\left\langle \mathbf{G},\mathbf{H}\right\rangle \quad
\text{and\quad}\left\Vert \mathbf{G}\right\Vert _{\mathbb{L}^{2}\left(
\Omega\right)  }:=\left(  \mathbf{G},\mathbf{G}\right)  _{\mathbb{L}%
^{2}\left(  \Omega\right)  }^{1/2},
\]
where $\left\langle \mathbf{G},\mathbf{H}\right\rangle =\sum_{i,j=1}%
^{3}G_{i,j}H_{i,j}$. We also need fractional order Sobolev norms on boundary
of triangles and introduce the relevant notation; for details see, e.g.,
\cite{Mclean00}. For a bounded Lipschitz domain $\omega\subset\mathbb{R}^{2}$
with boundary $\partial\omega$, let $L^{2}\left(  \partial\omega\right)  $ and
$H^{1}\left(  \partial\omega\right)  $ denote the usual Lebesgue and Sobolev
space on $\partial\omega$ with norm $\left\Vert \cdot\right\Vert
_{L^{2}\left(  \partial\omega\right)  }$ and $\left\Vert \cdot\right\Vert
_{H^{1}\left(  \partial\omega\right)  }$. For $0<s<1$ the fractional Sobolev
space on $\partial\omega$ of order $s$ is denoted by $H^{s}\left(
\partial\omega\right)  $ and equipped with the norm%
\[
\left\Vert v\right\Vert _{H^{s}\left(  \partial\omega\right)  }:=\left(
\left\Vert v\right\Vert _{L^{2}\left(  \partial\omega\right)  }^{2}+\left\vert
v\right\vert _{H^{s}\left(  \partial\omega\right)  }^{2}\right)  ^{1/2}%
\]
and seminorm%
\[
\left\vert v\right\vert _{H^{s}\left(  \partial\omega\right)  }:=\left(
\int_{\partial\omega}\int_{\partial\omega}\frac{\left\vert v\left(
\mathbf{x}\right)  -v\left(  \mathbf{y}\right)  \right\vert ^{2}}{\left\Vert
\mathbf{x}-\mathbf{y}\right\Vert ^{1+2s}}d\mathbf{y}d\mathbf{x}\right)
^{1/2}.
\]

We introduce the bilinear form $a:\mathbf{H}^{1}\left(  \Omega\right)
\times\mathbf{H}^{1}\left(  \Omega\right)  \rightarrow\mathbb{R}$ by%
\begin{equation}
a\left(  \mathbf{u},\mathbf{v}\right)  :=\left(  \nabla\mathbf{u}%
,\nabla\mathbf{v}\right)  _{\mathbb{L}^{2}\left(  \Omega\right)  },
\label{defabili}%
\end{equation}
where $\nabla\mathbf{u}$ and $\nabla\mathbf{v}$ denote the derivatives of
$\mathbf{u}$ and $\mathbf{v}$. The variational form of the Stokes problem is
given by: For given $\mathbf{F}\in\mathbf{H}^{-1}\left(  \Omega\right)  ,$%
\begin{equation}
\text{find }\left(  \mathbf{u},p\right)  \in\mathbf{H}_{0}^{1}\left(
\Omega\right)  \times L_{0}^{2}\left(  \Omega\right)  \;\text{s.t.\ }\left\{
\begin{array}
[c]{lll}%
a\left(  \mathbf{u},\mathbf{v}\right)  +\left(  p,\operatorname*{div}%
\mathbf{v}\right)  _{L^{2}\left(  \Omega\right)  } & =\mathbf{F}\left(
\mathbf{v}\right)  & \forall\mathbf{v}\in\mathbf{H}_{0}^{1}\left(
\Omega\right)  ,\\
\left(  \operatorname*{div}\mathbf{u},q\right)  _{L^{2}\left(  \Omega\right)
} & =0 & \forall q\in L_{0}^{2}\left(  \Omega\right)  .
\end{array}
\right.  \label{varproblemstokes}%
\end{equation}

It is well-known (see, e.g., \cite{Girault86}) that (\ref{varproblemstokes})
is well posed. Since we consider non-conforming discretizations we restrict
the space $\mathbf{H}^{-1}\left(  \Omega\right)  $ for the right-hand side to
a smaller space and assume from now on for simplicity that $\mathbf{F}\left(
\mathbf{v}\right)  =\left(  \mathbf{f},\mathbf{v}\right)  _{\mathbf{L}%
^{2}\left(  \Omega\right)  }$ for some $\mathbf{f}\in\mathbf{L}^{2}\left(
\Omega\right)  $; for a more general setting we refer to
\cite{Veeser_qo_part1}, \cite{Veeser_qo}.\bigskip

In the following a discretization for problem (\ref{varproblemstokes}) is
introduced. Let $\mathcal{T}=\left\{  K_{i}:1\leq i\leq n\right\}  $ denote a
triangulation of $\Omega$ consisting of closed triangles which are
\textit{conforming}: the intersection of two different triangles is either
empty, a common edge, or a common point. We also assume $\Omega
=\operatorname*{dom}\mathcal{T}$, where
\begin{equation}
\operatorname*{dom}\mathcal{T}:=\operatorname*{int}\left(
%TCIMACRO{\dbigcup \limits_{K\in\mathcal{T}}}%
%BeginExpansion
{\displaystyle\bigcup\limits_{K\in\mathcal{T}}}
%EndExpansion
K\right)  \label{defdomT}%
\end{equation}
and $\operatorname*{int}\left(  M\right)  :=\overset{\circ}{M}$ denotes the
interior of a set $M\subset\mathbb{R}^{2}$.

Piecewise versions of differential operators\ such as $\nabla$ and
$\operatorname*{div}$ are defined for functions $u$ and vector fields
$\mathbf{w}$ which are sufficiently smooth in the interior of the triangles
$K\subset\mathcal{T}$ by%
\[
\left.  \left(  \nabla_{\mathcal{T}}u\right)  \right\vert _{\overset{\circ}%
{K}}:=\nabla\left(  \left.  u\right\vert _{\overset{\circ}{K}}\right)
\quad\text{and\quad}\left.  \left(  \operatorname*{div}\nolimits_{\mathcal{T}%
}\mathbf{w}\right)  \right\vert _{\overset{\circ}{K}}:=\operatorname*{div}%
\left(  \left.  \mathbf{w}\right\vert _{\overset{\circ}{K}}\right)  .
\]
The values on $\partial K$ are arbitrary since $\partial K$ has measure zero.

An important measure for the quality of a finite element triangulation is the
shape-regularity constant given by%
\begin{equation}
\gamma_{\mathcal{T}}:=\max_{K\in\mathcal{T}}\frac{h_{K}}{\rho_{K}}
\label{defgammat}%
\end{equation}
with the local mesh width $h_{K}:=\operatorname*{diam}K$ and $\rho_{K}$
denoting the diameter of the largest inscribed ball in $K$. The global mesh
width is $h_{\mathcal{T}}:=\max\left\{  h_{K}:K\in\mathcal{T}\right\}  $.

\begin{remark}
\label{Remangle}It is well known that the shape-regularity implies that there
exists some minimal angle $\phi_{\mathcal{T}}>0$ depending only on
$\gamma_{\mathcal{T}}$ such that every triangle angle in $\mathcal{T}$ is
bounded from below by $\phi_{\mathcal{T}}$. In turn, every triangle angle in
$\mathcal{T}$ is bounded from above by $\pi-2\phi_{\mathcal{T}}$.
\end{remark}

The set of edges in $\mathcal{T}$ is denoted by $\mathcal{E}\left(
\mathcal{T}\right)  $, while the subset of boundary edges is $\mathcal{E}%
_{\partial\Omega}\left(  \mathcal{T}\right)  :=\left\{  E\in\mathcal{E}\left(
\mathcal{T}\right)  :E\subset\partial\left(  \operatorname*{dom}%
\mathcal{T}\right)  \right\}  $; the subset of inner edges is given by
$\mathcal{E}_{\Omega}\left(  \mathcal{T}\right)  :=\mathcal{E}\left(
\mathcal{T}\right)  \backslash\mathcal{E}_{\partial\Omega}\left(
\mathcal{T}\right)  $. For each edge $E\in\mathcal{E}$ we fix a unit vector
$\mathbf{n}_{E}$ orthogonal to $E$ with the convention that $\mathbf{n}_{E}$
is the outer normal vector for boundary edges $E\in\mathcal{E}_{\partial
\Omega}$.

The set of triangle vertices in $\mathcal{T}$ is denoted by $\mathcal{V}%
\left(  \mathcal{T}\right)  $, while the subset of inner vertices is
$\mathcal{V}_{\Omega}\left(  \mathcal{T}\right)  :=\left\{  \mathbf{V}%
\in\mathcal{V}\left(  \mathcal{T}\right)  :\mathbf{V}\notin\partial\left(
\operatorname*{dom}\mathcal{T}\right)  \right\}  $ and $\mathcal{V}%
_{\partial\Omega}\left(  \mathcal{T}\right)  :=\mathcal{V}\left(
\mathcal{T}\right)  \backslash\mathcal{V}_{\Omega}\left(  \mathcal{T}\right)
$. For $K\in\mathcal{T}$, the set of its vertices is denoted by $\mathcal{V}%
\left(  K\right)  $. For $E\in\mathcal{E}\left(  \mathcal{T}\right)  $, we
define the edge patch by%
\[
\mathcal{T}_{E}:=\left\{  K\in\mathcal{T}:E\subset K\right\}  \quad
\text{and\quad}\omega_{E}:=%
%TCIMACRO{\dbigcup \limits_{K\in\mathcal{T}_{E}}}%
%BeginExpansion
{\displaystyle\bigcup\limits_{K\in\mathcal{T}_{E}}}
%EndExpansion
K.
\]
For $\mathbf{z}\in\mathcal{V}\left(  \mathcal{T}\right)  $, the nodal patch is
defined by%
\begin{equation}
\mathcal{T}_{\mathbf{z}}:=\left\{  K\in\mathcal{T}:\mathbf{z}\in K\right\}
\quad\text{and\quad}\omega_{\mathbf{z}}:=%
%TCIMACRO{\dbigcup \limits_{K\in\mathcal{T}_{\mathbf{z}}}}%
%BeginExpansion
{\displaystyle\bigcup\limits_{K\in\mathcal{T}_{\mathbf{z}}}}
%EndExpansion
K \label{nodalpatch}%
\end{equation}
with local mesh width $h_{\mathbf{z}}:=\max\left\{  h_{K}:K\in\mathcal{T}%
_{\mathbf{z}}\right\}  $. For $K\in\mathcal{T}$, we set
\begin{equation}
\mathcal{T}_{K}:=\left\{  K^{\prime}\in\mathcal{T}\mid K\cap K^{\prime}%
\neq\emptyset\right\}  \quad\text{and\quad}\omega_{K}:=%
%TCIMACRO{\dbigcup \limits_{K^{\prime}\in\mathcal{T}_{K}}}%
%BeginExpansion
{\displaystyle\bigcup\limits_{K^{\prime}\in\mathcal{T}_{K}}}
%EndExpansion
K^{\prime}. \label{trianglepatch}%
\end{equation}

For a subset $M\subset\mathbb{R}^{2}$, we denote by $\left[  M\right]  $ its
convex hull; in this way an edge $E$ with endpoints $\mathbf{a},\mathbf{b}$
can be written as $E=\left[  \mathbf{a},\mathbf{b}\right]  =\left[
\mathbf{b},\mathbf{a}\right]  $.

Let $\mathbb{N}=\left\{  1,2,\ldots\right\}  $ and $\mathbb{N}_{0}%
:=\mathbb{N\cup}\left\{  0\right\}  $. For $m\in\mathbb{N}$, we employ the
usual multiindex notation for $%
%TCIMACRO{\TeXButton{boldmu}{\mbox{\boldmath$ \mu$}}}%
%BeginExpansion
\mbox{\boldmath$ \mu$}%
%EndExpansion
=\left(  \mu_{i}\right)  _{i=1}^{m}\in\mathbb{N}_{0}^{m}$ and points
$\mathbf{x}=\left(  x_{i}\right)  _{i=1}^{m}\in\mathbb{R}^{m}$%
\[
\left\vert
%TCIMACRO{\TeXButton{boldmu}{\mbox{\boldmath$ \mu$}}}%
%BeginExpansion
\mbox{\boldmath$ \mu$}%
%EndExpansion
\right\vert :=\mu_{1}+\ldots+\mu_{m},\quad\mathbf{x}^{%
%TCIMACRO{\TeXButton{boldmu}{\mbox{\boldmath$ \mu$}}}%
%BeginExpansion
\mbox{\boldmath$ \mu$}%
%EndExpansion
}:=%
%TCIMACRO{\dprod \limits_{j=1}^{m}}%
%BeginExpansion
{\displaystyle\prod\limits_{j=1}^{m}}
%EndExpansion
x_{j}^{\mu_{j}}.
\]

Let $\mathbb{P}_{m,k}$ denote the space of $m$-variate polynomials of maximal
degree $k$, consisting of functions of the form%
\[
\sum_{\substack{%
%TCIMACRO{\TeXButton{boldmu}{\mbox{\boldmath$ \mu$}}}%
%BeginExpansion
\mbox{\boldmath$ \mu$}%
%EndExpansion
\in\mathbb{N}_{0}^{m}\\\left\vert
%TCIMACRO{\TeXButton{boldmu}{\mbox{\boldmath$ \mu$}}}%
%BeginExpansion
\mbox{\boldmath$ \mu$}%
%EndExpansion
\right\vert \leq k}}a_{%
%TCIMACRO{\TeXButton{boldmu}{\mbox{\boldmath$ \mu$}}}%
%BeginExpansion
\mbox{\boldmath$ \mu$}%
%EndExpansion
}\mathbf{x}^{%
%TCIMACRO{\TeXButton{boldmu}{\mbox{\boldmath$ \mu$}}}%
%BeginExpansion
\mbox{\boldmath$ \mu$}%
%EndExpansion
}%
\]
for real coefficients $a_{%
%TCIMACRO{\TeXButton{boldmu}{\mbox{\boldmath$ \mu$}}}%
%BeginExpansion
\mbox{\boldmath$ \mu$}%
%EndExpansion
}$. Formally, we set $\mathbb{P}_{m,-1}:=\left\{  0\right\}  $. To indicate
the domain explicitly in notation we write sometimes $\mathbb{P}_{k}\left(
D\right)  $ for $D\subset\mathbb{R}^{m}$ and skip the index $m$ since it is
then clear from the argument $D$.

We introduce the following finite element spaces%
\begin{equation}%
\begin{array}
[c]{ll}
& \mathbb{P}_{k}\left(  \mathcal{T}\right)  :=\left\{  q\in L^{2}\left(
\Omega\right)  \mid\forall K\in\mathcal{T}:\left.  q\right\vert _{\overset
{\circ}{K}}\in\mathbb{P}_{k}\left(  \overset{\circ}{K}\right)  \right\}  ,\\
\text{and (cf. (\ref{defdomT}))} & \mathbb{P}_{k,0}\left(  \mathcal{T}\right)
:=\left\{  q\in\mathbb{P}_{k}\left(  \mathcal{T}\right)  :\int
_{\operatorname*{dom}\mathcal{T}}q=0\right\}  .
\end{array}
\label{Pkdefs}%
\end{equation}
Furthermore, let%
\[%
\begin{array}
[c]{ll}
& S_{k}\left(  \mathcal{T}\right)  :=\mathbb{P}_{k}\left(  \mathcal{T}\right)
\cap H^{1}\left(  \operatorname*{dom}\mathcal{T}\right)  ,\\
\text{and} & S_{k,0}\left(  \mathcal{T}\right)  :=S_{k}\left(  \mathcal{T}%
\right)  \cap H_{0}^{1}\left(  \operatorname*{dom}\mathcal{T}\right)  .
\end{array}
\]
The vector-valued versions are denoted by $\mathbf{S}_{k}\left(
\mathcal{T}\right)  :=S_{k}\left(  \mathcal{T}\right)  ^{2}$ and
$\mathbf{S}_{k,0}\left(  \mathcal{T}\right)  :=S_{k,0}\left(  \mathcal{T}%
\right)  ^{2}$. Finally, we define the Crouzeix-Raviart space by%
%TCIMACRO{\TeXButton{CRdeffull}{\begin{subequations}
%\label{CRdeffull}
%\end{subequations}}}%
%BeginExpansion
\begin{subequations}
\label{CRdeffull}
\end{subequations}%
%EndExpansion%
\begin{align}
\operatorname*{CR}\nolimits_{k}\left(  \mathcal{T}\right)   &  :=\left\{
v\in\mathbb{P}_{k}\left(  \mathcal{T}\right)  \mid\forall q\in\mathbb{P}%
_{k-1}\left(  E\right)  \quad\forall E\in\mathcal{E}_{\Omega}\left(
\mathcal{T}\right)  \quad\int_{E}\left[  v\right]  _{E}q=0\right\}  ,\tag{%
%TCIMACRO{\TeXButton{CRdeffull}{\ref{CRdeffull}}}%
%BeginExpansion
\ref{CRdeffull}%
%EndExpansion
a}\label{CRdeffulla}\\
\operatorname*{CR}\nolimits_{k,0}\left(  \mathcal{T}\right)   &  :=\left\{
v\in\operatorname*{CR}\nolimits_{k}\left(  \mathcal{T}\right)  \mid\forall
q\in\mathbb{P}_{k-1}\left(  E\right)  \quad\forall E\in\mathcal{E}%
_{\partial\Omega}\left(  \mathcal{T}\right)  \quad\int_{E}vq=0\right\}  .
\tag{%
%TCIMACRO{\TeXButton{CRdeffull}{\ref{CRdeffull}}}%
%BeginExpansion
\ref{CRdeffull}%
%EndExpansion
b}\label{PCR0RB}%
\end{align}
Here, $\left[  v\right]  _{E}$ denotes the jump of $v\in\mathbb{P}_{k}\left(
\mathcal{T}\right)  $ across an edge $E\in\mathcal{E}_{\Omega}\left(
\mathcal{T}\right)  $
\[
\left[  u\right]  _{E}\left(  \mathbf{x}\right)  :=\lim_{\varepsilon\searrow
0}\left(  u\left(  \mathbf{x}+\varepsilon\mathbf{n}_{E}\right)  -u\left(
\mathbf{x}-\varepsilon\mathbf{n}_{E}\right)  \right)  .
\]
and $\mathbb{P}_{k-1}\left(  E\right)  $ is the space of polynomials of
maximal degree $k-1$ with respect to the local variable in $E$.

We have collected all ingredients for defining the Crouzeix-Raviart
discretization for the Stokes equation. For $k\in\mathbb{N}$, let the discrete
velocity space and pressure space be defined by
\[
\mathbf{CR}_{k,0}\left(  \mathcal{T}\right)  :=\left(  \operatorname*{CR}%
\nolimits_{k,0}\left(  \mathcal{T}\right)  \right)  ^{2}\quad\text{and\quad
}M_{k-1}\left(  \mathcal{T}\right)  :=\mathbb{P}_{k-1,0}\left(  \mathcal{T}%
\right)  .
\]
Then, the discretization is given by: find $\left(  \mathbf{u}%
_{\operatorname*{CR}},p_{\operatorname*{disc}}\right)  \in\mathbf{CR}%
_{k,0}\left(  \mathcal{T}\right)  \times M_{k-1}\left(  \mathcal{T}\right)
\;$such that%
\begin{equation}
\left\{
\begin{array}
[c]{lll}%
a_{\mathcal{T}}\left(  \mathbf{u}_{\operatorname*{CR}},\mathbf{v}\right)
-b_{\mathcal{T}}\left(  \mathbf{v},p_{\operatorname*{disc}}\right)  & =\left(
\mathbf{f},\mathbf{v}\right)  _{\mathbf{L}^{2}\left(  \Omega\right)  } &
\forall\mathbf{v}\in\mathbf{CR}_{k,0}\left(  \mathcal{T}\right)  ,\\
b_{\mathcal{T}}\left(  \mathbf{u}_{\operatorname*{CR}},q\right)  & =0 &
\forall q\in M_{k-1}\left(  \mathcal{T}\right)  ,
\end{array}
\right.  \label{discrStokes}%
\end{equation}
where the bilinear forms $a_{\mathcal{T}}:\mathbf{CR}_{k,0}\left(
\mathcal{T}\right)  \times\mathbf{CR}_{k,0}\left(  \mathcal{T}\right)
\rightarrow\mathbb{R}$ and $b_{\mathcal{T}}:\mathbf{CR}_{k,0}\left(
\mathcal{T}\right)  \times M_{k-1}\left(  \mathcal{T}\right)  \rightarrow
\mathbb{R}$ are given by%
\[
a_{\mathcal{T}}\left(  \mathbf{u},\mathbf{v}\right)  :=\left(  \nabla
_{\mathcal{T}}\mathbf{u},\nabla_{\mathcal{T}}\mathbf{v}\right)  _{\mathbb{L}%
^{2}\left(  \Omega\right)  }\quad\text{and\quad}b_{\mathcal{T}}\left(
\mathbf{v},q\right)  :=\left(  \operatorname*{div}\nolimits_{\mathcal{T}%
}\mathbf{v},q\right)  _{L^{2}\left(  \Omega\right)  }.
\]

It is well known that problem (\ref{discrStokes}) is well-posed if (i): the
bilinear form $a_{\mathcal{T}}\left(  \cdot,\cdot\right)  $ is coercive and
(ii): $b_{\mathcal{T}}\left(  \cdot,\cdot\right)  $ satisfies the inf-sup condition.

To verify the condition (i) we introduce, for a conforming triangulation
$\mathcal{T}$ of the domain $\Omega$, the \textit{broken Sobolev space}%
\[
H^{1}\left(  \mathcal{T}\right)  :=\left\{  u\in L^{2}\left(  \Omega\right)
\mid\forall K\in\mathcal{T}:\left.  u\right\vert _{K}\in H^{1}\left(
K\right)  \right\}
\]
and define, for $u\in H^{1}\left(  \mathcal{T}\right)  $, the \textit{broken
}$H^{1}$\textit{--seminorm} by%
\[
\left\Vert u\right\Vert _{H^{1}\left(  \mathcal{T}\right)  }:=\left\Vert
\nabla_{\mathcal{T}}u\right\Vert _{\mathbf{L}^{2}\left(  \Omega\right)
}=\left(  \sum_{K\in\mathcal{T}}\left\Vert \nabla u\right\Vert _{\mathbf{L}%
^{2}\left(  K\right)  }^{2}\right)  ^{1/2}.
\]

In \cite[Lem. 2]{CrouzeixRaviart}) it is proved that $\left\Vert
\cdot\right\Vert _{\mathbf{H}^{1}\left(  \mathcal{T}\right)  }$ defines a norm
in $\mathbf{CR}_{k,0}\left(  \mathcal{T}\right)  +\mathbf{H}_{0}^{1}\left(
\Omega\right)  $ which is equivalent to the norm $\left(  \sum_{K\in
\mathcal{T}}\left\Vert \mathbf{u}\right\Vert _{\mathbf{H}^{1}\left(  K\right)
}^{2}\right)  ^{1/2}$ with equivalence constants independent of the polynomial
degree and the mesh width (see Theorem \ref{TheoDiscFried}). This directly
implies the coercivity of $a_{\mathcal{T}}\left(  \cdot,\cdot\right)  $:%
\[
a_{\mathcal{T}}\left(  \mathbf{u},\mathbf{u}\right)  \geq\left\Vert
\mathbf{u}\right\Vert _{\mathbf{H}^{1}\left(  \mathcal{T}\right)  }^{2}%
\quad\forall\mathbf{u}\in\mathbf{CR}_{k,0}\left(  \mathcal{T}\right)  .
\]
Hence, well-posedness of (\ref{discrStokes}) follows from the inf-sup
condition for $b_{\mathcal{T}}\left(  \cdot,\cdot\right)  $.

\begin{definition}
Let $\mathcal{T}$ denote a conforming triangulation for $\Omega$. The pair
$\mathbf{CR}_{k,0}\left(  \mathcal{T}\right)  \times M_{k-1}\left(
\mathcal{T}\right)  $ is \emph{inf-sup stable} if there exists a constant
$c_{\mathcal{T},k}$ such that%
\begin{equation}
\inf_{p\in M_{k-1}\left(  \mathcal{T}\right)  \backslash\left\{  0\right\}
}\sup_{\mathbf{v}\in\mathbf{CR}_{k,0}\left(  \mathcal{T}\right)
\backslash\left\{  \mathbf{0}\right\}  }\frac{\left(  p,\operatorname*{div}%
_{\mathcal{T}}\mathbf{v}\right)  _{L^{2}\left(  \Omega\right)  }}{\left\Vert
\mathbf{v}\right\Vert _{\mathbf{H}^{1}\left(  \mathcal{T}\right)  }\left\Vert
p\right\Vert _{L^{2}\left(  \Omega\right)  }}\geq c_{\mathcal{T},k}>0.
\label{infsupcond}%
\end{equation}

\end{definition}

We are now in the position to formulate our main theorem.

\begin{theorem}
\label{Theomain}Let $\Omega\subset\mathbb{R}^{2}$ be a bounded polygonal
Lipschitz domain and let $\mathcal{T}$ denote a conforming triangulation of
$\Omega$ consisting of more than a single triangle. Let $k\in\mathbb{N}$. If
$k\geq3$ is odd we assume that $\mathcal{T}$ contains at least one inner
vertex. Then, the inf-sup condition (\ref{infsupcond}) holds:%
\begin{equation}
c_{\mathcal{T},k}\geq c_{\mathcal{T}}\left(  \log\left(  k+1\right)  \right)
^{-\alpha} \label{infsupp_dep}%
\end{equation}
for a constant $c_{\mathcal{T}}>0$ depending only on the shape-regularity of
the mesh and on the maximal outer angle $\alpha_{\Omega}$. In particular
$c_{\mathcal{T}}$ is independent of the mesh width $h_{\mathcal{T}}$ and the
polynomial degree $k$. The value of $\alpha\geq0$ is given by%
\[
\alpha=\left\{
\begin{array}
[c]{ll}%
1/2 & \left\{
\begin{array}
[c]{l}%
\text{if }k\text{ is even,}\\
\text{or }k\geq3\text{ is odd and all triangles in }\mathcal{T}\text{ have at
least one inner vertex,}%
\end{array}
\right. \\
\left(  1+L\right)  /2 & \text{otherwise,}%
\end{array}
\right.
\]
where $L$ depends only on the mesh topology via the number of steps involved
in the step-by-step construction introduced in (\ref{stepbystepextension}).
\end{theorem}

%

%TCIMACRO{\TeXButton{Proof}{\proof}}%
%BeginExpansion
\proof
%EndExpansion
The estimate $c_{\mathcal{T},k}>0$ follows, for $k=1$ from
\cite{CrouzeixRaviart}, for $k=2$ from \cite[Thm. 3.1]{ChaLeeLee}, for even
$k\geq4$ from \cite{Baran_Stoyan}, for odd $k\geq5$ from \cite{CCSS_CR_1}, and
for $k=3$ from \cite{CCSS_CR_2}. We set%
\[
c_{\mathcal{T},\operatorname*{low}}:=\min\left\{  c_{\mathcal{T},k}:1\leq
k\leq3\right\}  .
\]
Estimate (\ref{infsupp_dep}) for some $c_{\mathcal{T}}:=c_{\mathcal{T}%
,\operatorname*{high}}^{\operatorname*{odd}}>0$ for odd $k\geq5$ is proved in
Section \ref{PTpge5}, Lem. \ref{LemmaFinOdd}, while the estimate for some
$c_{\mathcal{T}}:=c_{\mathcal{T},\operatorname*{high}}^{\operatorname*{even}%
}>0$ for even $k\geq4$ is proved in Section \ref{SecEvenk}. Both constants
$c_{\mathcal{T},\operatorname*{high}}^{\operatorname*{odd}}$, $c_{\mathcal{T}%
,\operatorname*{high}}^{\operatorname*{even}}$ depend only on the
shape-regularity of the mesh and $\alpha_{\Omega}$. Hence, $c_{\mathcal{T}%
}\geq\min\left\{  c_{\mathcal{T},\operatorname*{low}},c_{\mathcal{T}%
,\operatorname*{high}}^{\operatorname*{odd}},c_{\mathcal{T}%
,\operatorname*{high}}^{\operatorname*{even}}\right\}  $.%
%TCIMACRO{\TeXButton{End Proof}{\endproof}}%
%BeginExpansion
\endproof
%EndExpansion

We emphasize that the original definition in \cite{CrouzeixRaviart} allows for
slightly more general finite element spaces, more precisely, the spaces
$\operatorname*{CR}_{k}\left(  \mathcal{T}\right)  $ can be enriched by
locally supported functions. From this point of view, the definition
(\ref{CRdeffull}) describes a minimal Crouzeix-Raviart space.

The possibility for enrichment has been used frequently in the literature to
prove inf-sup stability for the arising finite element spaces (see, e.g.,
\cite{CrouzeixRaviart}, \cite{Guzman_divfree}, \cite{Matthies_nonconf_2005}).
In contrast, we will prove the $k$-explicit estimate of the inf-sup constant
for the Crouzeix-Raviart space $\operatorname*{CR}_{k}\left(  \mathcal{T}%
\right)  $.

\section{Proof of Theorem \ref{Theomain}\label{SecProof}}

In this section, we will analyse the $k$-dependence of the inf-sup constant in
the form (\ref{infsupp_dep}), first for odd polynomial degree $k\geq5$ and
then for even degree $k\geq4$.

\subsection{Barycentric coordinates and basis functions for the
velocity\label{VelBasis}}

In this section, we introduce basis functions for the finite element spaces in
Section \ref{NuDiscrete}. We begin with introducing some general notation.

\begin{notation}
\label{Notation}For vectors $\mathbf{a}_{i}\in\mathbb{R}^{n}$, $1\leq i\leq
m$, we write $\left[  \mathbf{a}_{1}\mid\mathbf{a}_{2}\mid\ldots\mid
\mathbf{a}_{m}\right]  $ for the $n\times m$ matrix with column vectors
$\mathbf{a}_{i}$. For $\mathbf{v}=\left(  v_{1},v_{2}\right)  ^{T}%
\in\mathbb{R}^{2}$ we set $\mathbf{v}^{\perp}:=\left(  v_{2},-v_{1}\right)
^{T}$. Let\ $\mathbf{e}_{k,i}\in\mathbb{R}^{k}$\ be the $i$-th canonical unit
vector in $\mathbb{R}^{k}$.

For $\mathbf{v}\in\mathbb{R}^{n}$, $\left\Vert \mathbf{v}\right\Vert $ is the
Euclidean vector norm while the induced matrix norm is given for
$\mathbf{B}\in\mathbb{R}^{n\times n}$ by $\left\Vert \mathbf{B}\right\Vert
:=\sup\left\{  \left\Vert \mathbf{Bx}\right\Vert /\left\Vert \mathbf{x}%
\right\Vert :\mathbf{x}\in\mathbb{R}^{n}\backslash\left\{  \mathbf{0}\right\}
\right\}  $.

Vertices in a triangle are numbered counterclockwise. In a triangle $K$ with
vertices $\mathbf{A}_{1}$, $\mathbf{A}_{2}$, $\mathbf{A}_{3}$ the angle at
$\mathbf{A}_{i}$ is called $\alpha_{i}$. If a triangle is numbered by an index
(e.g., $K_{\ell}$), the angle at $A_{\ell,i}$ is called $\alpha_{\ell,i}$. For
quantities in a triangle $K$ as, e.g., angles $\alpha_{j}$, $1\leq j\leq3$, we
use the cyclic numbering convention $\alpha_{3+1}:=\alpha_{1}$ and
$\alpha_{1-1}:=\alpha_{3}$.

For a $d$-dimensional measurable set $D,$ we write $\left\vert D\right\vert $
for its measure; for a discrete set, say $\mathcal{J}$, we denote by
$\left\vert \mathcal{J}\right\vert $ its cardinality.

In the proofs, we consider frequently nodal patches $\mathcal{T}_{\mathbf{z}}$
for inner vertices $\mathbf{z}\in\mathcal{V}_{\Omega}\left(  \mathcal{T}%
\right)  $. The number $m$ denotes the number of triangles in $\mathcal{T}%
_{\mathbf{z}}$. Various quantities in this patch such as, e.g., the triangles
in $\mathcal{T}_{\mathbf{z}}$, have an index which runs from $1$ to $m$. Here,
we use the cyclic numbering convention $K_{m+1}:=K_{1}$ and $K_{1-1}:=K_{m}$
and apply this analogously for other quantities in the nodal patch.
\end{notation}

Let the closed reference triangle $\widehat{K}$ be the triangle with vertices
$\mathbf{\hat{A}}_{1}:=\left(  0,0\right)  ^{T}$, $\mathbf{\hat{A}}%
_{2}:=\left(  1,0\right)  ^{T}$, $\mathbf{\hat{A}}_{3}:=\left(  0,1\right)
^{T}$. The nodal points on the reference element of order $k\in\mathbb{N}_{0}$
are given by%
\[
\widehat{\mathcal{N}}_{k}:=\left\{
\begin{array}
[c]{ll}%
\left\{  \dfrac{1}{k}%
%TCIMACRO{\TeXButton{boldmu}{\mbox{\boldmath$ \mu$}}}%
%BeginExpansion
\mbox{\boldmath$ \mu$}%
%EndExpansion
\mid%
%TCIMACRO{\TeXButton{boldmu}{\mbox{\boldmath$ \mu$}}}%
%BeginExpansion
\mbox{\boldmath$ \mu$}%
%EndExpansion
\in\mathbb{N}_{0}^{2}:\mathbb{\quad}\left\vert
%TCIMACRO{\TeXButton{boldmu}{\mbox{\boldmath$ \mu$}}}%
%BeginExpansion
\mbox{\boldmath$ \mu$}%
%EndExpansion
\right\vert \leq k\right\}  & k\geq1,\\
& \\
\left\{  \left(  \dfrac{1}{3},\dfrac{1}{3}\right)  \right\}  & k=0.
\end{array}
\right.
\]
For a triangle $K\subset\mathbb{R}^{2}$, we denote by $\chi_{K}:\widehat
{K}\rightarrow K$ an affine bijection. The mapped nodal points of order
$k\in\mathbb{N}_{0}$ on $K$ are given by%
\[
\mathcal{N}_{k}\left(  K\right)  :=\left\{  \chi_{K}\left(  \mathbf{z}\right)
:\mathbf{z}\in\widehat{\mathcal{N}}_{k}\right\}  .
\]
Nodal points of order $k$ on $\mathcal{T}$ are defined by%
\[
\mathcal{N}_{k}\left(  \mathcal{T}\right)  :=%
%TCIMACRO{\dbigcup \limits_{K\in\mathcal{T}}}%
%BeginExpansion
{\displaystyle\bigcup\limits_{K\in\mathcal{T}}}
%EndExpansion
\mathcal{N}_{k}\left(  K\right)  \text{,\quad}\mathcal{N}_{\partial\Omega}%
^{k}\left(  \mathcal{T}\right)  :=\mathcal{N}_{k}\left(  \mathcal{T}\right)
\cap\partial\Omega,\quad\text{and\quad}\mathcal{N}_{k,\Omega}\left(
\mathcal{T}\right)  :=\mathcal{N}_{k}\left(  \mathcal{T}\right)  \cap\Omega.
\]
We introduce the Lagrange basis for the space $S_{k}\left(  \mathcal{T}%
\right)  $, which is indexed by the nodal points $\mathbf{z}\in\mathcal{N}%
_{k}\left(  \mathcal{T}\right)  $ and characterized by
\begin{equation}
B_{k,\mathbf{z}}\in S_{k}\left(  \mathcal{T}\right)  \quad\text{and\quad
}\forall\mathbf{z}^{\prime}\in\mathcal{N}_{k}\left(  \mathcal{T}\right)
\qquad B_{k,\mathbf{z}}\left(  \mathbf{z}^{\prime}\right)  =\delta
_{\mathbf{z},\mathbf{z}^{\prime}}, \label{basisfunctions}%
\end{equation}
where $\delta_{\mathbf{z},\mathbf{z}^{\prime}}$ is the Kronecker delta. A
basis for the space $S_{k,0}\left(  \mathcal{T}\right)  $ is given by
$B_{k,\mathbf{z}}$, $\mathbf{z}\in\mathcal{N}_{k,\Omega}\left(  \mathcal{T}%
\right)  $.\bigskip

Let $K$ denote a triangle with vertices $\mathbf{A}_{i}$, $1\leq i\leq3$, and
let $\lambda_{K,\mathbf{A}_{i}}\in\mathbb{P}_{1}\left(  K\right)  $ be the
\textit{barycentric coordinate} for the node $\mathbf{A}_{i}$ defined by%
\begin{equation}
\lambda_{K,\mathbf{A}_{i}}\left(  \mathbf{A}_{j}\right)  =\delta_{i,j}%
\quad1\leq i,j\leq3. \label{lambdaintro1}%
\end{equation}
If the numbering of the vertices in $K$ is fixed, we write $\lambda_{K,i}$
short for $\lambda_{K,\mathbf{A}_{i}}$. For the barycentric coordinate on the
reference element $\widehat{K}$ for the vertex $\mathbf{\hat{A}}_{j}$ we write
$\widehat{\lambda}_{j}$, $j=1,2,3$. Elementary calculation yield (see, e.g.,
\cite[Appendix A]{CCSS_CR_1})%
\begin{equation}
\partial_{\mathbf{n}_{k}}\lambda_{K,\mathbf{A}_{i}}=\frac{\left\vert
E_{i}\right\vert }{2\left\vert K\right\vert }\times\left\{
\begin{array}
[c]{ll}%
-1 & i=k,\\
\cos\alpha_{\ell} & \ell\text{ s.t. }\left\{  \ell,i,k\right\}  =\left\{
1,2,3\right\}  ,
\end{array}
\right.  \label{normcomp}%
\end{equation}
where $E_{i}$ is the edge of $K$ opposite to $\mathbf{A}_{i}$, $\mathbf{n}%
_{k}$ the outward unit normal at $E_{k}$, and $\alpha_{\ell}$ the angle in $K$
at $\mathbf{A}_{\ell}$.

\begin{definition}
Let $L_{k}$ denote the usual univariate Legendre polynomial of degree $k$ (see
\cite[Table 18.3.1]{NIST:DLMF}). Let $k\in\mathbb{N}$ be even and
$K\in\mathcal{T}$. Then, the \emph{non-conforming triangle bubble }is given by%
\[
B_{k,K}^{\operatorname*{CR}}:=\left\{
\begin{array}
[c]{ll}%
\dfrac{1}{2}\left(  -1+%
%TCIMACRO{\dsum \limits_{i=1}^{3}}%
%BeginExpansion
{\displaystyle\sum\limits_{i=1}^{3}}
%EndExpansion
L_{k}\left(  1-2\lambda_{K,i}\right)  \right)  & \text{on }K,\\
0 & \text{on }\Omega\backslash K.
\end{array}
\right.
\]
For $k$ odd and $E\in\mathcal{E}\left(  \mathcal{T}\right)  $, the
\emph{non-conforming edge bubble} is given by%
\begin{equation}
B_{k,E}^{\operatorname*{CR}}:=\left\{
\begin{array}
[c]{ll}%
L_{k}\left(  1-2\lambda_{K,\mathbf{A}_{K,E}}\right)  & \text{on }K\text{ for
}K\in\mathcal{T}_{E},\\
0 & \text{on }\Omega\backslash\omega_{E},
\end{array}
\right.  \label{BkEdef}%
\end{equation}
where $\mathbf{A}_{K,E}$ denotes the vertex in $K$ opposite to $E$.
\end{definition}

Different representations of the functions $B_{k,E}^{\operatorname*{CR}}$,
$B_{k,K}^{\operatorname*{CR}}$ exist in the literature, see \cite{BaranCVD},
\cite{Ainsworth_Rankin}, \cite[for $p=4,6.$]{ChaLeeLee}, \cite{ccss_2012}
while the formula for $B_{k,K}^{\operatorname*{CR}}$ has been introduced in
\cite{Baran_Stoyan} and the one for $B_{k,E}^{\operatorname*{CR}}$ in
\cite{CCSS_CR_1}.

\begin{proposition}
\label{ThmBasisCRscalar}A basis for the space $\operatorname*{CR}_{k,0}\left(
\mathcal{T}\right)  $ is given

\begin{enumerate}
\item for even $k$ by%
\[
\left\{  B_{k,\mathbf{z}}\mid\mathbf{z}\in\mathcal{N}_{k,\Omega}\left(
\mathcal{T}\right)  \right\}  \cup\left\{  B_{k,K}^{\operatorname*{CR}}\mid
K\in\mathcal{T}\right\}  ,
\]

\item for odd $k$ by%
\[
\left\{  B_{k,\mathbf{z}}\mid\mathbf{z}\in\mathcal{N}_{k,\Omega}\left(
\mathcal{T}\right)  \backslash\mathcal{V}_{\Omega}\left(  \mathcal{T}\right)
\right\}  \cup\left\{  B_{k,E}^{\operatorname*{CR}}\mid E\in\mathcal{E}%
_{\Omega}\left(  \mathcal{T}\right)  \right\}  .
\]

\end{enumerate}
\end{proposition}

The proof of this proposition and the following corollary can be found, e.g.,
in \cite[Rem. 3]{BaranCVD}, \cite[Thm. 22]{ccss_2012}, \cite[Cor.
3.4]{CCSS_CR_1}.

\begin{corollary}
\label{CorBasis}A basis for the space $\mathbf{CR}_{k,0}\left(  \mathcal{T}%
\right)  $ is given

\begin{enumerate}
\item for even $k$ by%
\begin{equation}%
\begin{array}
[c]{l}%
\left\{  B_{k,\mathbf{z}}\mathbf{v}_{\mathbf{z}}\mid\mathbf{z}\in
\mathcal{N}_{k,\Omega}\left(  \mathcal{T}\right)  \right\}  \cup\left\{
B_{k,\mathbf{z}}\mathbf{w}_{\mathbf{z}}\mid\mathbf{z}\in\mathcal{N}_{k,\Omega
}\left(  \mathcal{T}\right)  \right\} \\
\quad\\
\quad\cup\left\{  B_{k,K}^{\operatorname*{CR}}\mathbf{v}_{K}\mid
K\in\mathcal{T}\right\}  \cup\left\{  B_{k,K}^{\operatorname*{CR}}%
\mathbf{w}_{K}\mid K\in\mathcal{T}\right\}  ,
\end{array}
\label{basisvel}%
\end{equation}

\item for odd $k$ by%
\begin{equation}%
\begin{array}
[c]{l}%
\left\{  B_{k,\mathbf{z}}\mathbf{v}_{\mathbf{z}}\mid\mathbf{z}\in
\mathcal{N}_{k,\Omega}\left(  \mathcal{T}\right)  \backslash\mathcal{V}%
_{\Omega}\left(  \mathcal{T}\right)  \right\}  \cup\left\{  B_{k,\mathbf{z}%
}\mathbf{w}_{\mathbf{z}}\mid\mathbf{z}\in\mathcal{N}_{k,\Omega}\left(
\mathcal{T}\right)  \backslash\mathcal{V}_{\Omega}\left(  \mathcal{T}\right)
\right\} \\
\quad\\
\quad\cup\left\{  B_{k,E}^{\operatorname*{CR}}\mathbf{v}_{E}\mid
E\in\mathcal{E}_{\Omega}\left(  \mathcal{T}\right)  \right\}  \cup\left\{
B_{k,E}^{\operatorname*{CR}}\mathbf{w}_{E}\mid E\in\mathcal{E}_{\Omega}\left(
\mathcal{T}\right)  \right\}  .
\end{array}
\label{basisvelodd}%
\end{equation}

\end{enumerate}

Here, for any nodal point $\mathbf{z}$, the linearly independent vectors
$\mathbf{v}_{\mathbf{z}},\mathbf{w}_{\mathbf{z}}\in\mathbb{R}^{2}$ can be
chosen arbitrarily. The same holds for any triangle $K$ for the vectors
$\mathbf{v}_{K},\mathbf{w}_{K}\in\mathbb{R}^{2}$ in (\ref{basisvel}) and for
any $E\in\mathcal{E}_{\Omega}\left(  \mathcal{T}\right)  $ for the vectors
$\mathbf{v}_{E},\mathbf{w}_{E}\in\mathbb{R}^{2}$ in (\ref{basisvelodd}).
\end{corollary}

\begin{remark}
The original definition of Crouzeix-Raviart spaces by \cite{CrouzeixRaviart}
is implicit and given for conforming simplicial finite element meshes in
$\mathbb{R}^{d}$, $d=2,3$. For their practical implementation, a basis is
needed and Corollary \ref{CorBasis} provides a simple definition. A basis for
Crouzeix-Raviart finite elements in $\mathbb{R}^{3}$ is introduced in
\cite{Fortin_d3} for $k=2$, a general construction is given in \cite{CDS}, and
a basis for a minimal Crouzeix-Raviart spaces in general dimension $d$ is
presented in \cite{SauterTorres_CR3D}.
\end{remark}

\subsection{The case of odd $k\geq5$\label{PTpge5}}

In this section, we assume for the following%
\begin{equation}%
\begin{array}
[c]{ll}%
\text{a)} & k\geq5\text{ is odd and}\\
\text{b)} & \mathcal{T}\text{ is a conforming triangulation and has at least
one inner vertex.}%
\end{array}
\label{Assumptab}%
\end{equation}
This section is structured as follows. In \S \ref{GeoPrel} we generalize the
concept of critical points (see \cite{vogelius1983right}, \cite{ScottVogelius}%
) to $\eta$-critical points which turn out to be essential for estimates with
constants depending on the mesh only via the shape-regularity constant and
$\alpha_{\Omega}$. We split these $\eta$-critical points into a set of
\textquotedblleft obtuse\textquotedblright\ $\eta-$critical points and
\textquotedblleft acute\textquotedblright\ $\eta-$critical points. In
\S \ref{Obtuse}, we provide the proof of Theorem \ref{Theomain} for a maximal
partial triangulation that does not contain acute $\eta-$critical points and
satisfies (\ref{Assumptab}). Finally, in \S \ref{SecAcute} we present the
argument to allow for acute $\eta-$critical points.

\subsubsection{Geometric preliminaries\label{GeoPrel}}

For the analysis of the inf-sup constant we start with the definition of
\textit{critical points} (see \cite{vogelius1983right}, \cite{ScottVogelius}).

\begin{definition}
\label{DefCritpoint}Let $\mathcal{T}$ denote a triangulation as in
\S \ref{NuDiscrete}. For $\mathbf{z}\in\mathcal{V}\left(  \mathcal{T}\right)
$, let
\[
\mathcal{E}_{\mathbf{z}}:=\left\{  E\in\mathcal{E}\left(  \mathcal{T}\right)
:\mathbf{z}\text{ is an endpoint of }E\right\}  .
\]
The point $\mathbf{z}\in\mathcal{V}\left(  \mathcal{T}\right)  $ is a
\emph{critical point} for $\mathcal{T}$ if there exist two straight infinite
lines $L_{1}$, $L_{2}$ in $\mathbb{R}^{2}$ such that all edges $E\in
\mathcal{E}_{\mathbf{z}}$ satisfy $E\subset L_{1}\cup L_{2}$. The set of all
critical points in $\mathcal{T}$ is $\mathcal{C}_{\mathcal{T}}$.
\end{definition}

%

%TCIMACRO{\FRAME{ftbpFU}{5.5737in}{2.9386in}{0pt}{\Qcb{Illustration of the four
%critical cases as in Remark \ref{RemCritGeom}. Left top: inner critical point,
%right top: acute critical point, left bottom: flat critical point, right
%bottom: obtuse critical point.}}{\Qlb{FigCritCases}}{critcases.eps}%
%width 5.5737in;  height 2.9386in;  depth 0pt;  original-width 5.5166in;
%original-height 2.8954in;  cropleft "0";  croptop "1";  cropright "1";
%cropbottom "0";  filename 'critcases.EPS';file-properties "XNPEU";}}}%
%BeginExpansion
\begin{figure}
[ptb]
\begin{center}
\includegraphics[
height=2.9386in,
width=5.5737in
]%
{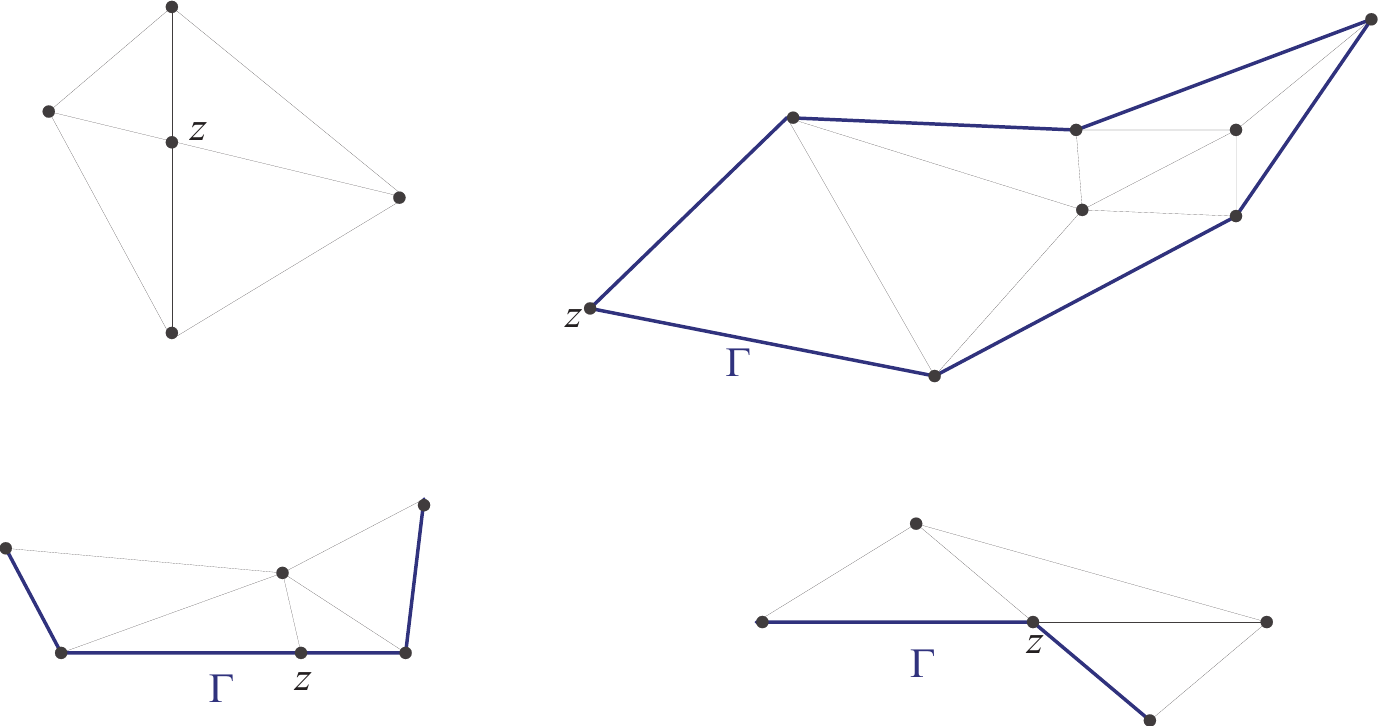}%
\caption{Illustration of the four critical cases as in Remark
\ref{RemCritGeom}. Left top: inner critical point, right top: acute critical
point, left bottom: flat critical point, right bottom: obtuse critical point.}%
\label{FigCritCases}%
\end{center}
\end{figure}
%EndExpansion

\begin{remark}
\label{RemCritGeom}Geometric configurations where critical points occur are
well studied in the literature (see, e.g., \cite{ScottVogelius}). Any critical
point $\mathbf{z}\in\mathcal{C}_{\mathcal{T}}$ belongs to one of the following
cases (see Fig. \ref{FigCritCases}):

\begin{enumerate}
\item $\mathbf{z}\in\mathcal{V}_{\Omega}\left(  \mathcal{T}\right)  $ and
$\mathcal{T}_{\mathbf{z}}$ consists of four triangles and $\mathbf{z}$ is the
intersections of the two diagonals in the quadrilateral $\omega_{\mathbf{z}}$.

\item $\mathbf{z}\in\mathcal{V}_{\partial\Omega}\left(  \mathcal{T}\right)  $
and $\operatorname*{card}\mathcal{E}_{\mathbf{z}}=2$, i.e., both edges
$E\in\mathcal{E}_{\mathbf{z}}$ are boundary edges with joint $\mathbf{z}$.

\item $\mathbf{z}\in\mathcal{V}_{\partial\Omega}\left(  \mathcal{T}\right)  $
and $\operatorname*{card}\mathcal{E}_{\mathbf{z}}=3$ and two edges
$E\in\mathcal{E}_{\mathbf{z}}$ are boundary edges which lie on a straight
boundary piece.

\item $\mathbf{z}\in\mathcal{V}_{\partial\Omega}\left(  \mathcal{T}\right)  $
and $\operatorname*{card}\mathcal{E}_{\mathbf{z}}=4$ and each of the two
boundary edges is aligned with one edge of $\mathcal{E}_{\mathbf{z}}%
\cap\mathcal{E}_{\Omega}\left(  \mathcal{T}\right)  $.
\end{enumerate}
\end{remark}

\begin{definition}
Let $\mathcal{T}$ denote a triangulation as in \S \ref{NuDiscrete}. Let
$\mathbf{z}\in\mathcal{V}\left(  \mathcal{T}\right)  $ and the nodal patch
$\mathcal{T}_{\mathbf{z}}$ as in (\ref{nodalpatch}). Let the triangles
$K_{\ell}$, $1\leq\ell\leq m$, in $\mathcal{T}_{\mathbf{z}}$ be numbered
counterclockwise and denote the angle in $K_{\ell}$ at $\mathbf{z}$ by
$\omega_{\ell}$. Then,
\[
\Theta\left(  \mathbf{z}\right)  :=\left\{
\begin{array}
[c]{ll}%
\max\left\{  \left\vert \sin\left(  \omega_{1}+\omega_{2}\right)  \right\vert
,\left\vert \sin\left(  \omega_{2}+\omega_{3}\right)  \right\vert
,\ldots,\left\vert \sin\left(  \omega_{m}+\omega_{1}\right)  \right\vert
\right\}  & \text{if }\mathbf{z}\in\mathcal{V}_{\Omega}\left(  \mathcal{T}%
\right)  ,\\
\max\left\{  \left\vert \sin\left(  \omega_{1}+\omega_{2}\right)  \right\vert
,\left\vert \sin\left(  \omega_{2}+\omega_{3}\right)  \right\vert
,\ldots,\left\vert \sin\left(  \omega_{m-1}+\omega_{m}\right)  \right\vert
\right\}  & \text{if }\mathbf{z}\in\Gamma\wedge m>1,\\
0 & \text{if }\mathbf{z\in}\Gamma\wedge m=1.
\end{array}
\right.
\]

\end{definition}

\begin{remark}
It is easy to see that $\mathbf{z}\in\mathcal{C}_{\mathcal{T}}$ if and only if
$\Theta\left(  \mathbf{z}\right)  =0$.
\end{remark}

\begin{lemma}
\label{Lemangle}Let $\phi_{\mathcal{T}}$ be as in Remark \ref{Remangle}. Set%
\[
\eta_{0}:=\min\left\{  \frac{1}{2},c_{1},\frac{3\phi_{\mathcal{T}}}{\pi}%
,\sin\phi_{\mathcal{T}}\right\}
\]
with%
\[
c_{1}:=\left\{
\begin{array}
[c]{ll}%
\min\left\{  \sin2\phi_{\mathcal{T}},\left\vert \sin\left(  2\pi
-4\phi_{\mathcal{T}}\right)  \right\vert \right\}  & \phi_{\mathcal{T}}\leq
\pi/8,\\
\sin2\phi_{\mathcal{T}} & \pi/8<\phi_{\mathcal{T}}\leq\pi/4,\\
1 & \phi_{\mathcal{T}}>\pi/4.
\end{array}
\right.
\]
Let $0\leq\eta<\eta_{0}$ be fixed. If, for $\mathbf{z}\in\mathcal{V}\left(
\mathcal{T}\right)  $, it holds $\Theta\left(  \mathbf{z}\right)  \leq\eta$,
then, for any edge $E=\left[  \mathbf{z},\mathbf{z}^{\prime}\right]
\in\mathcal{E}_{\Omega}\left(  \mathcal{T}\right)  $ it holds%
\[
\Theta\left(  \mathbf{z}^{\prime}\right)  \geq\eta_{0}.
\]

\end{lemma}

%

%TCIMACRO{\TeXButton{Proof}{\proof}}%
%BeginExpansion
\proof
%EndExpansion
Let $\mathbf{z}\in\mathcal{V}\left(  \mathcal{T}\right)  $ and consider an
edge $E=\left[  \mathbf{z},\mathbf{z}^{\prime}\right]  \in\mathcal{E}_{\Omega
}\left(  \mathcal{T}\right)  $. Then, there are two triangles $K,K^{\prime}%
\in\mathcal{T}$ which are adjacent to $E$. The angle in $K$ resp. $K^{\prime}$
at $\mathbf{z}$ is denoted by $\omega$ resp. $\omega^{\prime}$.

\textbf{1st case. }Let $\omega+\omega^{\prime}\leq\pi/2$ or $\omega
+\omega^{\prime}\geq\frac{3}{2}\pi$. Then, we conclude from Remark
\ref{Remangle} that%
\[
2\phi_{\mathcal{T}}\leq\omega+\omega^{\prime}\leq\frac{\pi}{2}\quad
\text{or\quad}\frac{3}{2}\pi\leq\omega+\omega^{\prime}\leq2\pi-4\phi
_{\mathcal{T}}.
\]
For the left inequality to hold, the minimal angle must satisfy $\phi
_{\mathcal{T}}\leq\pi/4$ while for the right inequality, it must hold
$\phi_{\mathcal{T}}\leq\pi/8$. For $\phi_{\mathcal{T}}>\pi/4$, the \textbf{1st
case} is empty. For $\phi_{\mathcal{T}}\leq\pi/4$ we get%
\[
\Theta\left(  \mathbf{z}\right)  \geq\left\{
\begin{array}
[c]{ll}%
\min\left\{  \sin2\phi_{\mathcal{T}},\left\vert \sin\left(  2\pi
-4\phi_{\mathcal{T}}\right)  \right\vert \right\}  & \phi_{\mathcal{T}}\leq
\pi/8\\
\sin2\phi_{\mathcal{T}} & \pi/8<\phi_{\mathcal{T}}\leq\pi/4
\end{array}
\right\}  \geq c_{1}\geq\eta_{0}.
\]
Since $\eta<\eta_{0}\leq c_{1}$ this case cannot appear.

\textbf{2nd case. }Let $\pi/2<\omega+\omega^{\prime}<3\pi/2$. The condition
$\left\vert \sin\left(  \omega+\omega^{\prime}\right)  \right\vert \leq\eta$
implies that $\omega+\omega^{\prime}=\pi+\delta$ with%
\begin{equation}
\left\vert \delta\right\vert \leq\arcsin\eta\overset{\text{\cite[4.24.1]%
{NIST:DLMF}}}{=}\eta\sum_{\ell=0}^{\infty}\eta^{2\ell}\frac{\left(
2\ell\right)  !}{\left(  \ell!\right)  ^{2}4^{\ell}\left(  2\ell+1\right)
}\overset{\eta_{0}\leq1/2}{\leq}\eta\sum_{\ell=0}^{\infty}2^{-2\ell}%
\frac{\left(  2\ell\right)  !}{\left(  \ell!\right)  ^{2}4^{\ell}\left(
2\ell+1\right)  }=\frac{\pi\eta}{3}. \label{deltaest}%
\end{equation}
Consequently the two angles $\alpha$ in $K$ and $\alpha^{\prime}$ in
$K^{\prime}$ at $\mathbf{z}^{\prime}$ satisfy%
\[
\alpha+\alpha^{\prime}=2\pi-\omega-\omega^{\prime}-\beta-\beta^{\prime}%
=\pi-\delta-\beta-\beta^{\prime}\leq\pi+\frac{\pi\eta}{3}-2\phi_{\mathcal{T}%
}\overset{\pi\eta/3\leq\phi_{\mathcal{T}}}{\leq}\pi-\phi_{\mathcal{T}},
\]
where $\beta$ (resp. $\beta^{\prime}$) denotes the third angle in $K$ (resp.
$K^{\prime}$). Hence, in this case%
\[
\Theta\left(  \mathbf{z}^{\prime}\right)  \geq\left\vert \sin\left(  \pi
-\phi_{\mathcal{T}}\right)  \right\vert =\sin\phi_{\mathcal{T}}\geq\eta_{0}.
\]%
%TCIMACRO{\TeXButton{End Proof}{\endproof}}%
%BeginExpansion
\endproof
%EndExpansion

\begin{definition}
\label{Defetacriticalpoint}Let $\eta_{0}$ be as in Lemma \ref{Lemangle}. For
$0\leq\eta<\eta_{0}$, the set of $\eta$\emph{-critical points }$\mathcal{C}%
_{\mathcal{T}}\left(  \eta\right)  $ is given by%
\[
\mathcal{C}_{\mathcal{T}}\left(  \eta\right)  :=\left\{  \mathbf{z}%
\in\mathcal{V}\left(  \mathcal{T}\right)  \mid\Theta\left(  \mathbf{z}\right)
\leq\eta\right\}  .
\]
A point $\mathbf{z}\in\mathcal{C}_{\mathcal{T}}\left(  \eta\right)
\backslash\mathcal{C}_{\mathcal{T}}\left(  0\right)  $ is called a
\emph{nearly critical point}. An $\eta$-critical point $\mathbf{z}%
\in\mathcal{C}_{\mathcal{T}}\left(  \eta\right)  $ is \emph{isolated} if all
edge $\left[  \mathbf{z,z}^{\prime}\right]  \in\mathcal{E}\left(
\mathcal{T}\right)  $ satisfy: $\mathbf{z}^{\prime}$ is not an $\eta$-critical point.
\end{definition}

\begin{figure}
[ptb]
\begin{center}
\includegraphics[
height=2.9386in,
width=5.5737in
]%
{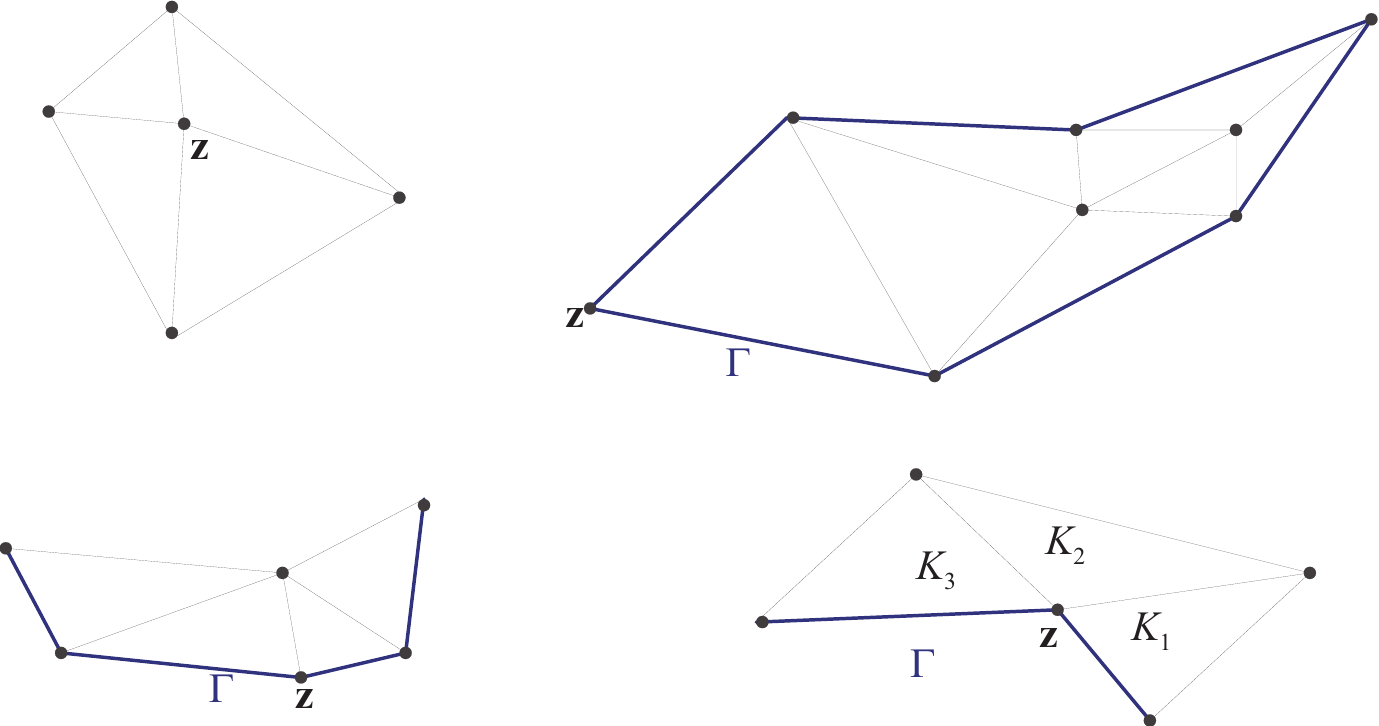}%
\caption{Illustration of the four $\eta$-critical cases as in Remark
\ref{Remetacritcases}. Left top: inner $\eta-$critical point, right top: acute
critical point, left bottom: flat $\eta$-critical point, right bottom: concave
$\eta$-critical point.}%
\label{RemCritGeom2}%
\end{center}
\end{figure}

By perturbing the geometric configurations in Remark \ref{RemCritGeom} we
obtain the following subcases (see Fig. \ref{RemCritGeom2}).

\begin{definition}
\label{Remetacritcases}Let $\eta_{0}$ be as in Lemma \ref{Lemangle} and
$0\leq\eta<\eta_{0}$. If $\mathbf{z}\in\mathcal{V}\left(  \mathcal{T}\right)
$ satisfies

\begin{enumerate}
\item $\mathbf{z}\in\mathcal{V}_{\Omega}\left(  \mathcal{T}\right)  $ and
$\operatorname*{card}\mathcal{T}_{\mathbf{z}}=4$ and $\Theta\left(
\mathbf{z}\right)  \leq\eta$. Then $\mathbf{z}$ is an \emph{inner} $\eta
$\emph{-critical point}. Let%
\[
\mathcal{C}_{\mathcal{T}}^{\operatorname*{inner}}\left(  \eta\right)
:=\left\{  \mathbf{z}\in\mathcal{C}_{\mathcal{T}}\left(  \eta\right)
:\mathbf{z}\text{ is an inner }\eta\text{-critical point}\right\}  .
\]

\item $\mathbf{z}\in\mathcal{V}_{\partial\Omega}\left(  \mathcal{T}\right)  $
and $\operatorname*{card}\mathcal{E}_{\mathbf{z}}=2$. Then $\mathbf{z}$ is an
\emph{acute critical point}. Let\footnote{Note that the set of acute critical
points is independent of $\eta$.}%
\[
\mathcal{C}_{\mathcal{T}}^{\operatorname{acute}}:=\left\{  \mathbf{z}%
\in\mathcal{C}_{\mathcal{T}}:\mathbf{z}\text{ is an acute critical
point}\right\}  .
\]

\item $\mathbf{z}\in\mathcal{V}_{\partial\Omega}\left(  \mathcal{T}\right)  $
and $\operatorname*{card}\mathcal{E}_{\mathbf{z}}=3$ and $\Theta\left(
\mathbf{z}\right)  \leq\eta$. Then $\mathbf{z}$ is \emph{flat }$\eta
$\emph{-critical point}. Let%
\[
\mathcal{C}_{\mathcal{T}}^{\operatorname*{flat}}\left(  \eta\right)
:=\left\{  \mathbf{z}\in\mathcal{C}_{\mathcal{T}}\left(  \eta\right)
:\mathbf{z}\text{ is a flat }\eta\text{-critical point}\right\}  .
\]

\item $\mathbf{z}\in\mathcal{V}_{\partial\Omega}\left(  \mathcal{T}\right)  $
and $\operatorname*{card}\mathcal{E}_{\mathbf{z}}=4$ and $\Theta\left(
\mathbf{z}\right)  \leq\eta$. Then $\mathbf{z}$ is a (locally) \emph{concave
}$\eta$-\emph{critical point.} Let%
\[
\mathcal{C}_{\mathcal{T}}^{\operatorname*{concave}}\left(  \eta\right)
:=\left\{  \mathbf{z}\in\mathcal{C}_{\mathcal{T}}\left(  \eta\right)
:\mathbf{z}\text{ is a concave }\eta\text{-critical point}\right\}  .
\]

\end{enumerate}
\end{definition}%

%TCIMACRO{\FRAME{ftbpFU}{5.5737in}{2.9386in}{0pt}{\Qcb{Illustration of the four
%$\eta$-critical cases as in Remark \ref{Remetacritcases}. Left top: inner
%$\eta-$critical point, right top: acute critical point, left bottom: flat
%$\eta$-critical point, right bottom: concave $\eta$-critical point.}%
%}{\Qlb{RemCritGeom2}}{etacritcases.eps}{\special{ language "Scientific Word";
%type "GRAPHIC";  maintain-aspect-ratio TRUE;  display "USEDEF";
%valid_file "F";  width 5.5737in;  height 2.9386in;  depth 0pt;
%original-width 5.5166in;  original-height 2.8954in;  cropleft "0";
%filename 'etacritcases.EPS';file-properties "XNPEU";}}}%
%BeginExpansion
%EndExpansion

The acute critical points require some special treatment and we denote the
union of the others by
\[
\mathcal{C}_{\mathcal{T}}^{\operatorname*{obtuse}}\left(  \eta\right)
:=\mathcal{C}_{\mathcal{T}}^{\operatorname*{inner}}\left(  \eta\right)
\cup\mathcal{C}_{\mathcal{T}}^{\operatorname*{flat}}\left(  \eta\right)
\cup\mathcal{C}_{\mathcal{T}}^{\operatorname*{concave}}\left(  \eta\right)  .
\]

The following lemma states that for a possibly adjusted $\eta_{0}$, still
depending only on the shape-regularity of the mesh and the maximal outer angle
$\alpha_{\Omega}$, the $\eta$-critical points belong to one of the four
categories described in Definition \ref{Remetacritcases}.

\begin{lemma}
\label{Lemeta0}Let $\mathcal{T}$ be a conforming triangulation such that
$D:=\operatorname*{dom}\mathcal{T}$ is a Lipschitz domain.

Then, there exists some $\eta_{0}^{\prime}\in\left]  0,\eta_{0}\right]  $
depending only on the shape-regularity of the mesh and the minimal outer angle
$\alpha_{D}$ such that for $0\leq\eta<\eta_{0}^{\prime}$ any $\eta$-critical
point belongs to one of the four categories described in Definition
\ref{Remetacritcases}.
\end{lemma}

%

%TCIMACRO{\TeXButton{Proof}{\proof}}%
%BeginExpansion
\proof
%EndExpansion
Let $\mathbf{z}$ be an $\eta$-critical point. We set $m:=\operatorname*{card}%
\mathcal{T}_{\mathbf{z}}$ and choose a counterclockwise numbering for the
triangles in $\mathcal{T}_{\mathbf{z}}$, i.e., $K_{i}$, $1\leq i\leq m$. The
shape-regularity of the mesh implies that there is some $m_{\max}$ depending
only on $\phi_{\mathcal{T}}$ such that $m\leq m_{\max}$. Denote by $\omega
_{i}$ the angle in $K_{i}$ at $\mathbf{z}$. Let $\eta_{0}^{\prime}\in\left]
0,\eta_{0}\right]  $ which will be fixed later and assume $0\leq\eta<\eta
_{0}^{\prime}$. Since $\mathbf{z}$ is an $\eta$-critical it holds%
\[
\left\vert \sin\left(  \omega_{i}+\omega_{i+1}\right)  \right\vert \leq
\eta<\eta_{0}^{\prime}\quad\forall1\leq i\leq m^{\prime}\quad\text{for
}m^{\prime}:=\left\{
\begin{array}
[c]{ll}%
m & \text{if }\mathbf{z}\in\mathcal{V}_{\Omega}\left(  \mathcal{T}\right)  ,\\
m-1 & \text{if }\mathbf{z}\in\mathcal{V}_{\partial\Omega}\left(
\mathcal{T}\right)  .
\end{array}
\right.
\]
The shape-regularity implies $\phi_{\mathcal{T}}\leq\omega_{i}\leq\pi
-2\phi_{\mathcal{T}}$ and, for $\delta=\arcsin\eta_{0}^{\prime}$, we get%
\begin{equation}
\omega_{i}+\omega_{i+1}\in\left[  2\phi_{\mathcal{T}},\delta\right]
\cup\left[  \pi-\delta,\pi+\delta\right]  \cup\left[  2\pi-\delta,2\pi
-4\phi_{\mathcal{T}}\right]  \quad\text{for all }1\leq i\leq m^{\prime}.
\label{intervals}%
\end{equation}
Since $\arcsin:\left[  0,1\right[  \rightarrow\mathbb{R}_{\geq0}$ is
monotonously increasing with $\arcsin0=0$ and $\lim_{x\rightarrow1}\arcsin
x=+\infty$, we can select $\eta_{0}^{\prime}$ such $0<\delta<2\phi
_{\mathcal{T}}$. In turn, the first and last interval in (\ref{intervals}) are
empty and%
\begin{equation}
\omega_{i}+\omega_{i+1}=:\pi+\delta_{i}\quad\text{for some }\delta_{i}\text{
with }\left\vert \delta_{i}\right\vert \leq\delta\text{ for all }1\leq i\leq
m^{\prime}. \label{defdeltai}%
\end{equation}

\textbf{Case 1: } $\mathbf{z}\in\mathcal{V}_{\Omega}\left(  \mathcal{T}%
\right)  $.

In this case we obtain%
\begin{equation}
4\pi=\sum_{i=1}^{m}\left(  \omega_{i}+\omega_{i+1}\right)  \overset
{\text{(\ref{defdeltai})}}{=}m\pi+\sum_{i=1}^{m}\delta_{i}. \label{4pieq}%
\end{equation}
By adjusting $\eta_{0}^{\prime}$ such that $m_{\max}\delta<\pi$ we conclude
that $m=4$ and $\sum_{i=1}^{m}\delta_{i}=0$. Hence, $\mathbf{z}$ is an inner
$\eta$-critical point according to Definition \ref{Remetacritcases}(1).

\textbf{Case 2a: }$\mathbf{z}\in\mathcal{V}_{\partial\Omega}\left(
\mathcal{T}\right)  $ and $m\leq3$.

These cases correspond to acute critical/flat $\eta$-critical/concave $\eta
$-critical points according to Definition \ref{Remetacritcases}(2-4).

\textbf{Case 2b: }$\mathbf{z}\in\mathcal{V}_{\partial\Omega}\left(
\mathcal{T}\right)  $ and $m\geq4$.

We argue as in Case \textbf{1} but take into account that the patch
$\mathcal{T}_{\mathbf{z}}$ is not \textquotedblleft closed\textquotedblright%
\ since $\mathbf{z}$ is a boundary point. Let $\alpha:=2\pi-\sum_{i=1}%
^{m}\omega_{i}$ be the \textquotedblleft outer angle\textquotedblright\ of the
domain at $\mathbf{z}$. Then%
\[
\omega_{1}+\omega_{m}+2\sum_{i=2}^{m-1}\omega_{i}=\sum_{i=1}^{m-1}\left(
\omega_{i}+\omega_{i+1}\right)  =\left(  m-1\right)  \pi+\sum_{i=1}%
^{m-1}\delta_{i}.
\]
By the definition of $\alpha$, we obtain%
\begin{align*}
\left(  m-1\right)  \pi+\sum_{i=1}^{m-1}\delta_{i}  &  =2\pi-\alpha+\sum
_{i=2}^{m-1}\omega_{i}=2\pi-\alpha+\sum_{\ell=1}^{\left\lfloor \frac{m-2}%
{2}\right\rfloor }\left(  \omega_{2\ell}+\omega_{2\ell+1}\right)  +Q\left(
m\right)  \omega_{m-1}\\
&  =2\pi-\alpha+\left\lfloor \frac{m-2}{2}\right\rfloor \pi+\sum_{\ell
=1}^{\left\lfloor \frac{m-2}{2}\right\rfloor }\delta_{2\ell}+Q\left(
m\right)  \omega_{m-1},
\end{align*}
where $Q\left(  m\right)  =0$ if $m$ is even and $Q\left(  m\right)  =1$ if
$m$ is odd. By rearranging the terms we get%
\begin{equation}
Q\left(  m\right)  \omega_{m-1}+\Delta_{m}=\left(  m-3-\left\lfloor \frac
{m-2}{2}\right\rfloor \right)  \pi+\alpha\quad\text{for\quad}\Delta_{m}%
:=\sum_{\ell=1}^{\left\lfloor \frac{m-2}{2}\right\rfloor }\delta_{2\ell}%
-\sum_{i=1}^{m-1}\delta_{i}. \label{QMOMEGAEQUAL}%
\end{equation}
We adjust $\eta_{0}^{\prime}$ such that $\delta=\arcsin\eta_{0}^{\prime}$
satisfies $m_{\max}\delta<\alpha$ and, in turn, $\left\vert \Delta
_{m}\right\vert \leq m_{\max}\delta<\alpha$. Then, it is easy to verify that%
\[
\left\vert Q\left(  m\right)  \omega_{m-1}+\Delta_{m}\right\vert <Q\left(
m\right)  \pi+\alpha\leq\left(  m-3-\left\lfloor \frac{m-2}{2}\right\rfloor
\right)  \pi+\alpha
\]
holds for all $m\geq4$. Hence, (\ref{QMOMEGAEQUAL}) cannot hold and there
exists no $\eta$-critical boundary point $\mathbf{z}$ for $m\geq4$.%

%TCIMACRO{\TeXButton{End Proof}{\endproof}}%
%BeginExpansion
\endproof
%EndExpansion

Next, we collect the $\eta-$critical points in pairwise disjoint,
edge-connected sets which we will define in the following.
%TCIMACRO{\FRAME{ftbpFU}{6.5077in}{1.6803in}{0pt}{\Qcb{Three types of obtuse
%$\eta$-critical points $\QTR{bf}{z}\in
%\QTR{cal}{C}_{\QTR{cal}{T}}^{\operatorname*{obtuse}}\left(  \eta\right)  $
%with associated inner edge $\QTR{frak}{E}\left(  \QTR{bf}{z}\right)  $, normal
%vector $\QTR{frak}{N}\left(  \QTR{bf}{z}\right)  $ and opposite endpoint
%$\QTR{frak}{V}\left(  \QTR{bf}{z}\right)  $ of $\QTR{frak}{E}\left(
%\QTR{bf}{z}\right)  $; left: inner $\eta$-critical point, middle: flat $\eta
%$-critical point, right: concave $\eta$-critical point.}}{\Qlb{Figobtuse}%
%}{fan4.eps}{\special{ language "Scientific Word";  type "GRAPHIC";
%width 6.5077in;  height 1.6803in;  depth 0pt;  original-width 7.7089in;
%original-height 1.97in;  cropleft "0";  croptop "1";  cropright "1";
%cropbottom "0";  filename 'fan4.EPS';file-properties "XNPEU";}}}%
%BeginExpansion
\begin{figure}
[ptb]
\begin{center}
\includegraphics[
height=1.6803in,
width=6.5077in
]%
{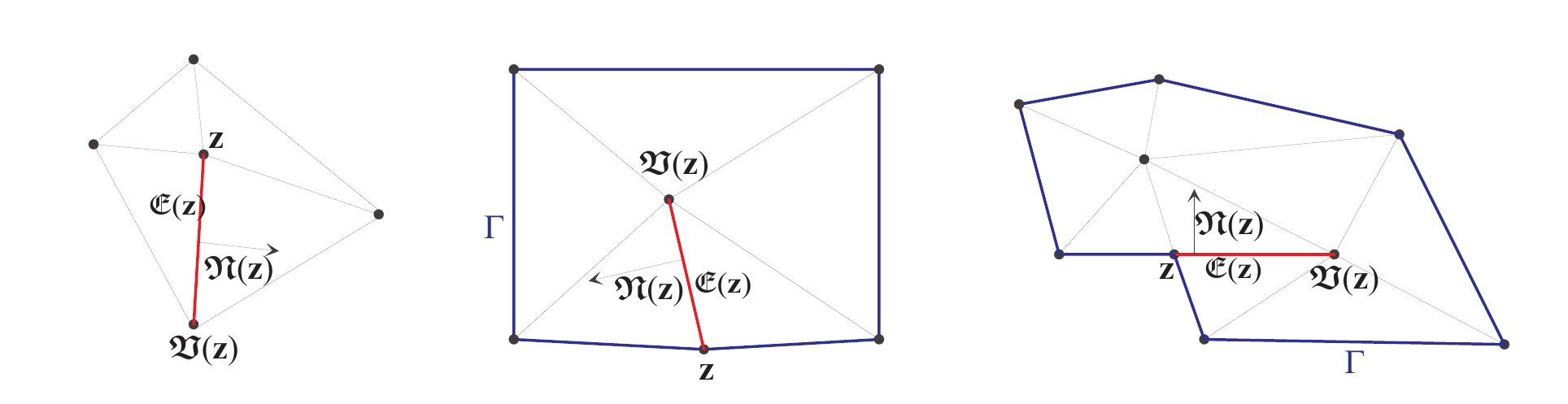}%
\caption{Three types of obtuse $\eta$-critical points $\mathbf{z}%
\in\mathcal{C}_{\mathcal{T}}^{\operatorname*{obtuse}}\left(  \eta\right)  $
with associated inner edge $\mathfrak{E}\left(  \mathbf{z}\right)  $, normal
vector $\mathfrak{N}\left(  \mathbf{z}\right)  $ and opposite endpoint
$\mathfrak{V}\left(  \mathbf{z}\right)  $ of $\mathfrak{E}\left(
\mathbf{z}\right)  $; left: inner $\eta$-critical point, middle: flat $\eta
$-critical point, right: concave $\eta$-critical point.}%
\label{Figobtuse}%
\end{center}
\end{figure}
%EndExpansion
We say two points $\mathbf{y},\mathbf{y}^{\prime}\in\mathcal{V}\left(
\mathcal{T}\right)  $ are \textit{edge-connected} if there is an edge
$E\in\mathcal{E}\left(  \mathcal{T}\right)  $ with endpoints $\mathbf{y,y}%
^{\prime}$. A subset $\mathcal{V}^{\prime}\subset\mathcal{V}\left(
\mathcal{T}\right)  $ is edge-connected if there is a numbering of the points
in $\mathcal{V}^{\prime}=\left\{  \mathbf{y}_{j}:1\leq j\leq n\right\}  $ such
that $\mathbf{y}_{j-1}$, $\mathbf{y}_{j}$ are edge-connected for all $2\leq
j\leq n$. A point $\mathbf{z}\in\mathcal{V}\left(  \mathcal{T}\right)  $ is
edge-connected to $\mathcal{V}^{\prime}$ if $\mathbf{z}\in\mathcal{V}^{\prime
}$ or there is $\mathbf{y\in}\mathcal{V}^{\prime}$ such that $\mathbf{z}$,
$\mathbf{y}$ are edge-connected.

From Lemma \ref{Lemangle} we know that two edge-connected points
$\mathbf{z},\mathbf{z}^{\prime}\in\mathcal{V}\left(  \mathcal{T}\right)  $ can
be both critical only if the connecting edge $E$ belongs to $\mathcal{E}%
_{\partial\Omega}\left(  \mathcal{T}\right)  $; in this case it holds
$\mathbf{z},\mathbf{z}^{\prime}\in\mathcal{V}_{\partial\Omega}\left(
\mathcal{T}\right)  $. Next, we will group the points in $\mathcal{C}%
_{\mathcal{T}}\left(  \eta\right)  $ into subsets called \textit{fans}.

From Lemma \ref{Lemangle} it follows that the points in $\mathcal{C}%
_{\mathcal{T}}^{\operatorname*{inner}}\left(  \eta\right)  $ are isolated (see
Def. \ref{Defetacriticalpoint}). All other $\eta$-critical points lie on the
boundary. Next, we define mappings $\mathfrak{E}:\mathcal{C}_{\mathcal{T}%
}^{\operatorname*{obtuse}}\left(  \eta\right)  \rightarrow\mathcal{E}_{\Omega
}\left(  \mathcal{T}\right)  $, $\mathfrak{N}:\mathcal{C}_{\mathcal{T}%
}^{\operatorname*{obtuse}}\left(  \eta\right)  \rightarrow\mathbb{S}_{2}$, and
$\mathfrak{V}:\mathcal{C}_{\mathcal{T}}^{\operatorname*{obtuse}}\left(
\eta\right)  \rightarrow\mathcal{V}\left(  \mathcal{T}\right)  \backslash
\mathcal{C}_{\mathcal{T}}\left(  \eta\right)  $. The construction is
illustrated in Figure \ref{Figobtuse}.

For $\mathbf{z}\in\mathcal{C}_{\mathcal{T}}^{\operatorname*{obtuse}}\left(
\eta\right)  $, Definition \ref{Remetacritcases} implies that $\left\vert
\mathcal{E}_{\mathbf{z}}\right\vert \geq3$ and hence $\mathcal{E}_{\mathbf{z}%
}\cap\mathcal{E}_{\Omega}\left(  \mathcal{T}\right)  \neq\emptyset$. We fix
one edge $E\in\mathcal{E}_{\mathbf{z}}\cap\mathcal{E}_{\Omega}\left(
\mathcal{T}\right)  $ and set $\mathfrak{E}\left(  \mathbf{z}\right)  :=E$.
Note that the choice of $E$ is unique for $\mathbf{z}\in\mathcal{C}%
_{\mathcal{T}}^{\operatorname*{flat}}\left(  \eta\right)  $. For
$\mathbf{z}\in\mathcal{C}_{\mathcal{T}}^{\operatorname*{inner}}\left(
\eta\right)  $ the choice is arbitrary. For $\mathbf{z}\in\mathcal{C}%
_{\mathcal{T}}^{\operatorname*{concave}}\left(  \eta\right)  $, the set
$\mathcal{E}_{\mathbf{z}}\cap\mathcal{E}_{\Omega}\left(  \mathcal{T}\right)  $
consists of two edges, say $E_{1}$, $E_{2}$. We fix one of them and set
$\mathfrak{E}\left(  \mathbf{z}\right)  :=E_{2}$. Let $\mathbf{z}^{\prime}%
\in\mathcal{V}\left(  \mathcal{T}\right)  $ be such that $\mathfrak{E}\left(
\mathbf{z}\right)  =\left[  \mathbf{z},\mathbf{z}^{\prime}\right]  $. Then
$\mathfrak{V}\left(  \mathbf{z}\right)  :=\mathbf{z}^{\prime}$. Lemma
\ref{Lemangle} implies that $\mathbf{z}^{\prime}$ is not an $\eta$-critical
point. A unit vector $\mathfrak{N}\left(  \mathbf{z}\right)  $ orthogonal to
$\mathfrak{E}\left(  \mathbf{z}\right)  $ is defined by the condition that
$\mathbf{z}^{\prime}-\mathbf{z}$ and $\mathfrak{N}\left(  \mathbf{z}\right)  $
form a right-handed system.

\begin{definition}
\label{DefFan}We decompose $C_{\mathcal{T}}^{\operatorname*{obtuse}}\left(
\eta\right)  $ into disjoint \emph{fans} $\mathcal{C}_{\mathcal{T},\ell
}\left(  \eta\right)  $, $\ell\in\mathcal{J}$, such that the following
conditions are satisfies

\begin{enumerate}
\item $\mathcal{C}_{\mathcal{T}}^{\operatorname*{obtuse}}\left(  \eta\right)
=%
%TCIMACRO{\dbigcup \limits_{\ell\in\mathcal{J}}}%
%BeginExpansion
{\displaystyle\bigcup\limits_{\ell\in\mathcal{J}}}
%EndExpansion
\mathcal{C}_{\mathcal{T},\ell}\left(  \eta\right)  ,$

\item for any $\ell\in\mathcal{J}$, the set $\mathcal{C}_{\mathcal{T},\ell
}\left(  \eta\right)  $ is edge-connected,

\item for any $\ell\in\mathcal{J}$, there is $\mathbf{z}_{\ell}\in
\mathcal{V}\left(  \mathcal{T}\right)  \backslash\mathcal{C}_{\mathcal{T}%
}\left(  \eta\right)  $ such that for all $\mathbf{z}\in\mathcal{C}%
_{\mathcal{T},\ell}\left(  \eta\right)  $ it holds $\mathfrak{V}\left(
\mathbf{z}\right)  =\mathbf{z}_{\ell}$ and, vice versa:

\item any $\mathbf{z}^{\prime}\in\mathcal{C}_{\mathcal{T}}%
^{\operatorname*{obtuse}}\left(  \eta\right)  $ which is edge-connected to
some $\mathcal{C}_{\mathcal{T},\ell}\left(  \eta\right)  $ and satisfies
$\mathfrak{V}\left(  \mathbf{z}^{\prime}\right)  =\mathbf{z}_{\ell}$ belongs
to $\mathcal{C}_{\mathcal{T},\ell}\left(  \eta\right)  $.
\end{enumerate}
\end{definition}

The following lemma will allow us to construct a right-inverse for the
divergence operator separately for each fan.

\begin{lemma}
\label{Lemomegaell}Let $\eta_{0}$ be as in in Lemma \ref{Lemangle} and let
$0\leq\eta<\eta_{0}$ be fixed.

\begin{enumerate}
\item[a.] Then, the mapping $\mathfrak{E}:\mathcal{C}_{\mathcal{T}%
}^{\operatorname*{obtuse}}\left(  \eta\right)  \rightarrow\mathcal{E}_{\Omega
}\left(  \mathcal{T}\right)  $ is injective.

\item[b.] For $\ell\in\mathcal{J}$, let $\omega_{\ell}:=%
%TCIMACRO{\dbigcup \limits_{\mathbf{z}\in\mathcal{C}_{\mathcal{T},\ell}\left(
%\eta\right)  }}%
%BeginExpansion
{\displaystyle\bigcup\limits_{\mathbf{z}\in\mathcal{C}_{\mathcal{T},\ell
}\left(  \eta\right)  }}
%EndExpansion
\omega_{\mathfrak{E}\left(  \mathbf{z}\right)  }$. The domains $\omega_{\ell}$
have pairwise disjoint interior.
\end{enumerate}
\end{lemma}

%

%TCIMACRO{\TeXButton{Proof}{\proof}}%
%BeginExpansion
\proof
%EndExpansion
\textbf{Part} \textbf{a. }The injectivity of the mapping $\mathfrak{E}%
:\mathcal{C}_{\mathcal{T}}^{\operatorname*{obtuse}}\left(  \eta\right)
\rightarrow\mathcal{E}_{\Omega}\left(  \mathcal{T}\right)  $ follows from
Lemma \ref{Lemangle}: if $\mathbf{z}\in\mathcal{C}_{\mathcal{T}}\left(
\eta\right)  $ and $\mathbf{z}^{\prime}\in\mathcal{V}\left(  \mathcal{T}%
\right)  $ is such that $E:=\left[  \mathbf{z},\mathbf{z}^{\prime}\right]
\in\mathcal{E}_{\Omega}\left(  \mathcal{T}\right)  $ then $\mathbf{z}^{\prime
}\notin\mathcal{C}_{\mathcal{T}}\left(  \eta\right)  $.

\begin{figure}
[ptb]
\begin{center}
\includegraphics[
height=3.0761in,
width=5.7432in
]%
{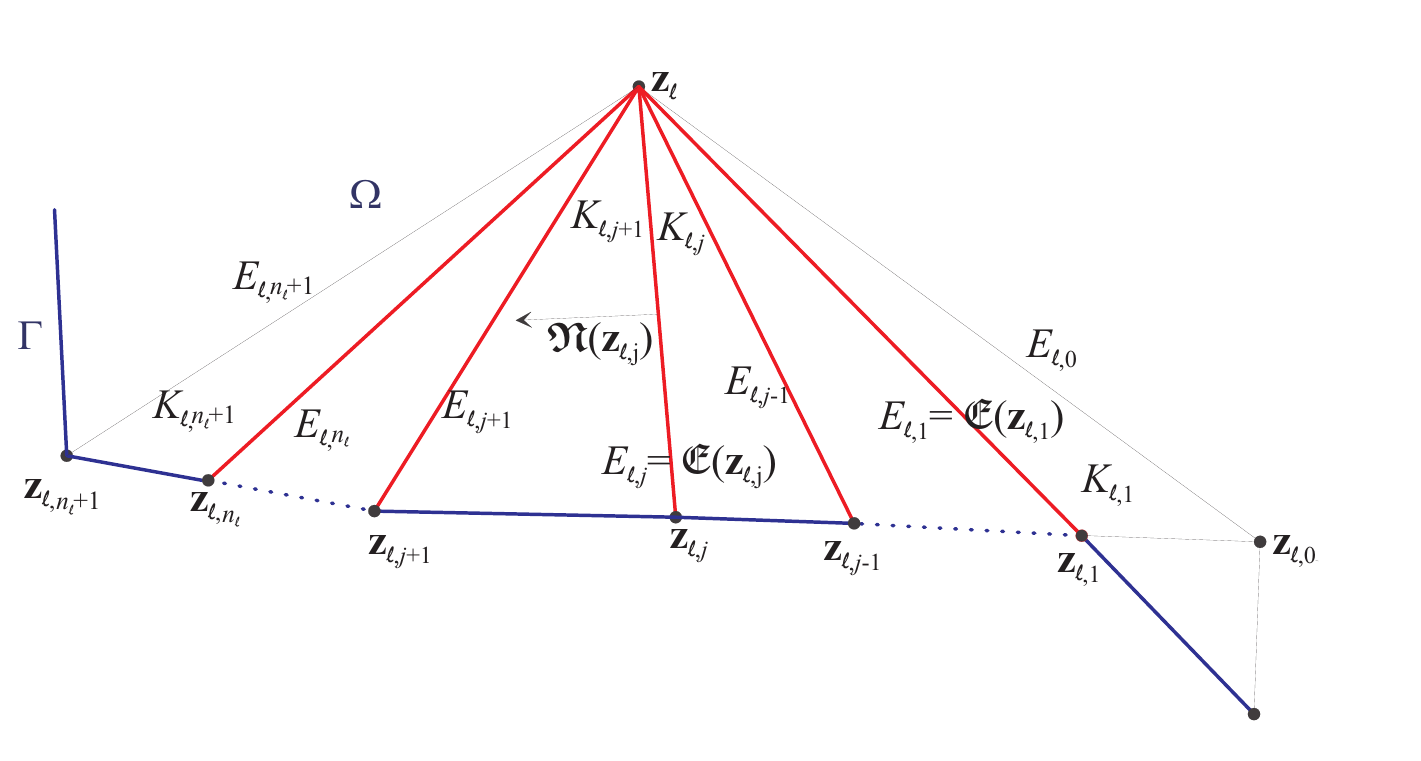}%
\caption{Nodal patch, illustrating edge-connected obtuse $\eta$-critical
points of a fan $\mathcal{C}_{\mathcal{T},\ell}\left(  \eta\right)  $. In this
example, the left-most $\eta$-critical point is $\mathbf{z}_{\ell,n_{\ell}}$
and of type \textquotedblleft flat\textquotedblright, the right-most is
$\mathbf{z}_{\ell,1}$ of type \textquotedblleft concave\textquotedblright%
\ with $\mathfrak{E}\left(  \mathbf{z}_{\ell,1}\right)  =\left[
\mathbf{z}_{\ell,1},\mathbf{z}\right]  $. The extremal points $\mathbf{z}%
_{\ell,0}$ and $\mathbf{z}_{\ell,n_{\ell}+1}$ do not belong to $\mathcal{C}%
_{\mathcal{T},\ell}\left(  \eta\right)  $. The edge connecting $\mathbf{z}%
_{\ell}$ with $\mathbf{z}_{\ell,j}$ is denoted by $E_{\ell,j}$. }%
\label{Fig_nod_patch}%
\end{center}
\end{figure}

\textbf{Part b. }The following construction is illustrated in Figure
\ref{Fig_nod_patch}.%
%TCIMACRO{\FRAME{ftbpFU}{5.7432in}{3.0761in}{0pt}{\Qcb{Nodal patch,
%illustrating edge-connected obtuse $\eta$-critical points of a fan
%$\QTR{cal}{C}_{\QTR{cal}{T},\ell}\left(  \eta\right)  $. In this example, the
%left-most $\eta$-critical point is $\QTR{bf}{z}_{\ell,n_{\ell}}$ and of type
%\textquotedblleft flat\textquotedblright, the right-most is $\QTR{bf}{z}_{\ell
%,1}$ of type \textquotedblleft concave\textquotedblright\ with
%$\QTR{frak}{E}\left(  \QTR{bf}{z}_{\ell,1}\right)  =\left[  \QTR{bf}{z}_{\ell
%,1},\QTR{bf}{z}\right]  $. The extremal points $\QTR{bf}{z}_{\ell,0}$ and
%$\QTR{bf}{z}_{\ell,n_{\ell}+1}$ do not belong to
%$\QTR{cal}{C}_{\QTR{cal}{T},\ell}\left(  \eta\right)  $. The edge connecting
%$\QTR{bf}{z}_{\ell}$ with $\QTR{bf}{z}_{\ell,j}$ is denoted by $E_{\ell,j}$.
%}}{\Qlb{Fig_nod_patch}}{fan5.eps}{\special{ language "Scientific Word";
%type "GRAPHIC";  maintain-aspect-ratio TRUE;  display "USEDEF";
%valid_file "F";  width 5.7432in;  height 3.0761in;  depth 0pt;
%original-width 5.6844in;  original-height 3.0329in;  cropleft "0";
%filename 'fan5.EPS';file-properties "XNPEU";}}}%
%BeginExpansion
%EndExpansion
For fixed $\ell\in\mathcal{J}$, we number the points in $\mathcal{C}%
_{\mathcal{T},\ell}$ by $\mathbf{z}_{\ell,j}$, $1\leq j\leq n_{\ell}$, such
that $\mathbf{z}_{\ell,j-1}$ and $\mathbf{z}_{\ell,j}$ are edge-connected for
all $2\leq j\leq n_{\ell}$ and $\mathbf{z}_{\ell,1}$ and $\mathbf{z}%
_{\ell,n_{\ell}}$ are the endpoints in the polygonal line through these
points. Let $K_{\ell,j}$ be the triangle with vertices $\mathbf{z}_{\ell,j-1}%
$, $\mathbf{z}_{\ell}$, $\mathbf{z}_{\ell,j}$, $2\leq j\leq n_{\ell}\ $and let
$\mathcal{T}_{\ell}:=\left\{  K_{\ell,j}:2\leq j\leq n_{\ell}\right\}  $. Note
that this set is empty if $\left\vert \mathcal{C}_{\mathcal{T},\ell
}\right\vert =1$. Let $K_{\ell,1},K_{\ell,n_{\ell}+1}\in\mathcal{T}%
\backslash\mathcal{T}_{\ell}$ be two different triangles such that $E_{\ell
,1}:=\left[  \mathbf{z}_{\ell,1},\mathbf{z}_{\ell}\right]  \subset\partial
K_{\ell,1}$ and $E_{\ell,n_{\ell}}:=\left[  \mathbf{z}_{\ell,n_{\ell}%
},\mathbf{z}_{\ell}\right]  \subset\partial K_{\ell,n_{\ell}+1}$. Let
$\mathbf{z}_{\ell,0}$ be the third vertex in $K_{\ell,1}$ and observe that it
does not belong to $\mathcal{C}_{\mathcal{T},\ell}$. Since $E_{\ell,1}$ is an
inner edge and $\mathbf{z}_{\ell,1}$ is an $\eta$-critical point, Lemma
\ref{Lemangle} implies that $\mathbf{z}_{\ell}$ is not an $\eta$-critical
point. Next, we show that $E_{\ell,0}:=\left[  \mathbf{z}_{\ell,0}%
,\mathbf{z}_{\ell}\right]  $ does not belong to $\mathcal{E}_{\mathcal{T}%
}^{\operatorname*{obtuse}}:=\left\{  \mathfrak{E}\left(  \mathbf{z}\right)
:\mathbf{z}\in\mathcal{C}_{\mathcal{T}}^{\operatorname*{obtuse}}\left(
\eta\right)  \right\}  $; from this the assertion follows. We assume
$E_{\ell,0}\in\mathcal{E}_{\mathcal{T}}^{\operatorname*{obtuse}}$ and derive a
contradiction. Since $\mathbf{z}_{\ell}\notin\mathcal{C}_{\mathcal{T}%
}^{\operatorname*{obtuse}}\left(  \eta\right)  $, this assumption implies that
$\mathbf{z}_{\ell,0}\in\mathcal{C}_{\mathcal{T}}^{\operatorname*{obtuse}%
}\left(  \eta\right)  $. If $\mathfrak{V}\left(  \mathbf{z}_{\ell,0}\right)
=\mathbf{z}_{\ell}$ then Definition \ref{DefFan}(4) implies that
$\mathbf{z}_{\ell,0}\in\mathcal{C}_{\mathcal{T},\ell}$ and this is a
contradiction. If $\mathfrak{V}\left(  \mathbf{z}_{\ell,0}\right)
\neq\mathbf{z}_{\ell}$, then $E_{\ell,0}\notin\mathcal{E}_{\mathcal{T}%
}^{\operatorname*{obtuse}}$.
%TCIMACRO{\TeXButton{End Proof}{\endproof}}%
%BeginExpansion
\endproof
%EndExpansion

Since the acute critical points need some special treatment we define a
sequence of triangulations $\mathcal{T}_{i}$, $1\leq i\leq L$, with the properties

\begin{enumerate}
\item
\begin{equation}
\mathcal{T}_{1}\subset\mathcal{T}_{2}\subset\ldots\subset\mathcal{T}%
_{L}=\mathcal{T}, \label{stepbystepextension}%
\end{equation}

\item $\mathcal{T}_{1}$ is a maximal subset of $\mathcal{T}$ such that
$\mathcal{C}_{\mathcal{T}_{1}}^{\operatorname{acute}}=\emptyset$,

\item for $j=1,2,\ldots,L$,%
\[
\mathcal{T}_{j}=\left\{  K\in\mathcal{T}\mid\text{at least one edge of
}K\text{ belongs to }\mathcal{E}\left(  \mathcal{T}_{j-1}\right)  \right\}  .
\]

\end{enumerate}

By this step-by-step procedure, triangles are attached to a previous
triangulation $\mathcal{T}_{j-1}$ which have an edge in common with the set of
edges in $\mathcal{T}_{j-1}$. The proof of Theorem \ref{Theomain} under
assumption (\ref{Assumptab}) then consists of first proving the inf-sup
stability for $\mathcal{T}_{1}$ and then to investigate the effect of
attaching a triangle to an inf-sup stable triangulation. A sufficient
condition for $L=0$ is that every triangle in $\mathcal{T}$ has an interior point.

\subsubsection{The case $\mathcal{C}_{\mathcal{T}}^{\operatorname{acute}%
}=\emptyset$\label{Obtuse}}

In this section, we prove the inf-sup stability for the triangulation
$\mathcal{T}_{1}$ in (\ref{stepbystepextension}) where $\mathcal{C}%
_{\mathcal{T}_{1}}^{\operatorname{acute}}=\emptyset$. For simplicity we skip
the index $1$ and write $\mathcal{T}$, $\mathcal{C}_{\mathcal{T}}\left(
\eta\right)  $, etc.

Next, we define some fundamental non-conforming Crouzeix-Raviart vector fields
which will be used to eliminate the critical pressures in the Stokes element
$\left(  \mathbf{S}_{k,0}\left(  \mathcal{T}\right)  ,\mathbb{P}%
_{k-1,0}\left(  \mathcal{T}\right)  \right)  $.

Essential properties of the Crouzeix-Raviart function $B_{k,E}%
^{\operatorname*{CR}}$ are: it is a polynomial of degree $k$ on each
$K\subset\mathcal{T}_{E}$ and a Legendre polynomial on each edge $E^{\prime
}\subset\partial\omega_{E}$ so that the jump relations in (\ref{PCR0RB}) are
satisfied. Furthermore, $\left[  B_{k,E}^{\operatorname*{CR}}\right]  _{E}=0$
and $\operatorname*{supp}B_{k,E}^{\operatorname*{CR}}\subset\omega_{E}$.

In the first step, we modify the function $B_{k,E}^{\operatorname*{CR}}$ by
adding a conforming edge bubble in $S_{k,0}\left(  \mathcal{T}\right)  $ such
that the $H^{1}\left(  \Omega\right)  $ norm of the modified function has an
improved behaviour with respect to $k$.

Let $E\in\mathcal{E}_{\Omega}\left(  \mathcal{T}\right)  $ with endpoints
$\mathbf{V}_{1}$, $\mathbf{V}_{2}$. Set $\mathbf{t}_{E}=\left(  \mathbf{V}%
_{2}-\mathbf{V}_{1}\right)  /\left\Vert \mathbf{V}_{2}-\mathbf{V}%
_{1}\right\Vert $ and consider a function $w_{E}\in\mathbb{P}_{k}\left(
\mathcal{T}\right)  $ with $\operatorname*{supp}w_{E}=\omega_{E}$ and%
%TCIMACRO{\TeXButton{Defweall}{\begin{subequations}
%\label{Defweall}
%\end{subequations}}}%
%BeginExpansion
\begin{subequations}
\label{Defweall}
\end{subequations}%
%EndExpansion%
\begin{align}
\left.  \left.  w_{E}\right\vert _{K}\right\vert _{E^{\prime}} &  =\left.
\left.  B_{k,E}^{\operatorname*{CR}}\right\vert _{K}\right\vert _{E^{\prime}%
}\quad\forall K\in\mathcal{T}_{E}\text{ and }E^{\prime}\subset\partial
K\cap\partial\omega_{E},\tag{%
%TCIMACRO{\TeXButton{Defweall}{\ref{Defweall}}}%
%BeginExpansion
\ref{Defweall}%
%EndExpansion
a}\label{Defwe}\\
\left[  w_{E}\right]  _{E} &  =0\quad\text{and\quad}\partial_{\mathbf{t}_{E}%
}w_{E}\left(  \mathbf{V}_{1}\right)  =\partial_{\mathbf{t}_{E}}w_{E}\left(
\mathbf{V}_{2}\right)  =0.\tag{%
%TCIMACRO{\TeXButton{Defweall}{\ref{Defweall}}}%
%BeginExpansion
\ref{Defweall}%
%EndExpansion
b}\label{Defweallb}%
\end{align}
Then, $w_{E}$ also belongs to the space $\operatorname*{CR}_{k,0}\left(
\mathcal{T}\right)  $ and%
\begin{equation}
\nabla\left(  \left.  w_{E}\right\vert _{K}\right)  \left(  \mathbf{z}\right)
=\nabla\left(  \left.  B_{k,E}^{\operatorname*{CR}}\right\vert _{K}\right)
\left(  \mathbf{z}\right)  \quad\forall K\in\mathcal{T\quad\forall}%
\mathbf{z}\in\mathcal{V}\left(  K\right)  .\label{gradrel}%
\end{equation}
The last relation can be derived from the following reasoning. For
$K\in\mathcal{T}$ and $\mathbf{z}\in\mathcal{V}\left(  K\right)  $, set
$\mathbf{t}_{\mathbf{y}}:=\left(  \mathbf{y}-\mathbf{z}\right)  /\left\Vert
\mathbf{y-z}\right\Vert $ for all $\mathbf{y}\in\mathcal{V}\left(  K\right)
\backslash\left\{  \mathbf{z}\right\}  $. Let $\mathbf{c}\in\mathbb{R}^{2}$ be
arbitrary. Clearly, $\mathbf{c}=\sum_{\mathbf{y}\in\mathcal{V}\left(
K\right)  \backslash\left\{  \mathbf{z}\right\}  }\alpha_{\mathbf{y}%
}\mathbf{t}_{\mathbf{y}}$ for some $\alpha_{\mathbf{y}}\in\mathbb{R}$. The
conditions in (\ref{Defweall}) imply that for $\mathbf{y}\in\mathcal{V}\left(
K\right)  \backslash\left\{  \mathbf{z}\right\}  $ it holds%
\[
\frac{\partial\left.  w_{E}\right\vert _{K}}{\partial\mathbf{t}_{\mathbf{y}}%
}=\frac{\partial\left.  B_{k,E}^{\operatorname*{CR}}\right\vert _{K}}%
{\partial\mathbf{t}_{\mathbf{y}}}.
\]
Hence,%
\begin{align}
\left\langle \nabla\left(  \left.  w_{E}\right\vert _{K}\right)
,\mathbf{c}\right\rangle \left(  \mathbf{z}\right)   &  =\sum_{\mathbf{y}%
\in\mathcal{V}\left(  K\right)  \backslash\left\{  \mathbf{z}\right\}  }%
\alpha_{\mathbf{y}}\frac{\partial\left.  w_{E}\right\vert _{K}}{\partial
\mathbf{t}_{\mathbf{y}}}\left(  \mathbf{z}\right)  \label{reasoningz}\\
&  =\sum_{\mathbf{y}\in\mathcal{V}\left(  K\right)  \backslash\left\{
\mathbf{z}\right\}  }\alpha_{\mathbf{y}}\frac{\partial\left.  B_{k,E}%
^{\operatorname*{CR}}\right\vert _{K}}{\partial\mathbf{t}_{\mathbf{y}}}\left(
\mathbf{z}\right)  =\left\langle \nabla\left(  \left.  \left(  B_{k,E}%
^{\operatorname*{CR}}\right)  \right\vert _{K}\right)  ,\mathbf{c}%
\right\rangle \left(  \mathbf{z}\right)  .\nonumber
\end{align}
Since $\mathbf{c}$ was arbitrary, (\ref{gradrel}) follows.

\begin{lemma}
\label{Lemtrianglemean}Let $k\geq5$ be odd and for $E\in\mathcal{E}_{\Omega
}\left(  \mathcal{T}\right)  $, let $B_{k,E}^{\operatorname*{CR}}$ be as in
(\ref{BkEdef}). Then, there exists a function $\tilde{B}_{k,E}%
^{\operatorname*{CR}}\in\operatorname*{CR}_{k,0}\left(  \mathcal{T}\right)  $ with

\begin{enumerate}
\item[i.] $\operatorname*{supp}\tilde{B}_{k,E}^{\operatorname*{CR}}=\omega
_{E},$

\item[ii.] for all $K\in\mathcal{T}$, for all $\mathbf{z}\in\mathcal{V}\left(
K\right)  $:%
\begin{equation}
\nabla\left(  \left.  \tilde{B}_{k,E}^{\operatorname*{CR}}\right\vert
_{K}\right)  \left(  \mathbf{z}\right)  =\nabla\left(  \left.  B_{k,E}%
^{\operatorname*{CR}}\right\vert _{K}\right)  \left(  \mathbf{z}\right)  ,
\label{neinher}%
\end{equation}

\item[iii.] for all $K\in\mathcal{T}_{E}$, for all $E^{\prime}\subset\partial
K\cap\mathcal{\partial}\omega_{E}$:%
\[
\left.  \left.  \tilde{B}_{k,E}^{\operatorname*{CR}}\right\vert _{K}%
\right\vert _{E^{\prime}}=\left.  \left.  B_{k,E}^{\operatorname*{CR}%
}\right\vert _{K}\right\vert _{E^{\prime}}\quad\text{and\quad}\left[
\tilde{B}_{k,E}^{\operatorname*{CR}}\right]  _{E}=0,
\]

\item[iv.] for all $K\in\mathcal{T}$%
\begin{equation}
\int_{K}\operatorname{div}_{\mathcal{T}}\left(  \tilde{B}_{k,E}%
^{\operatorname*{CR}}\mathbf{n}_{E}\right)  =0. \label{localintzero}%
\end{equation}

\item[v.] The piecewise gradient is bounded:%
\begin{equation}
\left\Vert \nabla_{\mathcal{T}}\tilde{B}_{k,E}^{\operatorname*{CR}}\right\Vert
_{\mathbf{L}^{2}\left(  \Omega\right)  }\leq C\sqrt{\log\left(  k+1\right)  }.
\label{pwgrad}%
\end{equation}

\end{enumerate}
\end{lemma}%

%TCIMACRO{\TeXButton{Proof}{\proof}}%
%BeginExpansion
\proof
%EndExpansion
We employ the reference triangle as in \cite{BabuskaMandel} in order to apply
the polynomial extension theorem therein. Let $\tilde{K}$ be the equilateral
triangle with vertices $\mathbf{\tilde{A}}_{1}:=\left(  -1,0\right)  ^{T}$,
$\mathbf{\tilde{A}}_{2}:=\left(  1,0\right)  ^{T}$, $\mathbf{\tilde{A}}%
_{3}:=\left(  0,\sqrt{3}\right)  ^{T}$ and let $\tilde{E}_{j}$ denote the edge
in $\tilde{K}$ opposite to $\mathbf{\tilde{A}}_{j}$, $1\leq j\leq3$.

Let $E\in\mathcal{E}_{\Omega}\left(  \mathcal{T}\right)  $ with endpoints
$\mathbf{V}_{1}$, $\mathbf{V}_{2}$, and let $K\in\mathcal{T}_{E}$. Choose an
affine pullback $\phi_{K}:\tilde{K}\rightarrow K$ such that $\phi_{K}\left(
\tilde{E}_{3}\right)  =E$. We employ the function $\tilde{\psi}_{k}^{\pm}%
\in\mathbb{P}_{k}\left(  \left[  -1,1\right]  \right)  $ given by%
\[
\tilde{\psi}_{k}^{-}\left(  x\right)  :=c_{k}^{-1}\frac{\left(  1+x\right)
\left(  1-x\right)  ^{2}}{4}P_{k-3}^{\left(  3,3\right)  }\left(  x\right)
\quad\text{and\quad}\tilde{\psi}_{k}^{+}\left(  x\right)  :=-\tilde{\psi}%
_{k}^{-}\left(  -x\right)
\]
with the Jacobi polynomials $P_{n}^{\left(  \alpha,\beta\right)  }$ (see,
e.g., \cite[\S 18.3]{NIST:DLMF}) and the normalisation factor $c_{k}:=\left(
-1\right)  ^{k-1}\binom{k}{3}$. These functions have been analysed in the
proof of Lemma A.1 in \cite[denoted by $F_{k}$]{Ainsworth_parker_I} and we
recall relevant properties. It holds $\tilde{\psi}_{k}^{\pm}\left(
\pm1\right)  =\left(  \tilde{\psi}_{k}^{-}\right)  ^{\prime}\left(  +1\right)
=\left(  \tilde{\psi}_{k}^{+}\right)  ^{\prime}\left(  -1\right)  =0$ and
$\left(  \tilde{\psi}_{k}^{-}\right)  ^{\prime}\left(  -1\right)  =\left(
\tilde{\psi}_{k}^{+}\right)  ^{\prime}\left(  +1\right)  =1$ (cf.
\cite[\S 18.3]{NIST:DLMF}). Their norms can be estimated by%
\[
\left\Vert \tilde{\psi}_{k}^{\pm}\right\Vert _{L^{2}\left(  \left[
-1,1\right]  \right)  }\leq Ck^{-3}\quad\text{and\quad}\left\Vert \left(
\tilde{\psi}_{k}^{\pm}\right)  ^{\prime}\right\Vert _{L^{2}\left(  \left[
-1,1\right]  \right)  }\leq Ck^{-1}.
\]
We set
\begin{equation}
\tilde{\varphi}_{k}\left(  x\right)  =L_{k-1}\left(  x\right)  -L_{k-1}%
^{\prime}\left(  -1\right)  \tilde{\psi}_{k}^{-}\left(  x\right)
-L_{k-1}^{\prime}\left(  1\right)  \tilde{\psi}_{k}^{+}\left(  x\right)  .
\label{defphiktilde}%
\end{equation}
Clearly, it holds%
\[
\tilde{\varphi}_{k}\left(  \pm1\right)  =1,\quad\tilde{\varphi}_{k}^{\prime
}\left(  \pm1\right)  =0.
\]
By using Lemma \ref{LemExInt}, we get%
\[
\left\Vert \tilde{\varphi}_{k}\right\Vert _{L^{2}\left(  \left[  -1,1\right]
\right)  }\leq Ck^{-1/2},\quad\left\Vert \tilde{\varphi}_{k}\right\Vert
_{H^{1}\left(  \left[  -1,1\right]  \right)  }\leq Ck.
\]
From \cite[Thm. 7.4]{BabuskaMandel} we conclude that there is $\tilde{w}%
_{E}\in\mathbb{P}_{k}\left(  \tilde{K}\right)  $ with
\[
\left.  \tilde{w}_{E}\right\vert _{\tilde{E}_{j}}=\left.  \left.
B_{k,E}^{\operatorname*{CR}}\right\vert _{K}\circ\phi_{K}\right\vert
_{\tilde{E}_{j}},\quad j=1,2\quad\text{and\quad}\left.  \tilde{w}%
_{E}\right\vert _{\tilde{E}_{3}}=\tilde{\varphi}_{k}%
\]
which satisfies%
\[
\left\Vert \tilde{w}_{E}\right\Vert _{H^{1}\left(  \tilde{K}\right)  }\leq
C\left\Vert \tilde{w}_{E}\right\Vert _{H^{1/2}\left(  \partial\tilde
{K}\right)  }%
\]
for a constant $C$ independent of $k$. Lemma \ref{LemBEH1/2norm} implies the
following estimate of the $H^{1/2}$ norm of $\tilde{w}_{E}$:%
\begin{equation}
\left\Vert \tilde{w}_{E}\right\Vert _{H^{1/2}\left(  \partial\tilde{K}\right)
}\leq C\sqrt{\log\left(  k+1\right)  }. \label{wtildeE}%
\end{equation}
In turn, we get%
\begin{equation}
\left\Vert \tilde{w}_{E}\right\Vert _{H^{1}\left(  \tilde{K}\right)  }\leq
C\sqrt{\log\left(  k+1\right)  }. \label{wtildeest}%
\end{equation}
By using the affine lifting $\phi_{K}$ to the triangle $K$ we define the
function $w_{E}$ by%
\[
\left.  w_{E}\right\vert _{K}:=\left\{
\begin{array}
[c]{ll}%
\tilde{w}_{E}\circ\phi_{K}^{-1} & \text{if }K\in\mathcal{T}_{E}\text{,}\\
0 & \text{otherwise.}%
\end{array}
\right.
\]
This function is continuous across $E$ (with value $\tilde{\varphi}_{k}%
\circ\left.  \phi_{K}^{-1}\right\vert _{E}$) and, on $E^{\prime}%
\subset\partial\omega_{E}$, it is a lifted Legendre polynomial. This implies
property (iii) for $w_{E}$. The function $w_{E}$ vanishes outside $\omega_{E}$
so that (i) holds. Since the construction implies that the derivative of
$w_{E}$ in the direction of $E$, evaluated at the endpoints $\mathbf{V}_{1}$,
$\mathbf{V}_{2}$ of $E$, is zero, we may apply the reasoning in
(\ref{reasoningz}) to obtain property (ii) for $w_{E}$.

From (\ref{wtildeest}) we obtain by the transformation rule for integrals and
the chain rule for differentiation%
\[
\left\Vert \nabla w_{E}\right\Vert _{\mathbf{L}^{2}\left(  K\right)  }\leq
C\left\Vert \tilde{w}_{E}\right\Vert _{H^{1}\left(  \tilde{K}\right)  }\leq
C\sqrt{\log\left(  k+1\right)  }\quad\forall K\in\mathcal{T}_{E}\text{.}%
\]
Next, we modify $w_{E}$ such that property (iv) holds without affecting the
other properties. Let $\psi_{E}\in S_{4,0}\left(  \mathcal{T}\right)  $ with
$\operatorname*{supp}\psi_{E}=\omega_{E}$ and
\begin{equation}
\left.  \psi_{E}\right\vert _{K}:=\alpha_{K}\lambda_{K,\mathbf{V}_{1}}%
^{2}\lambda_{K,\mathbf{V}_{2}}^{2}\quad\text{with\quad}\alpha_{K}:=\left(
\int_{K}\partial_{\mathbf{n}_{E}}w_{E}\right)  /\left(  \int_{K}%
\partial_{\mathbf{n}_{E}}\left(  \lambda_{K,\mathbf{V}_{1}}^{2}\lambda
_{K,\mathbf{V}_{2}}^{2}\right)  \right)  \quad\forall K\in\mathcal{T}_{E}.
\label{DefPsieta}%
\end{equation}
The modified function $\tilde{B}_{k,E}^{\operatorname*{CR}}$ finally is
defined by%
\begin{equation}
\tilde{B}_{k,E}^{\operatorname*{CR}}=w_{E}-\psi_{E}. \label{defboldpsi}%
\end{equation}
Since $\left.  \operatorname*{div}_{\mathcal{T}}\left(  \tilde{B}%
_{k,E}^{\operatorname*{CR}}\mathbf{n}_{E}\right)  \right\vert _{K}%
=\partial_{\mathbf{n}_{E}}\left(  \tilde{B}\left.  _{k,E}^{\operatorname*{CR}%
}\right\vert _{K}\right)  $ property (iv) follows by construction. The
gradient $\nabla_{\mathcal{T}}\psi_{E}$ vanishes in the vertices of $K$ so
that $\left(  \partial_{\mathbf{n}_{E}}\psi_{E}\right)  \left(  \mathbf{z}%
\right)  =0$ for all $\mathbf{z}\in\mathcal{V}\left(  K\right)  $ and (ii) is
inherited from $w_{E}$. Properties (i), (iii) are obvious. Next, we verify
(v). Let $\mathbf{V}_{3}$ denote the vertex in $K$ opposite to $E$. We first
compute%
\begin{align*}
\int_{K}\partial_{\mathbf{n}_{E}}\left(  \lambda_{K,\mathbf{V}_{1}}^{2}%
\lambda_{K,\mathbf{V}_{2}}^{2}\right)   &  =\sum_{j=1}^{2}2\partial
_{\mathbf{n}_{E}}\lambda_{K,\mathbf{V}_{j}}\int_{K}\lambda_{K,\mathbf{V}_{1}%
}\lambda_{K,\mathbf{V}_{2}}\lambda_{K,\mathbf{V}_{3-j}}=\int_{K}%
\lambda_{K,\mathbf{V}_{1}}^{2}\lambda_{K,\mathbf{V}_{2}}\sum_{j=1}%
^{2}2\partial_{\mathbf{n}_{E}}\lambda_{K,\mathbf{V}_{j}}\\
&  =-\frac{1}{15}\partial_{\mathbf{n}_{E}}\lambda_{K,\mathbf{V}_{3}}\left\vert
K\right\vert \overset{\text{(\ref{normcomp})}}{=}\frac{\left\vert E\right\vert
}{30},\\
\left\vert \int_{K}\partial_{\mathbf{n}_{E}}w_{E}\right\vert  &
\leq\left\vert K\right\vert ^{1/2}\left\Vert \nabla w_{E}\right\Vert
_{\mathbf{L}^{2}\left(  K\right)  }\leq C\left\vert K\right\vert ^{1/2}%
\sqrt{\log\left(  k+1\right)  }.
\end{align*}
In this way, $\left\vert \alpha_{K}\right\vert \leq C\sqrt{\log\left(
k+1\right)  }$ and an inverse inequality for quartic polynomials gives us%
\[
\left\Vert \nabla\psi_{E}\right\Vert _{\mathbf{L}^{2}\left(  K\right)  }\leq
Ch_{K}^{-1}\left\Vert \psi_{E}\right\Vert _{L^{2}\left(  K\right)  }\leq
Ch_{K}^{-1}\sqrt{\log\left(  k+1\right)  }\left\Vert 1\right\Vert
_{L^{2}\left(  K\right)  }\leq C\sqrt{\log\left(  k+1\right)  }.
\]
Hence, property (iv) follows.%
%TCIMACRO{\TeXButton{End Proof}{\endproof}}%
%BeginExpansion
\endproof
%EndExpansion

Next we recall a result which goes back to Vogelius \cite{vogelius1983right}
and Scott-Vogelius \cite{ScottVogelius}, see also \cite[Proof of Thm.
1]{GuzmanScott2019}.

\begin{definition}
\label{DefAz}Let $\eta_{0}$ be as in Lemma \ref{Lemangle}. For $0\leq\eta
<\eta_{0}$, the subspace $M_{\eta,k-1}^{\operatorname*{SV}}\left(
\mathcal{T}\right)  $ of the pressure space $M_{k-1}\left(  \mathcal{T}%
\right)  $ is given by%
\begin{equation}
M_{\eta,k-1}^{\operatorname*{SV}}\left(  \mathcal{T}\right)  :=\left\{  q\in
M_{k-1}\left(  \mathcal{T}\right)  \mid\forall\mathbf{z}\in\mathcal{C}%
_{\mathcal{T}}\left(  \eta\right)  :A_{\mathcal{T},\mathbf{z}}\left(
q\right)  =0\right\}  , \label{defMSV}%
\end{equation}
where, for $\mathbf{z}\in\mathcal{C}_{\mathcal{T}}\left(  \eta\right)  $, the
functional $A_{\mathcal{T},\mathbf{z}}\left(  q\right)  $ is as follows: fix
the counterclockwise numbering $K_{\ell}$, $1\leq\ell\leq m$, of the triangles
in the patch $\mathcal{T}_{\mathbf{z}}$ by the condition $K_{1}\cap
K_{2}=\mathfrak{E}\left(  \mathbf{z}\right)  $ and set%
\begin{equation}
A_{\mathcal{T},\mathbf{z}}\left(  q\right)  =\sum_{\ell=1}^{m}\left(
-1\right)  ^{\ell}\left(  \left.  q\right\vert _{K_{\ell}}\right)  \left(
\mathbf{z}\right)  . \label{DefATz}%
\end{equation}

\end{definition}

Note that $M_{0,k-1}^{\operatorname*{SV}}\left(  \mathcal{T}\right)  $ is the
pressure space introduced by Vogelius \cite{vogelius1983right} and
Scott-Vogelius \cite{ScottVogelius} and the following inclusions hold: for
$0\leq\eta\leq\eta^{\prime}\leq\eta_{0}$%
\[
M_{\eta^{\prime},k-1}^{\operatorname*{SV}}\left(  \mathcal{T}\right)  \subset
M_{\eta,k-1}^{\operatorname*{SV}}\left(  \mathcal{T}\right)  \subset
M_{0,k-1}^{\operatorname*{SV}}\left(  \mathcal{T}\right)  =Q_{h}^{k-1}%
\]
with the pressure space $Q_{h}^{k-1}$ in \cite[p. 517]{GuzmanScott2019}.

For the Scott-Vogelius pressure space $M_{0,k-1}^{\operatorname*{SV}}\left(
\mathcal{T}\right)  $, the existence of a continuous right-inverse of the
divergence operator into $\mathbf{S}_{k,0}\left(  \mathcal{T}\right)  $ was
proved in \cite{vogelius1983right} and \cite{ScottVogelius}.

\begin{proposition}
[Scott-Vogelius]For any $p\in M_{0,k-1}^{\operatorname*{SV}}\left(
\mathcal{T}\right)  $ there exists some $\mathbf{v}\in\mathbf{S}_{k,0}\left(
\mathcal{T}\right)  $ such that%
\[
\operatorname*{div}\mathbf{v}=q\mathbf{\quad}\text{and\quad}\left\Vert
\mathbf{v}\right\Vert _{\mathbf{H}^{1}\left(  \Omega\right)  }\leq C\left\Vert
q\right\Vert _{L^{2}\left(  \Omega\right)  },
\]
for a constant which only depends on the shape-regularity of the mesh, the
polynomial degree $k$, and on $\Theta_{\min}^{-1}$, where%
\begin{equation}
\Theta_{\min}:=\min_{\mathbf{z}\in\mathcal{V}\left(  \mathcal{T}\right)
\backslash\mathcal{C}_{\mathcal{T}}}\Theta\left(  \mathbf{z}\right)  .
\label{thetamin}%
\end{equation}
In particular, the constant $C$ is independent of $h$.
\end{proposition}

In Lemma \ref{LemProject}, we will show that, by subtracting the divergence of
a suitable Crouzeix-Raviart velocity from a given pressure in $M_{k-1}\left(
\mathcal{T}\right)  $, the resulting modified pressure belongs to the reduced
pressure space $M_{\eta,k-1}^{\operatorname*{SV}}\left(  \mathcal{T}\right)
$. As a preliminary, we need a bound of the functional $A_{\mathcal{T}%
,\mathbf{z}}$ in (\ref{DefATz}) which is explicit with respect to the local
mesh size and polynomial degree.

\begin{lemma}
\label{LemFunctional}There exists a constant $C$ which only depends on the
shape-regularity of the mesh such that%
\[
\left\vert A_{\mathcal{T},\mathbf{z}}\left(  q\right)  \right\vert \leq
C\frac{k^{2}}{h_{\mathbf{z}}}\left\Vert q\right\Vert _{L^{2}\left(
\omega_{\mathbf{z}}\right)  }\quad\forall q\in\mathbb{P}_{k-1}\left(
\mathcal{T}\right)
\]
for any $k\in\mathbb{N}$.
\end{lemma}

%

%TCIMACRO{\TeXButton{Proof}{\proof}}%
%BeginExpansion
\proof
%EndExpansion
Let $\mathbf{z}\in\mathcal{V}\left(  \mathcal{T}\right)  $ and $K\in
\mathcal{T}_{\mathbf{z}}$. The affine pullback to the reference triangle is
denoted by $\chi_{K}:\widehat{K}\rightarrow K$. For $q\in\mathbb{P}_{k}\left(
K\right)  $, let $\widehat{q}:=q\circ\chi_{K}$ and $\mathbf{\hat{z}}:=\chi
_{K}^{-1}\left(  \mathbf{z}\right)  $. Then%
\begin{equation}
\left\vert q\left(  \mathbf{z}\right)  \right\vert =\left\vert \widehat
{q}\left(  \widehat{\mathbf{z}}\right)  \right\vert \overset
{\text{\cite{Hesthaven_trace_2003}, \cite[Lem. 6.1]{Ainsworth_Jiang_hp}}}%
{\leq}\frac{\left(  k+1\right)  \left(  k+2\right)  }{\sqrt{2}}\left\Vert
\widehat{q}\right\Vert _{L^{2}\left(  \widehat{K}\right)  }=\binom{k+2}%
{2}\left\vert K\right\vert ^{-1/2}\left\Vert q\right\Vert _{L^{2}\left(
K\right)  }. \label{estfunctform}%
\end{equation}
A summation over all $K\in\mathcal{T}_{\mathbf{z}}$ leads to%
\[
\left\vert A_{\mathcal{T},\mathbf{z}}\left(  q\right)  \right\vert \leq
\binom{k+2}{2}\sum_{\ell=1}^{m}\left\vert K_{\ell}\right\vert ^{-1/2}%
\left\Vert q\right\Vert _{L^{2}\left(  K_{\ell}\right)  }\leq C\binom{k+2}%
{2}h_{\mathbf{z}}^{-1}\left\Vert q\right\Vert _{L^{2}\left(  \omega
_{\mathbf{z}}\right)  },
\]
where $C$ only depends on the shape-regularity of the mesh.%
%TCIMACRO{\TeXButton{End Proof}{\endproof}}%
%BeginExpansion
\endproof
%EndExpansion

The following lemma shows that the non-conforming Crouzeix-Raviart elements
allow us to modify a general pressure $q\in\mathbb{P}_{k-1,0}\left(
\mathcal{T}\right)  $ in such a way that the result belongs to $M_{\eta
,k-1}^{\operatorname*{SV}}\left(  \mathcal{T}\right)  $ \textit{provided
}$\mathcal{C}_{\mathcal{T}}^{\operatorname{acute}}=\emptyset$.

\begin{lemma}
\label{LemProject}Let assumption (\ref{Assumptab}) be satisfied and
$\mathcal{C}_{\mathcal{T}}^{\operatorname{acute}}=\emptyset$. There exists a
constant $\eta_{2}>0$ which only depends on the shape-regularity of the mesh
and $\alpha_{\Omega}$ (see (\ref{defalphatau})) such that for any fixed
$0\leq\eta<\eta_{2}$ and any $q\in\mathbb{P}_{k-1,0}\left(  \mathcal{T}%
\right)  $, there exists some $\mathbf{v}_{q}\in\mathbf{CR}_{k,0}\left(
\mathcal{T}\right)  $ such that%
\begin{align}
\int_{K}\operatorname{div}\mathbf{v}_{q}  &  =0\quad\forall K\in
\mathcal{T},\label{trianglemean}\\
q-\operatorname*{div}\nolimits_{\mathcal{T}}\mathbf{v}_{q}  &  \in
M_{\eta,k-1}^{\operatorname*{SV}}\left(  \mathcal{T}\right)
\label{inclusionMV}%
\end{align}
and%
\begin{equation}
\left\Vert \nabla_{\mathcal{T}}\mathbf{v}_{q}\right\Vert _{\mathbb{L}%
^{2}\left(  \Omega\right)  }\leq C_{\operatorname*{CR}}\sqrt{\log\left(
k+1\right)  }\left\Vert q\right\Vert _{L^{2}\left(  \Omega\right)  }.
\label{CRlifting1}%
\end{equation}
The constant $C_{\operatorname*{CR}}$ depends only on the shape-regularity of
the mesh and $\alpha_{\Omega}$.
\end{lemma}

%

%TCIMACRO{\TeXButton{Proof}{\proof}}%
%BeginExpansion
\proof
%EndExpansion
Let $q\in\mathbb{P}_{k-1,0}\left(  \mathcal{T}\right)  $. Let the fans
$\mathcal{C}_{\mathcal{T},\ell}\left(  \eta\right)  $, $\ell\in\mathcal{J}$,
be as in Definition \ref{DefFan}. For each fan $\mathcal{C}_{\mathcal{T},\ell
}\left(  \eta\right)  $ we employ an ansatz%
\begin{equation}
\mathbf{v}_{\ell}:=\sum_{\mathbf{z}\in\mathcal{C}_{\mathcal{T},\ell}\left(
\eta\right)  }\alpha_{\ell,\mathbf{z}}\tilde{B}_{k,\mathfrak{E}\left(
\mathbf{z}\right)  }^{\operatorname*{CR}}\mathfrak{N}\left(  \mathbf{z}%
\right)  \label{defvl}%
\end{equation}
for $\tilde{B}_{k,\mathfrak{E}\left(  \mathbf{z}\right)  }^{\operatorname*{CR}%
}$ as in (\ref{defboldpsi}), where the coefficients $\alpha_{\ell,\mathbf{z}%
}\in\mathbb{R}$ are defined next. The global function $\mathbf{v}_{q}$ is then
given by%
\[
\mathbf{v}_{q}=\sum_{\ell=1}^{N}\mathbf{v}_{\ell}.
\]
Property (\ref{trianglemean}) follows from this ansatz by using Lemma
\ref{Lemtrianglemean}. Next, we define the coefficients $\alpha_{\ell
,\mathbf{z}}$ in (\ref{defvl}) such that (\ref{inclusionMV}) holds and prove
the norm estimates for $\mathbf{v}_{q}$. Our construction of the fans implies
that open interiors of the supports of $\mathbf{v}_{\ell}$ are pairwise
disjoint (see Lem. \ref{Lemomegaell}); as a consequence the definition of
$\left(  \alpha_{\ell,\mathbf{z}}\right)  _{\mathbf{z}\in\mathcal{C}%
_{\mathcal{T}}\left(  \eta\right)  }$ and the estimate of $\nabla
_{\mathcal{T}}\mathbf{v}_{q}$ can be performed for each fan separately.

\begin{figure}
[ptb]
\begin{center}
\includegraphics[
height=2.2805in,
width=4.2004in
]%
{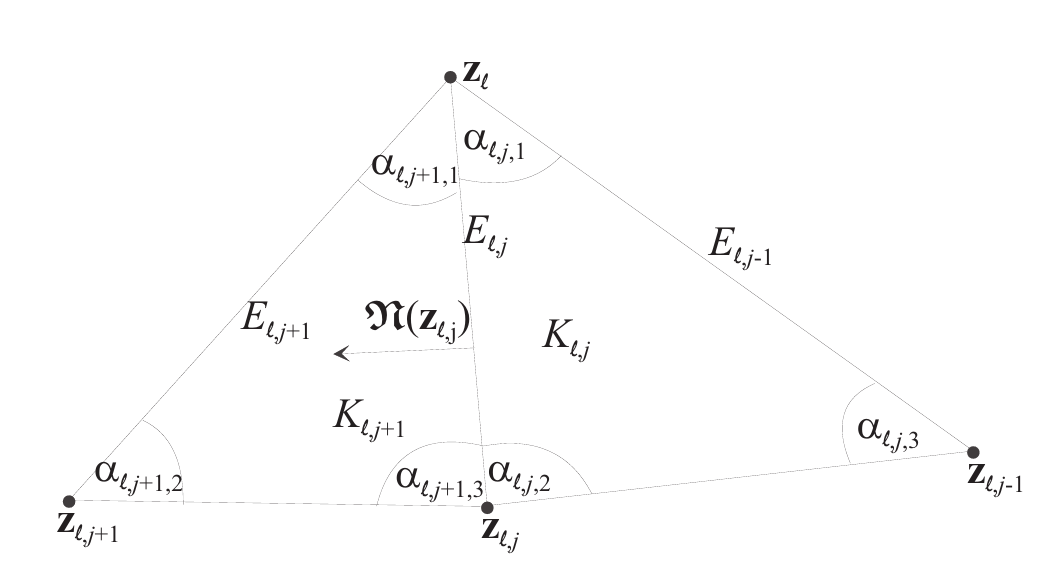}%
\caption{Local numbering convention of the angles in $K_{\ell,j}$ and
$K_{\ell,j+1}$. The angle in $K_{\ell,j}$ at $\mathbf{z}_{\ell}$ is denoted by
$\alpha_{\ell,j,1}$, at $\mathbf{z}_{\ell,j}$ by $\alpha_{\ell,j,2}$, at
$\mathbf{z}_{\ell,j-1}$ by $\alpha_{\ell,j,3}$ and in $K_{\ell,j+1}$
accordingly.}%
\label{Figfan2}%
\end{center}
\end{figure}

For $\mathbf{z}\in\mathcal{C}_{\mathcal{T},\ell}\left(  \eta\right)  $, let
$E:=\mathfrak{E}\left(  \mathbf{z}\right)  $, $\mathbf{n}_{E}:=\mathfrak{N}%
\left(  \mathbf{z}\right)  $. Let $K_{\mathbf{z}}^{-}$, $K_{\mathbf{z}}^{+}$,
denote the triangles in $\mathcal{T}_{E}$ with the convention that
$\mathbf{n}_{E}$ points into $K_{\mathbf{z}}^{+}$. The vertex in
$K_{\mathbf{z}}^{\pm}$ opposite to $E$ is denoted by $\mathbf{A}^{\pm}$. We
use%
\[
\operatorname*{div}\left(  \left.  \tilde{B}_{k,E}^{\operatorname*{CR}%
}\mathbf{n}_{E}\right\vert _{K}\right)  \left(  \mathbf{y}\right)
\overset{\text{(\ref{neinher})}}{=}\operatorname*{div}\left(  \left.
B_{k,E}^{\operatorname*{CR}}\mathbf{n}_{E}\right\vert _{K}\right)  \left(
\mathbf{y}\right)  \quad\forall\mathbf{y}\in\mathcal{V}\left(  K\right)
\]
and compute the divergence of $B_{k,E}^{\operatorname*{CR}}\mathbf{n}_{E}$%
\begin{equation}
\operatorname*{div}\left(  \left.  B_{k,E}^{\operatorname*{CR}}\mathbf{n}%
_{E}\right\vert _{K}\right)  =\left\{
\begin{array}
[c]{ll}%
\mp\frac{\left\vert E\right\vert }{\left\vert K\right\vert }L_{k}^{\prime
}\left(  1-2\lambda_{K,\mathbf{A}^{\pm}}\right)  & \text{on }K=K_{\mathbf{z}%
}^{\pm},\quad i=1,2,\\
0 & \text{otherwise.}%
\end{array}
\right.  \label{divrep}%
\end{equation}
Well-known properties of Legendre polynomials applied to (\ref{divrep}) imply
that for any vertex $\mathbf{y}$ of $K$ and odd polynomial degree $k$%
\begin{equation}
\operatorname*{div}\left(  \left.  \tilde{B}_{k,E}^{\operatorname*{CR}%
}\mathbf{n}_{E}\right\vert _{K}\right)  \left(  \mathbf{y}\right)  =\mp
\binom{k+1}{2}\times\left\{
\begin{array}
[c]{ll}%
\frac{\left\vert E\right\vert }{\left\vert K\right\vert } & \forall
\mathbf{y}\in\mathcal{V}\left(  K\right)  ,\quad\text{if }K=K_{\mathbf{z}%
}^{\pm}\\
0 & \text{otherwise.}%
\end{array}
\right.  \label{divcomp}%
\end{equation}
Hence the condition $A_{\mathcal{T},\mathbf{y}}\left(  q-\operatorname*{div}%
_{\mathcal{T}}\mathbf{v}_{\ell}\right)  =0$ for all $\mathbf{y}\in
\mathcal{C}_{\mathcal{T},\ell}\left(  \eta\right)  $ is equivalent to the
system of linear equation%
\begin{equation}
\mathbf{M}_{\ell}%
%TCIMACRO{\TeXButton{boldalpha}{\mbox{\boldmath$ \alpha$}}}%
%BeginExpansion
\mbox{\boldmath$ \alpha$}%
%EndExpansion
_{\ell}=\mathbf{r}_{\ell} \label{defalphaell}%
\end{equation}
with%
\begin{equation}
\mathbf{M}_{\ell}:=\left(  A_{\mathcal{T},\mathbf{y}}\left(
\operatorname*{div}\nolimits_{\mathcal{T}}\left(  B_{k,\mathfrak{E}\left(
\mathbf{z}\right)  }^{\operatorname*{CR}}\mathfrak{N}\left(  \mathbf{z}%
\right)  \right)  \right)  \right)  _{\substack{\mathbf{y}\in\mathcal{C}%
_{\mathcal{T},\ell}\left(  \eta\right)  \\\mathbf{z}\in\mathcal{C}%
_{\mathcal{T},\ell}\left(  \eta\right)  }},\quad%
%TCIMACRO{\TeXButton{boldalpha}{\mbox{\boldmath$ \alpha$}}}%
%BeginExpansion
\mbox{\boldmath$ \alpha$}%
%EndExpansion
_{\ell}:=\left(  \alpha_{\ell,\mathbf{z}}\right)  _{\mathbf{z}\in
\mathcal{C}_{\mathcal{T},\ell}\left(  \eta\right)  },\quad\mathbf{r}_{\ell
}:=\left(  A_{\mathcal{T},\mathbf{y}}\left(  q\right)  \right)  _{\mathbf{y}%
\in\mathcal{C}_{\mathcal{T},\ell}\left(  \eta\right)  }. \label{DefMl}%
\end{equation}
The matrix $\mathbf{M}_{\ell}$ is explicitly given by%
\[
\mathbf{M}_{\ell}:=-\binom{k+1}{2}\left[
\begin{array}
[c]{ccccc}%
\frac{\left\vert E_{\ell,1}\right\vert }{\left\vert K_{\ell,1}\right\vert
}+\frac{\left\vert E_{\ell,1}\right\vert }{\left\vert K_{\ell,2}\right\vert }
& \frac{\left\vert E_{\ell,2}\right\vert }{\left\vert K_{\ell,2}\right\vert }
& 0 & \ldots & 0\\
\frac{\left\vert E_{\ell,1}\right\vert }{\left\vert K_{\ell,2}\right\vert } &
\frac{\left\vert E_{\ell,2}\right\vert }{\left\vert K_{\ell,2}\right\vert
}+\frac{\left\vert E_{\ell,2}\right\vert }{\left\vert K_{\ell,3}\right\vert }
& \ddots & \ddots & \vdots\\
0 & \ddots & \ddots &  & 0\\
\vdots & \ddots &  &  & \frac{\left\vert E_{\ell,n_{\ell}}\right\vert
}{\left\vert K_{\ell,n_{\ell}}\right\vert }\\
0 & \ldots & 0 & \frac{\left\vert E_{\ell,n_{\ell}-1}\right\vert }{\left\vert
K_{\ell,n_{\ell}}\right\vert } & \frac{\left\vert E_{\ell,n_{\ell}}\right\vert
}{\left\vert K_{\ell,n_{\ell}}\right\vert }+\frac{\left\vert E_{\ell,n_{\ell}%
}\right\vert }{\left\vert K_{\ell,n_{\ell+1}}\right\vert }%
\end{array}
\right]  .
\]
We use (cf. Fig. \ref{Figfan2})%
%TCIMACRO{\FRAME{ftbpFU}{4.2004in}{2.2805in}{0pt}{\Qcb{Local numbering
%convention of the angles in $K_{\ell,j}$ and $K_{\ell,j+1}$. The angle in
%$K_{\ell,j}$ at $\QTR{bf}{z}_{\ell}$ is denoted by $\alpha_{\ell,j,1}$, at
%$\QTR{bf}{z}_{\ell,j}$ by $\alpha_{\ell,j,2}$, at $\QTR{bf}{z}_{\ell,j-1}$ by
%$\alpha_{\ell,j,3}$ and in $K_{\ell,j+1}$ accordingly.}}{\Qlb{Figfan2}%
%}{fan2.eps}{\special{ language "Scientific Word";  type "GRAPHIC";
%width 4.2004in;  height 2.2805in;  depth 0pt;  original-width 4.1502in;
%original-height 2.2407in;  cropleft "0";  croptop "1";  cropright "1";
%cropbottom "0";  filename 'fan2.eps';file-properties "XNPEU";}}}%
%BeginExpansion
%EndExpansion%
\[
\frac{\left\vert E_{\ell,j}\right\vert }{\left\vert K_{\ell,j}\right\vert
}=\frac{2\sin\left(  \alpha_{\ell,j,1}+\alpha_{\ell,j,2}\right)  }{\left\vert
E_{\ell,j}\right\vert \sin\alpha_{\ell,j,1}\sin\alpha_{\ell,j,2}}%
\quad\text{and\quad}\frac{\left\vert E_{\ell,j}\right\vert }{\left\vert
K_{\ell,j+1}\right\vert }=\frac{2\sin\left(  \alpha_{\ell,j+1,1}+\alpha
_{\ell,j+1,3}\right)  }{\left\vert E_{\ell,j}\right\vert \sin\alpha
_{\ell,j+1,1}\sin\alpha_{\ell,j+1,3}}%
\]
(cf. \cite[formula before (3.29)]{CCSS_CR_1}) and obtain
\begin{equation}
\mathbf{M}_{\ell}=\mathbf{D}_{\ell}\left(  \mathbf{T}_{\ell}+%
%TCIMACRO{\TeXButton{boldcapdelta}{\mbox{\boldmath$ \Delta$}}}%
%BeginExpansion
\mbox{\boldmath$ \Delta$}%
%EndExpansion
_{\ell}\right)  \label{DefMlTlDl}%
\end{equation}
with $\mathbf{D}_{\ell}=-k\left(  k+1\right)  \operatorname*{diag}\left[
\left\vert E_{\ell,j}\right\vert ^{-1}:1\leq j\leq n_{\ell}\right]  $ and%
\begin{align}
\mathbf{T}_{\ell}  &  :=\left[
\begin{array}
[c]{ccccc}%
\frac{\sin\left(  \alpha_{\ell,1,1}+\alpha_{\ell,2,1}\right)  }{\sin
\alpha_{\ell,1,1}\sin\alpha_{\ell,2,1}} & \frac{1}{\sin\alpha_{\ell,2,1}} &
0 & \ldots & 0\\
\frac{1}{\sin\alpha_{\ell,2,1}} & \frac{\sin\left(  \alpha_{\ell,2,1}%
+\alpha_{\ell,3,1}\right)  }{\sin\alpha_{\ell,2,1}\sin\alpha_{\ell,3,1}} &
\ddots & \ddots & \vdots\\
0 & \ddots & \ddots &  & 0\\
\vdots & \ddots &  &  & \frac{1}{\sin\alpha_{\ell,n_{\ell},1}}\\
0 & \ldots & 0 & \frac{1}{\sin\alpha_{\ell,n_{\ell},1}} & \frac{\sin\left(
\alpha_{\ell,n_{\ell},1}+\alpha_{\ell,n_{\ell}+1,1}\right)  }{\sin\alpha
_{\ell,n_{\ell},1}\sin\alpha_{\ell,n_{\ell}+1,1}}%
\end{array}
\right]  ,\label{defTl}\\%
%TCIMACRO{\TeXButton{boldcapdelta}{\mbox{\boldmath$ \Delta$}}}%
%BeginExpansion
\mbox{\boldmath$ \Delta$}%
%EndExpansion
_{\ell}  &  :=\operatorname*{diag}\left[  \frac{\sin\left(  \alpha_{\ell
,j,2}+\alpha_{\ell,j+1,3}\right)  }{\sin\alpha_{\ell,j,2}\sin\alpha
_{\ell,j+1,3}}:1\leq j\leq n_{\ell}\right]  .\nonumber
\end{align}

In Lemma \ref{Tlest}, we will prove that the matrix $\mathbf{T}_{\ell}+%
%TCIMACRO{\TeXButton{boldcapdelta}{\mbox{\boldmath$ \Delta$}}}%
%BeginExpansion
\mbox{\boldmath$ \Delta$}%
%EndExpansion
_{\ell}$ is invertible and the inverse is bounded by a constant independent of
$h_{\mathcal{T}}$ and $k$. Hence,%
\begin{equation}
\left\Vert
%TCIMACRO{\TeXButton{boldalpha}{\mbox{\boldmath$ \alpha$}}}%
%BeginExpansion
\mbox{\boldmath$ \alpha$}%
%EndExpansion
_{\ell}\right\Vert \leq\tilde{C}\frac{h_{\mathbf{z}_{\ell}}}{k\left(
k+1\right)  }\left\Vert \mathbf{r}_{\ell}\right\Vert . \label{estboldalpha}%
\end{equation}
Let $\mathcal{T}_{\ell}:=\left\{  K_{\ell,j}:1\leq j\leq n_{\ell}+1\right\}  $
and $D_{\ell}:=\operatorname*{dom}\left(  \mathcal{T}_{\ell}\right)  $. We
estimate the function $\mathbf{v}_{\ell}$ in (\ref{defvl}) by%
\begin{align*}
\left\Vert \nabla_{\mathcal{T}}\mathbf{v}_{\ell}\right\Vert _{\mathbb{L}%
^{2}\left(  D_{\ell}\right)  }  &  \leq\left(  \sum_{\mathbf{z}\in
\mathcal{C}_{\mathcal{T},\ell}\left(  \eta\right)  }\left\vert \alpha
_{\ell,\mathbf{z}}\right\vert ^{2}\left\Vert \nabla_{\mathcal{T}}\tilde
{B}_{k,\mathfrak{E}\left(  \mathbf{z}\right)  }^{\operatorname*{CR}%
}\mathfrak{N}\left(  \mathbf{z}\right)  \right\Vert _{\mathbb{L}^{2}\left(
D_{\ell}\right)  }^{2}\right)  ^{1/2}\\
&  \leq\max_{\mathbf{z}\in\mathcal{C}_{\mathcal{T},\ell}\left(  \eta\right)
}\left\Vert \nabla_{\mathcal{T}}\tilde{B}_{k,\mathfrak{E}\left(
\mathbf{z}\right)  }^{\operatorname*{CR}}\mathfrak{N}\left(  \mathbf{z}%
\right)  \right\Vert _{\mathbb{L}^{2}\left(  D_{\ell}\right)  }\left\Vert
%TCIMACRO{\TeXButton{boldalpha}{\mbox{\boldmath$ \alpha$}}}%
%BeginExpansion
\mbox{\boldmath$ \alpha$}%
%EndExpansion
_{\ell}\right\Vert \\
&  \overset{\text{(\ref{pwgrad})}}{\leq}Ch_{\mathbf{z}_{\ell}}\frac{\sqrt
{\log\left(  k+1\right)  }}{\left(  k+1\right)  ^{2}}\left\Vert \mathbf{r}%
_{\ell}\right\Vert .
\end{align*}
The constant $C$ only depends on the shape-regularity of the mesh. We use
Lemma \ref{LemFunctional} and conclude that%
\[
\left\Vert \nabla_{\mathcal{T}}\mathbf{v}_{\ell}\right\Vert _{\mathbb{L}%
^{2}\left(  D_{\ell}\right)  }\leq C\sqrt{\log\left(  k+1\right)  }\left\Vert
q\right\Vert _{L^{2}\left(  D_{\ell}\right)  }\text{.}%
\]
Since the interiors of the supports $D_{\ell}$ have pairwise empty
intersection the estimate%
\[
\left\Vert \nabla_{\mathcal{T}}\mathbf{v}_{q}\right\Vert _{\mathbb{L}%
^{2}\left(  \Omega\right)  }\leq C\sqrt{\log\left(  k+1\right)  }\left\Vert
q\right\Vert _{L^{2}\left(  \Omega\right)  }%
\]
follows.%
%TCIMACRO{\TeXButton{End Proof}{\endproof}}%
%BeginExpansion
\endproof
%EndExpansion

\begin{definition}
\label{DefPiCR}Let assumption (\ref{Assumptab}) be satisfied and $\eta_{2}>0$
as in Lemma \ref{LemProject}. Fix $\eta\in\left[  0,\eta_{2}\right[  $. For
$q\in\mathbb{P}_{k,0}\left(  \mathcal{T}\right)  $, the linear map $\Pi
_{k}^{\operatorname*{CR}}:\mathbb{P}_{k-1,0}\left(  \mathcal{T}\right)
\rightarrow\mathbf{CR}_{k,0}\left(  \mathcal{T}\right)  $ is given by%
\[
\Pi_{k}^{\operatorname*{CR}}q:=\sum_{\ell\in\mathcal{J}}\sum_{\mathbf{z}%
\in\mathcal{C}_{\mathcal{T},\ell}\left(  \eta\right)  }\alpha_{\ell
,\mathbf{z}}\tilde{B}_{k,\mathfrak{E}\left(  \mathbf{z}\right)  }%
^{\operatorname*{CR}}\mathfrak{N}\left(  \mathbf{z}\right)
\]
with $%
%TCIMACRO{\TeXButton{boldalpha}{\mbox{\boldmath$ \alpha$}}}%
%BeginExpansion
\mbox{\boldmath$ \alpha$}%
%EndExpansion
_{\ell}$ as in (\ref{defalphaell}).
\end{definition}

For the proof of the following lemma we recall the definitions of some linear
maps from the literature which are related to the right inverse of the
divergence operator acting on some polynomial spaces.

Bernardi and Raugel introduced in \cite[Lem. II.4]{Bernardi_Raugel_1985} a
linear mapping $\Pi^{\operatorname*{BR}}:M_{k-1}\left(  \mathcal{T}\right)
\rightarrow\mathbf{S}_{2,0}\left(  \mathcal{T}\right)  $ with the property:
for any $q\in M_{k-1}\left(  \mathcal{T}\right)  $, the function
$\Pi^{\operatorname*{BR}}q\in\mathbf{S}_{2,0}\left(  \mathcal{T}\right)  $
satisfies%
\begin{equation}
\int_{K}q=\int_{K}\operatorname*{div}\left(  \Pi^{\operatorname*{BR}}q\right)
\mathbf{\quad}\forall K\in\mathcal{T} \label{defPiBR}%
\end{equation}
and
\begin{equation}
\left\Vert \Pi^{\operatorname*{BR}}q\right\Vert _{\mathbf{H}^{1}\left(
\Omega\right)  }\leq C_{\operatorname*{BR}}\left\Vert q\right\Vert
_{L^{2}\left(  \Omega\right)  } \label{CBR}%
\end{equation}
for a constant $C_{\operatorname*{BR}}$ which is independent of the mesh width
and the polynomial degree.

Next we consider some right inverse of the divergence operator on the space%
\[
M_{k-1}^{\operatorname*{V}}\left(  \mathcal{T}\right)  :=\left\{  q\in
M_{k-1}\left(  \mathcal{T}\right)  \mid\left(
\begin{array}
[c]{ll}%
\int_{K}q=0 & \forall K\in\mathcal{T}\\
\left.  q\right\vert _{K}\left(  \mathbf{y}\right)  =0 & \forall\mathbf{y}%
\in\mathcal{V}\left(  K\right)
\end{array}
\right)  \right\}  .
\]
Let%
\[
\mathbf{S}_{k}^{\operatorname*{V}}\left(  \mathcal{T}\right)  :=\left\{
\mathbf{u}\in\mathbf{S}_{k}\left(  \mathcal{T}\right)  \mid\forall
K\in\mathcal{T}:\quad\left.  \mathbf{u}\right\vert _{\partial K}%
=\mathbf{0}\right\}  .
\]
There exists a linear operator $\Pi^{\operatorname*{V}}:M_{k-1}%
^{\operatorname*{V}}\left(  \mathcal{T}\right)  \rightarrow\mathbf{S}%
_{k}^{\operatorname*{V}}\left(  \mathcal{T}\right)  $ such that for all $q\in
M_{k-1}^{\operatorname*{V}}\left(  \mathcal{T}\right)  $%
\begin{align}
q  &  =\operatorname*{div}\Pi^{\operatorname*{V}}q,\nonumber\\
\left\Vert \Pi^{\operatorname*{V}}q\right\Vert _{\mathbf{H}^{1}\left(
\Omega\right)  }  &  \leq C_{\operatorname*{V}}\left\Vert q\right\Vert
_{L^{2}\left(  \Omega\right)  }, \label{optests}%
\end{align}
where the constant $C_{\operatorname*{V}}$ is independent of the mesh width
and the polynomial degree. Note that in the original paper \cite[Lem.
2.5]{vogelius1983right} by M. Vogelius, the right-hand side in the estimate
(\ref{optests}) contains an additional factor $k^{\beta_{\operatorname*{V}}}$
for some positive $\beta_{\operatorname*{V}}$ (independent of the mesh width).
In \cite[Thm. 3.4]{Ainsworth_parker_I}, the operator in \cite[Lem.
2.5]{vogelius1983right} is modified and the estimate in the form
(\ref{optests}) is proved for the modified operator.

Finally, we reconsider the linear operator $\Pi^{\operatorname*{GS}}%
:M_{\eta,k-1}^{\operatorname*{SV}}\left(  \mathcal{T}\right)  \rightarrow
\mathbf{S}_{4,0}\left(  \mathcal{T}\right)  $ introduced by Guzm\'{a}n and
Scott in \cite[Proof of Lem. 6 and Lem. 7]{GuzmanScott2019} with the property
that, for any $q\in M_{\eta,k-1}^{\operatorname*{SV}}\left(  \mathcal{T}%
\right)  $, it holds%
%TCIMACRO{\TeXButton{Thetamin_eta}{\begin{subequations}
%\label{Thetamin_eta}
%\end{subequations}}}%
%BeginExpansion
\begin{subequations}
\label{Thetamin_eta}
\end{subequations}%
%EndExpansion%
\begin{align}
\left(  I-\operatorname*{div}\Pi^{\operatorname*{GS}}\right)  q  &  \in
M_{k-1}^{\operatorname*{V}}\left(  \mathcal{T}\right)  ,\tag{%
%TCIMACRO{\TeXButton{Thetamin_{e}ta}{\ref{Thetamin_eta}}}%
%BeginExpansion
\ref{Thetamin_eta}%
%EndExpansion
a}\label{Thetamin_etaa}\\
\left\Vert \nabla\Pi^{\operatorname*{GS}}q\right\Vert _{\mathbb{L}^{2}\left(
\Omega\right)  }  &  \leq C_{\operatorname*{GS}}k^{\kappa}\left(  \theta
_{\min}+\eta\right)  ^{-1}\left\Vert q\right\Vert _{L^{2}\left(
\Omega\right)  }\quad\text{for }\kappa=2. \tag{%
%TCIMACRO{\TeXButton{Thetamin_{e}ta}{\ref{Thetamin_eta}}}%
%BeginExpansion
\ref{Thetamin_eta}%
%EndExpansion
b}\label{Thetamin_etac}%
\end{align}

We emphasize that in \cite[Lemma 7]{GuzmanScott2019} the constant
$\Theta_{\min}^{-1}$ (cf. (\ref{thetamin})) instead of $\left(  \Theta_{\min
}+\eta\right)  ^{-1}$ appears in (\ref{Thetamin_etac}) so that the estimate of
$\left\Vert \nabla\Pi^{\operatorname*{GS}}q\right\Vert _{\mathbb{L}^{2}\left(
\Omega\right)  }$ for $q\in M_{0,k-1}^{\operatorname*{SV}}\left(
\mathcal{T}\right)  $ deteriorates in cases where the $\mathbf{z}$ is a
\textit{nearly critical} point, i.e., very close to the geometric situations
described in Remark \ref{RemCritGeom}. The proof of \cite[Lemma 7]%
{GuzmanScott2019} is split into an estimate related to points $\mathbf{z}$
with $A_{\mathcal{T},\mathbf{z}}\left(  q\right)  =0$ (cf. (\ref{DefATz})) and
an estimate for the remaining points $\mathbf{z}$ with $A_{\mathcal{T}%
,\mathbf{z}}\left(  q\right)  \neq0$. Only in this second part, the constant
$\Theta_{\min}^{-1}$ is involved. The result has been improved in \cite[Lem.
4.5]{Sauter_eta_wired} and it was shown that there is an operator $\Pi
_{\eta,k-1}:M_{\eta,k-1}^{\operatorname*{SV}}\left(  \mathcal{T}\right)
\rightarrow\mathbf{S}_{k,0}\left(  \mathcal{T}\right)  $ such that the
properties in (\ref{Thetamin_eta}) hold for $\kappa=0$: for any $q\in
M_{\eta,k-1}^{\operatorname*{SV}}\left(  \mathcal{T}\right)  $ it holds%
\begin{align*}
\left(  I-\operatorname*{div}\Pi_{\eta,k-1}\right)  q  &  \in M_{k-1}%
^{\operatorname*{V}}\left(  \mathcal{T}\right)  ,\\
\left\Vert \nabla\Pi_{\eta,k-1}q\right\Vert _{\mathbb{L}^{2}\left(
\Omega\right)  }  &  \leq C_{\pi}\left(  \Theta_{\min}+\eta\right)
^{-1}\left\Vert q\right\Vert _{L^{2}\left(  \Omega\right)  }.
\end{align*}

From Lemma \ref{LemProject}, we conclude that $q-\operatorname*{div}%
_{\mathcal{T}}\left(  \Pi_{k}^{\operatorname*{CR}}q\right)  \in M_{\eta
,k-1}^{\operatorname*{SV}}\left(  \mathcal{T}\right)  $ and the second part of
the proof in \cite[Lemma 7]{GuzmanScott2019} is applied only to points with%
\[
\min_{\mathbf{z}\in\mathcal{V}\left(  \mathcal{T}\right)  \backslash
\mathcal{C}_{\mathcal{T}}\left(  \eta\right)  }\Theta\left(  \mathbf{z}%
\right)  \geq\max\left\{  \eta,\Theta_{\min}\right\}  .
\]
Hence, (\ref{Thetamin_etac}) follows for $\eta$ depending only on the
shape-regularity of the mesh.

\begin{lemma}
\label{LemmaFinOdd}Let assumption (\ref{Assumptab}) be satisfied and let
$\mathcal{C}_{\mathcal{T}}^{\operatorname{acute}}=\emptyset$. There exists a
constant $\eta_{2}>0$ which only depends on the shape-regularity of the mesh
and $\alpha_{\Omega}$ as in (\ref{defalphatau}) such that for any fixed
$0\leq\eta<\eta_{2}$ and any $q\in\mathbb{P}_{k-1,0}\left(  \mathcal{T}%
\right)  $, there exists some $\mathbf{w}_{q}\in\mathbf{CR}_{k,0}\left(
\mathcal{T}\right)  $ such that%
\[
q=\operatorname*{div}\nolimits_{\mathcal{T}}\mathbf{w}_{q}%
\]
and%
\[
\left\Vert \mathbf{w}_{q}\right\Vert _{\mathbf{H}^{1}\left(  \mathcal{T}%
\right)  }\leq C\sqrt{\log\left(  k+1\right)  }\left(  \Theta_{\min}%
+\eta\right)  ^{-1}\left\Vert q\right\Vert _{L^{2}\left(  \Omega\right)  }.
\]
The constant $C$ only depends on the shape-regularity of the mesh and
$\alpha_{\Omega}$ but is independent of the mesh width and the polynomial
degree $k$.
\end{lemma}

%

%TCIMACRO{\TeXButton{Proof}{\proof}}%
%BeginExpansion
\proof
%EndExpansion
For the construction of $\mathbf{w}_{q}$ we follow and modify the lines of
proof in \cite[Thm. 1]{GuzmanScott2019} by a) involving the operator $\Pi
_{k}^{\operatorname*{CR}}$ and b) employing the concept of $\eta$-critical points.

For given $q\in\mathbb{P}_{k-1,0}\left(  \mathcal{T}\right)  $, we employ the
operators $\Pi_{k}^{\operatorname*{CR}}$, $\Pi^{\operatorname*{BR}}$,
$\Pi_{\eta,k-1}$, $\Pi^{\operatorname*{V}}$ in the definition of the function
$\mathbf{w}_{q}$
\begin{align}
\mathbf{w}_{q}  &  =\mathbf{T}_{1}+\mathbf{T}_{2}+\mathbf{T}_{3}%
+\mathbf{T}_{4},\label{defwq}\\
\mathbf{T}_{1}  &  :=\Pi^{\operatorname*{BR}}q,\nonumber\\
\mathbf{T}_{2}  &  :=\Pi_{k}^{\operatorname*{CR}}\left(  I-\operatorname*{div}%
\Pi^{\operatorname*{BR}}\right)  q,\nonumber\\
\mathbf{T}_{3}  &  :=\Pi_{\eta,k-1}\left(  I-\operatorname*{div}%
\nolimits_{\mathcal{T}}\Pi_{k}^{\operatorname*{CR}}\right)  \left(
I-\operatorname*{div}\Pi^{\operatorname*{BR}}\right)  q,\nonumber\\
\mathbf{T}_{4}  &  :=\Pi^{\operatorname*{V}}\left(  I-\operatorname*{div}%
\Pi_{\eta,k-1}\right)  \left(  I-\operatorname*{div}\nolimits_{\mathcal{T}}%
\Pi_{k}^{\operatorname*{CR}}\right)  \left(  I-\operatorname*{div}%
\Pi^{\operatorname*{BR}}\right)  q.\nonumber
\end{align}
By construction we have%
\[
\operatorname{div}_{\mathcal{T}}\mathbf{w}_{q}=q.
\]
The first two summands in (\ref{defwq}) satisfy%
\begin{align}
\left\Vert \mathbf{T}_{1}\right\Vert _{\mathbf{H}^{1}\left(  \Omega\right)  }
&  \overset{\text{(\ref{CBR})}}{\leq}C_{\operatorname*{BR}}\left\Vert
q\right\Vert _{L^{2}\left(  \Omega\right)  },\label{estPi1}\\
\left\Vert \mathbf{T}_{2}\right\Vert _{\mathbf{H}^{1}\left(  \mathcal{T}%
\right)  }  &  \overset{\text{(\ref{CRlifting1})}}{\leq}C_{\operatorname*{CR}%
}\sqrt{\log\left(  k+1\right)  }\left(  \left\Vert q\right\Vert _{L^{2}\left(
\Omega\right)  }+\left\Vert \Pi^{\operatorname*{BR}}q\right\Vert
_{\mathbf{H}^{1}\left(  \Omega\right)  }\right) \label{estPi2}\\
&  \overset{\text{(\ref{estPi1})}}{\leq}C_{\operatorname*{CR}}\sqrt
{\log\left(  k+1\right)  }\left(  1+C_{\operatorname*{BR}}\right)  \left\Vert
q\right\Vert _{L^{2}\left(  \Omega\right)  }.\nonumber
\end{align}
For the third term in (\ref{defwq}) we get%
\begin{align}
\left\Vert \nabla\mathbf{T}_{3}\right\Vert _{\mathbb{L}^{2}\left(
\Omega\right)  }  &  \leq C_{\pi}\left(  \Theta_{\min}+\eta\right)
^{-1}\left\Vert \left(  I-\operatorname*{div}\nolimits_{\mathcal{T}}\Pi
_{k}^{\operatorname*{CR}}\right)  \left(  I-\operatorname*{div}\Pi
^{\operatorname*{BR}}\right)  q\right\Vert _{L^{2}\left(  \Omega\right)
}\nonumber\\
&  \leq C_{\pi}\left(  \Theta_{\min}+\eta\right)  ^{-1}\left(  \left\Vert
\left(  I-\operatorname*{div}\Pi^{\operatorname*{BR}}\right)  q\right\Vert
_{L^{2}\left(  \Omega\right)  }+\left\Vert \operatorname*{div}%
\nolimits_{\mathcal{T}}\Pi_{k}^{\operatorname*{CR}}\left(
I-\operatorname*{div}\Pi^{\operatorname*{BR}}\right)  q\right\Vert
_{L^{2}\left(  \Omega\right)  }\right)  . \label{CGSest2}%
\end{align}
The combination with (\ref{CBR}), (\ref{CRlifting1}) leads to%
\begin{equation}
\left\Vert \nabla\mathbf{T}_{3}\right\Vert _{\mathbb{L}^{2}\left(
\Omega\right)  }\leq CC_{\pi}\left(  1+C_{\operatorname*{BR}}\right)  \left(
1+C_{\operatorname*{CR}}\right)  \frac{\sqrt{\log\left(  k+1\right)  }}%
{\Theta_{\min}+\eta}\left\Vert q\right\Vert _{L^{2}\left(  \Omega\right)  }.
\label{estPi3}%
\end{equation}
For the fourth term we get in a similar way%
\begin{align}
\left\Vert \nabla\mathbf{T}_{4}\right\Vert _{\mathbb{L}^{2}\left(
\Omega\right)  }  &  \overset{\text{(\ref{optests})}}{\leq}%
C_{\operatorname*{V}}\left(  \left\Vert \left(  I-\operatorname*{div}%
\nolimits_{\mathcal{T}}\Pi_{k}^{\operatorname*{CR}}\right)  \left(
I-\operatorname*{div}\Pi^{\operatorname*{BR}}\right)  q\right\Vert
_{L^{2}\left(  \Omega\right)  }+\left\Vert \operatorname*{div}\mathbf{T}%
_{3}\right\Vert _{L^{2}\left(  \Omega\right)  }\right) \nonumber\\
&  \leq C_{\operatorname*{V}}\left(  \left\Vert q\right\Vert _{L^{2}\left(
\Omega\right)  }+\left\Vert \operatorname{div}\mathbf{T}_{1}\right\Vert
_{L^{2}\left(  \Omega\right)  }+\left\Vert \operatorname{div}_{\mathcal{T}%
}\mathbf{T}_{2}\right\Vert _{L^{2}\left(  \Omega\right)  }+\left\Vert
\operatorname*{div}\mathbf{T}_{3}\right\Vert _{L^{2}\left(  \Omega\right)
}\right) \nonumber\\
&  \overset{\text{(\ref{estPi1}), (\ref{estPi2}), (\ref{estPi3})}}{\leq
}C_{\operatorname*{V}}\left(  1+C_{\operatorname*{BR}}\right)  \left(
2+C_{\operatorname*{CR}}\right)  \sqrt{\log\left(  k+1\right)  }\left(
1+CC_{\pi}\left(  \Theta_{\min}+\eta\right)  ^{-1}\right)  \left\Vert
q\right\Vert _{L^{2}\left(  \Omega\right)  }. \label{estPi4}%
\end{align}
The combination of (\ref{estPi1}), (\ref{estPi2}), (\ref{estPi3}),
(\ref{estPi4}) with (\ref{defwq}) leads to the assertion.
%TCIMACRO{\TeXButton{End Proof}{\endproof}}%
%BeginExpansion
\endproof
%EndExpansion

Lemma \ref{LemmaFinOdd} implies that for conforming triangulations
$\mathcal{T}$ which satisfy (\ref{Assumptab}) and $\mathcal{C}_{\mathcal{T}%
}^{\operatorname{acute}}=\emptyset$, there exists a bounded linear operator%
\[
\Pi_{\mathcal{T},k}^{\operatorname*{inv}}:\mathbb{P}_{k-1,0}\left(
\mathcal{T}\right)  \rightarrow\mathbf{CR}_{k,0}\left(  \mathcal{T}\right)
\]
such that $\operatorname*{div}_{\mathcal{T}}\circ\Pi_{\mathcal{T}%
,k}^{\operatorname*{inv}}$ is the identity on $\mathbb{P}_{k-1,0}\left(
\mathcal{T}\right)  $ and%
\[
\left\Vert \Pi_{\mathcal{T},k}^{\operatorname*{inv}}q\right\Vert
_{\mathbf{H}^{1}\left(  \mathcal{T}\right)  }\leq C_{\operatorname*{inv}}%
\sqrt{\log\left(  k+1\right)  }\left\Vert q\right\Vert _{L^{2}\left(
\Omega\right)  }%
\]
for a constant $C_{\operatorname*{inv}}$ which only depends on the
shape-regularity of the mesh and $\alpha_{\Omega}$ (cf. (\ref{defalphatau})).

\subsubsection{The case $\mathcal{C}_{\mathcal{T}}^{\operatorname*{acute}%
}\left(  \eta\right)  \neq\emptyset$\label{SecAcute}}

In this section, we remove the condition $\mathcal{C}_{\mathcal{T}%
}^{\operatorname{acute}}=\emptyset$ and construct a bounded right-inverse of
the piecewise divergence operator for odd $k\geq5$ and conforming
triangulations which contain at least one inner point. The construction is
based on the step-by-step procedure (cf. (\ref{stepbystepextension})) from the
triangulation $\mathcal{T}_{1}$ to $\mathcal{T}$. Inductively, we assume that
there is a triangulation $\mathcal{T}_{j}$ along a bounded right-inverse
$\Pi_{j,k}^{\operatorname*{inv}}:\mathbb{P}_{k-1,0}\left(  \mathcal{T}%
_{j}\right)  \rightarrow\mathbf{CR}_{k,0}\left(  \mathcal{T}_{j}\right)  $ of
the piecewise divergence operator. A single extension step is analysed by the
following lemma.

\begin{lemma}
Let $\mathcal{T}$ denote a conforming triangulation for the domain
$\Omega:=\operatorname*{dom}\mathcal{T}$ and let $\mathcal{T}^{\prime}%
\subset\mathcal{T}$ be a subset such that every triangle $K\in\mathcal{T}%
\backslash\mathcal{T}^{\prime}$ has one edge, say $E$, which belongs to
$\mathcal{E}\left(  \mathcal{T}^{\prime}\right)  $. We assume that
$\mathcal{T}^{\prime}$ has at least one inner vertex and set $\Omega^{\prime
}:=\operatorname*{dom}\mathcal{T}^{\prime}$. Assume that there exists a
bounded linear operator $\Pi_{\mathcal{T}^{\prime},k}^{\operatorname*{inv}%
}:\mathbb{P}_{k-1,0}\left(  \mathcal{T}^{\prime}\right)  \rightarrow
\mathbf{CR}_{k,0}\left(  \mathcal{T}^{\prime}\right)  $ with
$\operatorname*{div}_{\mathcal{T}^{\prime}}\circ\Pi_{\mathcal{T}^{\prime}%
,k}^{\operatorname*{inv}}=\operatorname*{Id}$ on $\mathbb{P}_{k-1,0}\left(
\mathcal{T}^{\prime}\right)  $ and%
\[
\left\Vert \Pi_{\mathcal{T}^{\prime},k}^{\operatorname*{inv}}q\right\Vert
_{\mathbf{H}^{1}\left(  \Omega^{\prime}\right)  }\leq C_{\mathcal{T}^{\prime}%
}\left\Vert q\right\Vert _{L^{2}\left(  \Omega^{\prime}\right)  }.
\]
Then, there exists a linear operator $\Pi_{\mathcal{T},k}^{\operatorname*{inv}%
}:\mathbb{P}_{k-1,0}\left(  \mathcal{T}\right)  \rightarrow\mathbf{CR}%
_{k,0}\left(  \mathcal{T}\right)  $ with $\operatorname*{div}_{\mathcal{T}%
}\circ\Pi_{\mathcal{T},k}^{\operatorname*{inv}}=\operatorname*{Id}$ on
$\mathbb{P}_{k-1,0}\left(  \mathcal{T}\right)  $ and%
\[
\left\Vert \Pi_{\mathcal{T},k}^{\operatorname*{inv}}q\right\Vert
_{\mathbf{H}^{1}\left(  \Omega\right)  }\leq C_{\mathcal{T}}\left\Vert
q\right\Vert _{L^{2}\left(  \Omega\right)  }\quad\text{with\quad
}C_{\mathcal{T}}:=C_{3}\sqrt{\log\left(  k+1\right)  }C_{\mathcal{T}^{\prime}}%
\]
for a constant $C_{3}$ which depends only on the shape-regularity of the mesh
and on $\alpha_{\Omega}$.
\end{lemma}

%

%TCIMACRO{\TeXButton{Proof}{\proof}}%
%BeginExpansion
\proof
%EndExpansion
Let $q\in\mathbb{P}_{k-1,0}\left(  \mathcal{T}\right)  $. We set
$\mathbf{v}_{0}:=\Pi^{\operatorname*{BR}}q$ where the operator $\Pi
^{\operatorname*{BR}}$ is as in (\ref{defPiBR}) and satisfies%
\[
\left\Vert \mathbf{v}_{0}\right\Vert _{\mathbf{H}^{1}\left(  \Omega\right)
}\leq C_{\operatorname*{BR}}\left\Vert q\right\Vert _{L^{2}\left(
\Omega\right)  }.
\]
Hence, $q_{1}:=q-\operatorname*{div}\mathbf{v}_{0}$ belongs to $\mathbb{P}%
_{k-1,0}\left(  \mathcal{T}\right)  $ and has trianglewise integral mean zero.

\begin{figure}
[ptb]
\begin{center}
\includegraphics[
height=2.84in,
width=5.7276in
]%
{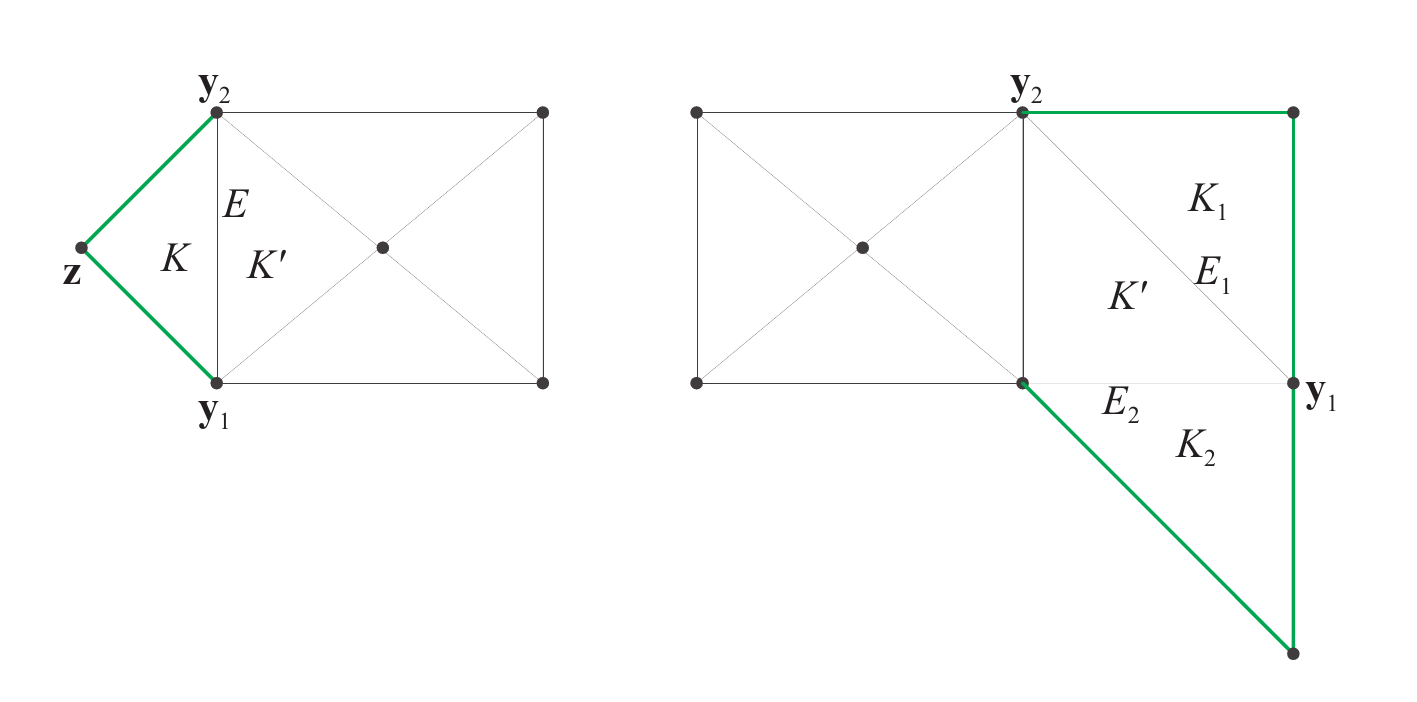}%
\caption{The black triangles form the triangulation $\mathcal{T}^{\prime}$.
Left: One triangle $K$ is attached to $\mathcal{T}^{\prime}$ having a common
side $E$ with $K^{\prime}\in\mathcal{T}^{\prime}$ and $\mathcal{T}%
_{\operatorname*{out}}\left(  K^{\prime}\right)  =\left\{  K\right\}  $.
Right: Two triangles $K_{1}$, $K_{2}\notin\mathcal{T}^{\prime}$ are attached
to a triangle $K^{\prime}\in\mathcal{T}^{\prime}$ and $\mathcal{T}%
_{\operatorname*{out}}\left(  K^{\prime}\right)  =\left\{  K_{1}%
,K_{2}\right\}  $.}%
\label{Fig_glue1}%
\end{center}
\end{figure}

The following construction is illustrated in Figure \ref{Fig_glue1}.%
%TCIMACRO{\FRAME{ftbpFU}{5.7276in}{2.84in}{0pt}{\Qcb{The black triangles form
%the triangulation $\QTR{cal}{T}^{\prime}$. Left: One triangle $K$ is attached
%to $\QTR{cal}{T}^{\prime}$ having a common side $E$ with $K^{\prime}%
%\in\QTR{cal}{T}^{\prime}$ and $\QTR{cal}{T}_{\operatorname*{out}}\left(
%K^{\prime}\right)  =\left\{  K\right\}  $. Right: Two triangles $K_{1}$,
%$K_{2}\notin\QTR{cal}{T}^{\prime}$ are attached to a triangle $K^{\prime}%
%\in\QTR{cal}{T}^{\prime}$ and $\QTR{cal}{T}_{\operatorname*{out}}\left(
%K^{\prime}\right)  =\left\{  K_{1},K_{2}\right\}  $.}}{\Qlb{Fig_glue1}%
%}{fan_glue1.eps}{\special{ language "Scientific Word";  type "GRAPHIC";
%width 5.7276in;  height 2.84in;  depth 0pt;  original-width 5.6688in;
%original-height 2.7968in;  cropleft "0";  croptop "1";  cropright "1";
%cropbottom "0";  filename 'fan_glue1.EPS';file-properties "XNPEU";}}}%
%BeginExpansion
%EndExpansion
Let $K^{\prime}\in\mathcal{T}^{\prime}$ be such that there exists a non-empty
subset $\mathcal{T}_{\operatorname*{out}}\left(  K^{\prime}\right)
\subset\mathcal{T}\backslash\mathcal{T}^{\prime}$ having the property that any
$K\in\mathcal{T}_{\operatorname*{out}}\left(  K^{\prime}\right)  $ shares an
edge with $K^{\prime}$. We have $\left\vert \mathcal{T}_{\operatorname*{out}%
}\left(  K^{\prime}\right)  \right\vert \leq2$; indeed, if $\left\vert
\mathcal{T}_{\operatorname*{out}}\left(  K^{\prime}\right)  \right\vert =3$,
then all three edges of $K^{\prime}$ are boundary edges which implies
$\mathcal{T}^{\prime}=\left\{  K^{\prime}\right\}  $ and violates the
condition that $\mathcal{T}^{\prime}$ must contain an inner vertex.

For $K\in\mathcal{T}_{\operatorname*{out}}\left(  K^{\prime}\right)  $, let
$\mathbf{z}$ denote the vertex in $K$ opposite to $E$ and set $\omega
_{E}=K\cup K^{\prime}$. The endpoints of $E$ are denoted by $\mathbf{y}_{1}$,
$\mathbf{y}_{2}$. We employ the ansatz (cf. (\ref{defboldpsi}))%
\begin{equation}
\mathbf{v}_{1}:=\alpha\tilde{B}_{k,E}^{\operatorname*{CR}}\mathbf{n}_{E}
\label{ansatzv1}%
\end{equation}
with the convention that $\mathbf{n}_{E}$ is the unit vector orthogonal to $E$
and directed into $K^{\prime}$. By construction it holds $\mathbf{v}_{1}%
\in\mathbf{CR}_{k,0}\left(  \mathcal{T}\right)  $ and $\operatorname*{supp}%
\mathbf{v}_{1}\subset\omega_{E}$. We determine $\alpha$ in (\ref{ansatzv1})
such that $\operatorname*{div}\left(  \left.  \mathbf{v}_{1}\right\vert
_{K}\right)  \left(  \mathbf{z}\right)  =q_{1}\left(  \mathbf{z}\right)  $ and
employ (\ref{divcomp}) to get%
\[
\left\vert \operatorname{div}\left(  \tilde{B}_{k,E}^{\operatorname*{CR}%
}\mathbf{n}_{E}\right)  \left(  \mathbf{z}\right)  \right\vert =\binom{k+1}%
{2}\frac{\left\vert E\right\vert }{\left\vert K\right\vert }.
\]
Hence $\left\vert \alpha\right\vert =\left\vert q_{1}\left(  \mathbf{z}%
\right)  \right\vert \frac{\left\vert K\right\vert }{\left\vert E\right\vert
}/\binom{k+1}{2}$ and we conclude as in the proof of Lemma \ref{LemProject}
that%
\[
\left\Vert \nabla_{\mathcal{T}}\mathbf{v}_{1}\right\Vert _{\mathbb{L}%
^{2}\left(  \omega_{E}\right)  }\leq C_{\operatorname*{CR}}\sqrt{\log\left(
k+1\right)  }\left\Vert q_{1}\right\Vert _{L^{2}\left(  K\right)  }.
\]
We set%
\[
q_{2}:=q_{1}-\operatorname*{div}\nolimits_{\mathcal{T}}\mathbf{v}_{1}%
\quad\text{so that }q_{1}=\operatorname*{div}\nolimits_{\mathcal{T}}%
\mathbf{v}_{1}+q_{2}%
\]
and note that%
\begin{equation}
\left\Vert q_{2}\right\Vert _{L^{2}\left(  K\right)  }\leq\left(
1+C_{\operatorname*{CR}}\sqrt{\log\left(  k+1\right)  }\right)  \left\Vert
q_{1}\right\Vert _{L^{2}\left(  K\right)  }. \label{q2est}%
\end{equation}
The construction implies $q_{2}\in\mathbb{P}_{k-1,0}\left(  \mathcal{T}%
\right)  $, $q_{2}$ has trianglewise integral mean zero, and $q_{2}\left(
\mathbf{z}\right)  =0$. Next, we employ the vector field defined in \cite[Lem.
4.9]{Sauter_eta_wired} which is a modification of the cubic vector field
defined in \cite[(3.5)]{GuzmanScott2019} but allows for better $k$-explicit
estimates. We recall the relevant lemma from \cite{Sauter_eta_wired} for the
existence of such vector fields and collect important properties.

\begin{lemma}
[{\cite[Lem. 4.9]{Sauter_eta_wired}}]\label{LemFundVF}Let $\mathcal{T}$ be a
conforming triangulation of $\Omega$ and let $k\geq3$. Let $E\in
\mathcal{E}\left(  \mathcal{T}\right)  $ with endpoints $\mathbf{y}_{1}$,
$\mathbf{y}_{2}$. Then there exist vector fields $\mathbf{v}_{E,j}$,
$j\in\left\{  1,2\right\}  $, with the following properties%
\begin{equation}%
\begin{array}
[c]{ll}%
\mathbf{v}_{E,j}\in\mathbf{S}_{k}\left(  \mathcal{T}\right)  \text{,} &
\operatorname*{supp}\mathbf{v}_{E,j}\subset\omega_{E},\\
\int_{K}\operatorname*{div}\mathbf{v}_{E,j}=0 & \forall K\in\mathcal{T},\quad
j\in\left\{  1,2\right\} \\
\left(  \operatorname{div}\left.  \mathbf{v}_{E,j}\right\vert _{K}\right)
\left(  \mathbf{v}\right)  =\left\{
\begin{array}
[c]{cl}%
1 & \text{if }K\in\mathcal{T}\left(  E\right)  \wedge\mathbf{v}=\mathbf{v}%
_{j}\\
0 & \text{otherwise,}%
\end{array}
\right.  & \forall K\in\mathcal{T},\quad\forall\mathbf{v}\in\mathcal{V}\left(
K\right)  ,\quad\forall j\in\left\{  1,2\right\} \\
\left\Vert \nabla\mathbf{v}_{E,j}\right\Vert _{\mathbb{L}^{2}\left(
\omega_{E}\right)  }\leq Ch_{E}k^{-2}. &
\end{array}
\label{vEj}%
\end{equation}

\end{lemma}

We employ this vector field to the edge $E=K\cap K^{\prime}$, set%
\[
\mathbf{v}_{2}:=\sum_{j=1}^{2}\left.  q_{2}\right\vert _{K}\left(
\mathbf{y}_{j}\right)  \mathbf{v}_{E,j}%
\]
and define%
\begin{equation}
q_{3}:=q_{2}-\operatorname{div}\mathbf{v}_{2}\quad\text{so that }%
q_{2}=\operatorname{div}\mathbf{v}_{2}+q_{3}. \label{defq3}%
\end{equation}
The function $q_{3}\in\mathbb{P}_{k-1,0}\left(  \mathcal{T}\right)  $ has
trianglewise integral mean zero and $\left.  q_{3}\right\vert _{K}$ vanishes
in all vertices of $K$. The norm $\left\Vert \nabla\mathbf{v}_{2}\right\Vert
_{\mathbb{L}^{2}\left(  K^{\prime}\cup K\right)  }$ can be estimated in the
same way as the function $\mathbf{T}_{3}$ in (\ref{estPi3}); however the
factor $\left(  \Theta_{\min}+\eta\right)  ^{-1}$ does not appear as in
(\ref{estPi3}) since the last estimate in (\ref{vEj}) does not depend on these
quantities. In this way, we get%
\begin{align}
\left\Vert \nabla_{\mathcal{T}}\mathbf{v}_{2}\right\Vert _{\mathbb{L}%
^{2}\left(  \omega_{E}\right)  }  &  \leq C_{1}\left\Vert q_{2}\right\Vert
_{L^{2}\left(  K\right)  }\label{v2est1}\\
&  \overset{\text{(\ref{q2est})}}{\leq}C_{1}\left(  1+C_{\operatorname*{CR}%
}\sqrt{\log\left(  k+1\right)  }\right)  \left\Vert q_{1}\right\Vert
_{L^{2}\left(  K\right)  }. \label{v2est2}%
\end{align}

Hence, from \cite[Thm. 3.4]{Ainsworth_parker_I} we deduce that there exists
$\mathbf{v}_{3}\in\mathbf{S}_{k,0}\left(  \mathcal{T}\right)  $ with
$\operatorname*{supp}\mathbf{v}_{3}=K$ such that $\operatorname*{div}%
\mathbf{v}_{3}=q_{3}$ on $K$ and%
\begin{align*}
h_{K}^{-1}\left\Vert \mathbf{v}_{3}\right\Vert _{\mathbf{L}^{2}\left(
K\right)  }+\left\Vert \nabla\mathbf{v}_{3}\right\Vert _{\mathbb{L}^{2}\left(
K\right)  }\leq &  C_{\operatorname*{V}}\left\Vert q_{3}\right\Vert
_{L^{2}\left(  K\right)  }\overset{\text{(\ref{defq3}), (\ref{v2est1})}}{\leq
}C_{\operatorname*{V}}\left(  1+C_{1}\right)  \left\Vert q_{2}\right\Vert
_{L^{2}\left(  K\right)  }\\
&  \overset{\text{(\ref{q2est})}}{\leq}C_{\operatorname*{V}}\left(
1+C_{1}\right)  \left(  1+C_{\operatorname*{CR}}\sqrt{\log\left(  k+1\right)
}\right)  \left\Vert q_{1}\right\Vert _{L^{2}\left(  K\right)  }.
\end{align*}
In this way we have constructed the function $\mathbf{v}_{K}\in\mathbf{CR}%
_{k,0}\left(  \mathcal{T}\right)  $ by%
\[
\mathbf{v}_{K}=\mathbf{v}_{1}+\mathbf{v}_{2}+\mathbf{v}_{3}%
\]
such that $\operatorname*{div}\mathbf{v}_{K}=q_{1}$ on $K$,
$\operatorname*{supp}\mathbf{v}_{K}\subset\omega_{E}$ and%
\begin{align*}
\left\Vert \nabla_{\mathcal{T}}\mathbf{v}_{K}\right\Vert _{\mathbb{L}%
^{2}\left(  \Omega\right)  }  &  =\left\Vert \nabla_{\mathcal{T}}%
\mathbf{v}_{K}\right\Vert _{\mathbb{L}^{2}\left(  \omega_{E}\right)  }\leq
\sum_{\ell=1}^{3}\left\Vert \nabla_{\mathcal{T}}\mathbf{v}_{\ell}\right\Vert
_{\mathbb{L}^{2}\left(  \Omega\right)  }\\
&  \leq C_{2}\sqrt{\log\left(  k+1\right)  }\left\Vert q_{1}\right\Vert
_{L^{2}\left(  K\right)  }\leq C_{2}\sqrt{\log\left(  k+1\right)  }\left(
\left\Vert q\right\Vert _{L^{2}\left(  K\right)  }+\left\Vert \nabla
\mathbf{v}_{0}\right\Vert _{\mathbb{L}^{2}\left(  K\right)  }\right)  ,
\end{align*}
where $C_{2}$ only depends on the shape-regularity of the mesh and
$\alpha_{\Omega}$ through the constants $C_{\operatorname*{V}}$,
$C_{\operatorname*{CR}}$, $C_{1}$. Let $\mathbf{v}_{q}:=\mathbf{v}_{0}%
+\sum_{K\in\mathcal{T}\backslash\mathcal{T}^{\prime}}\mathbf{v}_{K}$ and note
that by construction%
\[
\operatorname*{div}\nolimits_{\mathcal{T}}\mathbf{v}_{q}=q\quad\text{on
}\Omega\backslash\overline{\Omega^{\prime}}%
\]
and%
\begin{align*}
\left\Vert \nabla_{\mathcal{T}}\mathbf{v}_{q}\right\Vert _{\mathbb{L}%
^{2}\left(  \Omega\right)  }  &  \leq\left\Vert \nabla\mathbf{v}%
_{0}\right\Vert _{\mathbb{L}^{2}\left(  \Omega\right)  }+\sum_{K\in
\mathcal{T}\backslash\mathcal{T}^{\prime}}\left\Vert \nabla_{\mathcal{T}%
}\mathbf{v}_{K}\right\Vert _{\mathbb{L}^{2}\left(  \Omega\right)  }\\
&  \leq\left\Vert \nabla\mathbf{v}_{0}\right\Vert _{\mathbb{L}^{2}\left(
\Omega\right)  }+2C_{2}\sum_{K\in\mathcal{T}\backslash\mathcal{T}^{\prime}%
}\sqrt{\log\left(  k+1\right)  }\left(  \left\Vert q\right\Vert _{L^{2}\left(
K\right)  }+\left\Vert \nabla\mathbf{v}_{0}\right\Vert _{\mathbb{L}^{2}\left(
K\right)  }\right) \\
&  \leq C\sqrt{\log\left(  k+1\right)  }\left\Vert q\right\Vert _{L^{2}\left(
\Omega\right)  }.
\end{align*}
Finally, the linear map $\Pi_{\mathcal{T},k}^{\operatorname*{inv}}%
:\mathbb{P}_{k-1,0}\left(  \mathcal{T}\right)  \rightarrow\mathbf{CR}%
_{k,0}\left(  \mathcal{T}\right)  $ is defined by%
\[
\Pi_{\mathcal{T},k}^{\operatorname*{inv}}q=\mathbf{v}_{q}+\Pi_{\mathcal{T}%
^{\prime},k}^{\operatorname*{inv}}\left(  \left.  \left(
q-\operatorname*{div}\mathbf{v}_{q}\right)  \right\vert _{\Omega^{\prime}%
}\right)
\]
and satisfies $\operatorname*{div}\circ\Pi_{\mathcal{T},k}%
^{\operatorname*{inv}}=\operatorname*{Id}$ on $\mathbb{P}_{k-1,0}\left(
\mathcal{T}\right)  $ and%
\[
\left\Vert \nabla_{\mathcal{T}}\Pi_{\mathcal{T},k}^{\operatorname*{inv}%
}q\right\Vert _{\mathbb{L}^{2}\left(  \Omega\right)  }\leq C_{3}\sqrt
{\log\left(  k+1\right)  }C_{\mathcal{T}^{\prime}}\left\Vert q\right\Vert
_{L^{2}\left(  \Omega\right)  }%
\]
for some $C_{3}$ which only depends on the shape-regularity of the mesh and
$\alpha_{\Omega}$.
%TCIMACRO{\TeXButton{End Proof}{\endproof}}%
%BeginExpansion
\endproof
%EndExpansion

By iterating this argument we can prove Theorem \ref{Theomain} for the case
(\ref{Assumptab}).

\begin{theorem}
Let $\mathcal{T}$ be a conforming triangulation which contains at least one
interior vertex. Let $k\geq5$ be odd and let $L\in\mathbb{N}_{0}$ the number
of steps in the construction (\ref{stepbystepextension}). Then, the inf-sup
constant for the corresponding Crouzeix-Raviart discretization satisfies%
\[
c_{\mathcal{T},k}\geq c_{\mathcal{T}}\left(  \log\left(  k+1\right)  \right)
^{-\left(  L+1\right)  /2}%
\]
for some constant $c_{\mathcal{T}}$ which depends only on the shape-regularity
of the mesh and $\alpha_{\Omega}$. If every triangle in $\mathcal{T}$ has an
interior vertex, then $L=0$.
\end{theorem}

\subsection{The case of even $k$\label{SecEvenk}}

This case is slightly simpler that the case of odd $k$ since the
non-conforming Crouzeix-Raviart functions for even $k$ have smaller support
(i.e., one triangle) compared to two triangles (which share an edge) for odd
$k$.

In this section, we assume%
\begin{equation}%
\begin{array}
[c]{ll}%
\text{a)} & k\geq4\text{ is even and}\\
\text{b)} & \mathcal{T}\text{ is a conforming triangulation and contains more
than a single triangle.}%
\end{array}
\label{asseven}%
\end{equation}

\begin{remark}
\label{RemTle1}It is easy to verify that $\left\vert \mathcal{T}\right\vert
>1$ implies that there exists a mapping $\mathfrak{K}:\mathcal{C}%
_{\mathcal{T}}\left(  \eta\right)  $ $\rightarrow\mathcal{T}$ with
$\mathfrak{K}\left(  \mathbf{z}\right)  \in\mathcal{T}_{\mathbf{z}}$ and not
all vertices of $\mathfrak{K}\left(  \mathbf{z}\right)  $ are $\eta$-critical.
\end{remark}

For an $\eta$-critical point $\mathbf{z}\in\mathcal{C}_{\mathcal{T}}\left(
\eta\right)  $, let $n_{\mathbf{z}}:=\left\vert \mathcal{T}_{\mathbf{z}%
}\right\vert $ and fix a counterclockwise numbering of the triangles in
\begin{equation}
\mathcal{T}_{\mathbf{z}}=\left\{  K_{j}^{\mathbf{z}}:1\leq j\leq
n_{\mathbf{z}}\right\}  \label{tznumbconv}%
\end{equation}
such that $K_{j}^{\mathbf{z}}$ and $K_{j+1}^{\mathbf{z}}$ share an edge for
all $1\leq j\leq n_{\mathbf{z}}-1$. With this notation at hand, the functional
$A_{\mathcal{T},\mathbf{z}}$ is given by%
\[
A_{\mathcal{T},\mathbf{z}}q:=\sum_{j=1}^{n_{\mathbf{z}}}\left(  -1\right)
^{j}\left.  q\right\vert _{K_{j}^{\mathbf{z}}}\left(  \mathbf{z}\right)
\quad\forall q\in\mathbb{P}_{k-1,0}\left(  \mathcal{T}\right)  .
\]

\begin{lemma}
Let assumption (\ref{asseven}) be satisfied. There exists a constant $\eta
_{2}>0$ which only depends on the shape-regularity of the mesh and
$\alpha_{\Omega}$ such that for any fixed $0\leq\eta<\eta_{2}$ and any
$q\in\mathbb{P}_{k-1,0}\left(  \mathcal{T}\right)  $, there exists some
$\mathbf{v}_{q}\in\mathbf{CR}_{k,0}\left(  \mathcal{T}\right)  $ such that%
\begin{align}
\int_{K}\operatorname{div}\mathbf{v}_{q}  &  =0\quad\forall K\in
\mathcal{T},\label{intcondeven}\\
q-\operatorname*{div}\nolimits_{\mathcal{T}}\mathbf{v}_{q}  &  \in
M_{\eta,k-1}^{\operatorname*{SV}}\left(  \mathcal{T}\right)
\label{qcorecteven}%
\end{align}
and%
\begin{equation}
\left\Vert \mathbf{v}_{q}\right\Vert _{\mathbf{H}^{1}\left(  \mathcal{T}%
\right)  }\leq C_{\operatorname*{CR}}\sqrt{\log\left(  k+1\right)  }\left\Vert
q\right\Vert _{L^{2}\left(  \Omega\right)  }. \label{vq1}%
\end{equation}
The constant $C_{\operatorname*{CR}}$ depends only on the shape-regularity of
the mesh and $\alpha_{\Omega}$.
\end{lemma}

%

%TCIMACRO{\TeXButton{Proof}{\proof}}%
%BeginExpansion
\proof
%EndExpansion
For a triangle $K\in\mathcal{T}$, we set $\mathcal{C}_{K}\left(  \eta\right)
:=\left\{  \mathbf{z}\in\mathcal{V}\left(  K\right)  \cap\mathcal{C}%
_{\mathcal{T}}\left(  \eta\right)  \right\}  $ and $\mathcal{C}_{K}%
^{\operatorname*{active}}\left(  \eta\right)  :=\left\{  \mathbf{z}%
\in\mathcal{C}_{K}\left(  \eta\right)  :K=\mathfrak{K}\left(  \mathbf{z}%
\right)  \right\}  $. Their cardinalities are denoted by $n_{K}:=\left\vert
\mathcal{C}_{K}\left(  \eta\right)  \right\vert $ and $n_{K}%
^{\operatorname*{active}}:=\left\vert \mathcal{C}_{K}^{\operatorname*{active}%
}\left(  \eta\right)  \right\vert $. Note that $0\leq n_{K}%
^{\operatorname*{active}}\leq n_{K}\leq2$ (cf. Rem. \ref{RemTle1}). We number
the vertices $\mathbf{V}_{j}$ in $K$ with the convention $\left\{
\mathbf{V}_{j}:1\leq j\leq n_{K}\right\}  =\mathcal{C}_{K}\left(  \eta\right)
$ and $\left\{  \mathbf{V}_{j}:1\leq j\leq n_{K}^{\operatorname*{active}%
}\right\}  =\mathcal{C}_{K}^{\operatorname*{active}}\left(  \eta\right)  $.
The angle in $K$ at $\mathbf{V}_{j}$ is denoted by $\alpha_{j}$. Let $E_{j}$
be the edge in $K$ opposite to $\mathbf{V}_{j}$ and let $\mathbf{n}_{j}$
denote the outward unit normal vector at $E_{j}$.

In a similar way as for the construction of $\tilde{B}_{k,E}%
^{\operatorname*{CR}}$ in Lemma \ref{Lemtrianglemean} (and employing Lemma
\ref{LemBKH1/2norm} instead of Lemma \ref{LemBEH1/2norm} for (\ref{wtildeE}))
there exists a function $\tilde{B}_{k,K}^{\operatorname*{CR}}\in
\operatorname*{CR}_{k,0}\left(  \mathcal{T}\right)  $ with

\begin{enumerate}
\item $\operatorname*{supp}\tilde{B}_{k,K}^{\operatorname*{CR}}=K,$

\item for all $K\in\mathcal{T}$, for all $\mathbf{z}\in\mathcal{V}\left(
K\right)  $%
\begin{equation}
\nabla\left(  \left.  \tilde{B}_{k,K}^{\operatorname*{CR}}\right\vert
_{K}\right)  \left(  \mathbf{z}\right)  =\nabla\left(  \left.  B_{k,K}%
^{\operatorname*{CR}}\right\vert _{K}\right)  \left(  \mathbf{z}\right)  ,
\label{divtrianglerel}%
\end{equation}

\item for all $K\in\mathcal{T}$%
\begin{equation}
\left.  \left.  \tilde{B}_{k,K}^{\operatorname*{CR}}\right\vert _{K}%
\right\vert _{\partial K}=\left.  \left.  B_{k,K}^{\operatorname*{CR}%
}\right\vert _{K}\right\vert _{\partial K}, \label{traceKequal}%
\end{equation}

\item for all $K\in\mathcal{T}$ and any $\mathbf{c}\in\mathbb{R}^{2}$%
\begin{equation}
\int_{K}\operatorname{div}_{\mathcal{T}}\left(  \tilde{B}_{k,K}%
^{\operatorname*{CR}}\mathbf{c}\right)  =0. \label{divKeven=0}%
\end{equation}

\item The piecewise gradient is bounded by%
\begin{equation}
\left\Vert \tilde{B}_{k,K}^{\operatorname*{CR}}\right\Vert _{H^{1}\left(
\mathcal{T}\right)  }\leq C\sqrt{\log\left(  k+1\right)  }. \label{Btildeest}%
\end{equation}

\end{enumerate}

For $j=1,2$, we define%
\begin{equation}%
%TCIMACRO{\TeXButton{boldpsi}{\mbox{\boldmath$ \psi$}}}%
%BeginExpansion
\mbox{\boldmath$ \psi$}%
%EndExpansion
_{k,K}^{\operatorname*{CR},j}:=\left\{
\begin{array}
[c]{ll}%
\tilde{B}_{k,K}^{\operatorname*{CR}}\mathbf{n}_{j} & \text{in }K,\\
0 & \text{in }\Omega\backslash K,
\end{array}
\right.  \label{defboldpsieven}%
\end{equation}
i.e., we fix $\mathbf{v}_{K}:=\mathbf{n}_{1}$ and $\mathbf{w}_{K}%
=\mathbf{n}_{2}$ in (\ref{basisvel}). The divergence of $%
%TCIMACRO{\TeXButton{boldpsi}{\mbox{\boldmath$ \psi$}}}%
%BeginExpansion
\mbox{\boldmath$ \psi$}%
%EndExpansion
_{k,K}^{\operatorname*{CR},j}$ evaluated at a vertex $\mathbf{V}_{s}$,
$s=1,2$, is given by%
\begin{align*}
\operatorname{div}\left(  \left.
%TCIMACRO{\TeXButton{boldpsi}{\mbox{\boldmath$ \psi$}}}%
%BeginExpansion
\mbox{\boldmath$ \psi$}%
%EndExpansion
_{k,K}^{\operatorname*{CR},j}\right\vert _{K}\right)  \left(  \mathbf{V}%
_{s}\right)   &  \overset{\text{(\ref{divtrianglerel})}}{=}\operatorname{div}%
\left(  \left.  B_{k,K}^{\operatorname*{CR}}\mathbf{n}_{j}\right\vert
_{K}\right)  \left(  \mathbf{V}_{s}\right)  =-%
%TCIMACRO{\dsum \limits_{i=1}^{3}}%
%BeginExpansion
{\displaystyle\sum\limits_{i=1}^{3}}
%EndExpansion
L_{k}^{\prime}\left(  1-2\lambda_{K,i}\left(  \mathbf{V}_{s}\right)  \right)
\partial_{\mathbf{n}_{j}}\lambda_{K,i}\\
&  =k\left(  k+1\right)  \partial_{\mathbf{n}_{j}}\lambda_{K,s}\overset
{\text{(\ref{normcomp})}}{=}\binom{k+1}{2}\frac{\left\vert E_{s}\right\vert
}{\left\vert K\right\vert }\times\left\{
\begin{array}
[c]{ll}%
-1 & j=s,\\
\cos\alpha_{3} & j\neq s.
\end{array}
\right.
\end{align*}

Let $q\in\mathbb{P}_{k-1,0}\left(  \mathcal{T}\right)  $. We choose $%
%TCIMACRO{\TeXButton{bolddelta}{\mbox{\boldmath$ \delta$}}}%
%BeginExpansion
\mbox{\boldmath$ \delta$}%
%EndExpansion
_{K}:=\left(  \delta_{K,j}\right)  _{j=1}^{2}$ by the conditions for $s=1,2$%
\begin{equation}
A_{\mathcal{T},\mathbf{V}_{s}}\left(  \operatorname*{div}%
\nolimits_{\mathcal{T}}\left(  \sum_{j=1}^{2}\delta_{K,j}%
%TCIMACRO{\TeXButton{boldpsi}{\mbox{\boldmath$ \psi$}}}%
%BeginExpansion
\mbox{\boldmath$ \psi$}%
%EndExpansion
_{k,K}^{\operatorname*{CR},j}\right)  \right)  \overset{!}{=}\left\{
\begin{array}
[c]{ll}%
A_{\mathcal{T},\mathbf{V}_{s}}\left(  q\right)  & \text{if }\mathbf{V}_{s}%
\in\mathcal{C}_{K}^{\operatorname*{active}}\left(  \eta\right)  ,\\
0 & \text{otherwise.}%
\end{array}
\right.  \label{condkeven}%
\end{equation}
For $s=1,2$, let $\ell_{s}$ be defined by $K_{\ell_{s}}^{\mathbf{V}_{s}}=K$
(cf. (\ref{tznumbconv})). Then,%
\begin{align*}
A_{\mathcal{T},\mathbf{V}_{s}}\left(  \operatorname*{div}%
\nolimits_{\mathcal{T}}\left(  \sum_{j=1}^{2}\delta_{K,j}%
%TCIMACRO{\TeXButton{boldpsi}{\mbox{\boldmath$ \psi$}}}%
%BeginExpansion
\mbox{\boldmath$ \psi$}%
%EndExpansion
_{k,K}^{\operatorname*{CR},j}\right)  \right)   &  =\left(  -1\right)
^{\ell_{s}}\sum_{j=1}^{2}\delta_{K,j}\left(  \operatorname*{div}\left.
%TCIMACRO{\TeXButton{boldpsi}{\mbox{\boldmath$ \psi$}}}%
%BeginExpansion
\mbox{\boldmath$ \psi$}%
%EndExpansion
_{k,K}^{\operatorname*{CR},j}\right\vert _{K}\right)  \left(  \mathbf{V}%
_{s}\right) \\
&  =\mathbf{M}_{K}%
%TCIMACRO{\TeXButton{bolddelta}{\mbox{\boldmath$ \delta$}}}%
%BeginExpansion
\mbox{\boldmath$ \delta$}%
%EndExpansion
_{K}%
\end{align*}
for%
\[
\mathbf{M}_{K}=\left(  -1\right)  ^{\ell_{s}+1}\binom{k+1}{2}\left\vert
K\right\vert ^{-1}\left[
\begin{array}
[c]{cc}%
\left\vert E_{1}\right\vert  & -\left\vert E_{1}\right\vert \cos\alpha_{3}\\
-\left\vert E_{2}\right\vert \cos\alpha_{3} & \left\vert E_{2}\right\vert
\end{array}
\right]  .
\]
We define $\mathbf{r}_{K}=\left(  r_{K,s}\right)  _{s=1}^{2}$ by%
\[
r_{K,s}:=\left\{
\begin{array}
[c]{ll}%
A_{\mathcal{T},\mathbf{V}_{s}}\left(  q\right)  & \text{if }\mathbf{V}_{s}%
\in\mathcal{C}_{K}^{\operatorname*{active}}\left(  \eta\right)  ,\\
0 & \text{otherwise}%
\end{array}
\right.
\]
so that $%
%TCIMACRO{\TeXButton{bolddelta}{\mbox{\boldmath$ \delta$}}}%
%BeginExpansion
\mbox{\boldmath$ \delta$}%
%EndExpansion
_{K}$ is the solution of%
\[
\mathbf{M}_{K}%
%TCIMACRO{\TeXButton{bolddelta}{\mbox{\boldmath$ \delta$}}}%
%BeginExpansion
\mbox{\boldmath$ \delta$}%
%EndExpansion
_{K}=\mathbf{r}_{K}.
\]
Observe that%
\[
\det\left[
\begin{array}
[c]{cc}%
\left\vert E_{1}\right\vert  & -\left\vert E_{1}\right\vert \cos\alpha_{3}\\
-\left\vert E_{2}\right\vert \cos\alpha_{3} & \left\vert E_{2}\right\vert
\end{array}
\right]  =\left\vert E_{1}\right\vert \left\vert E_{2}\right\vert \sin
^{2}\alpha_{3}=2\left\vert K\right\vert \sin\alpha_{3}%
\]
and $\sin\alpha_{3}\geq\sin\phi_{\mathcal{T}}>0$ due to the shape-regularity
of the mesh. For the coefficient $%
%TCIMACRO{\TeXButton{bolddelta}{\mbox{\boldmath$ \delta$}}}%
%BeginExpansion
\mbox{\boldmath$ \delta$}%
%EndExpansion
_{K}$ we get explicitly%
\begin{equation}%
%TCIMACRO{\TeXButton{bolddelta}{\mbox{\boldmath$ \delta$}}}%
%BeginExpansion
\mbox{\boldmath$ \delta$}%
%EndExpansion
_{K}=\frac{\left(  -1\right)  ^{\ell_{s}+1}}{\left(  k+1\right)  k\sin
\alpha_{3}}\left[
\begin{array}
[c]{cc}%
\left\vert E_{2}\right\vert  & \left\vert E_{1}\right\vert \cos\alpha_{3}\\
\left\vert E_{2}\right\vert \cos\alpha_{3} & \left\vert E_{1}\right\vert
\end{array}
\right]  \mathbf{r}_{K} \label{Defalphazeven}%
\end{equation}
with an estimate%
\[
\left\Vert
%TCIMACRO{\TeXButton{bolddelta}{\mbox{\boldmath$ \delta$}}}%
%BeginExpansion
\mbox{\boldmath$ \delta$}%
%EndExpansion
_{K}\right\Vert \leq C\frac{h_{K}}{k\left(  k+1\right)  }\left\Vert
\mathbf{r}_{K}\right\Vert \overset{\text{Lem. \ref{LemFunctional}}}{\leq
}C\left\Vert q\right\Vert _{L^{2}\left(  \omega_{K}\right)  }\quad\text{if
}\mathcal{C}_{K}^{\operatorname*{active}}\left(  \eta\right)  \neq\emptyset,
\]
where $C$ only depends on the shape-regularity of the mesh. Note that this is
the analogue for even $k$ to (\ref{estboldalpha}). If $\mathcal{C}%
_{K}^{\operatorname*{active}}\left(  \eta\right)  =\emptyset$, it holds $%
%TCIMACRO{\TeXButton{bolddelta}{\mbox{\boldmath$ \delta$}}}%
%BeginExpansion
\mbox{\boldmath$ \delta$}%
%EndExpansion
_{K}=\mathbf{0}$.

We define the global function%
\begin{equation}
\mathbf{v}_{q}:=\sum_{\substack{K\in\mathcal{T}\\\mathcal{C}_{K}%
^{\operatorname*{active}}\left(  \eta\right)  \neq\emptyset}}\sum_{j=1}%
^{2}\delta_{K,j}%
%TCIMACRO{\TeXButton{boldpsi}{\mbox{\boldmath$ \psi$}}}%
%BeginExpansion
\mbox{\boldmath$ \psi$}%
%EndExpansion
_{k,K}^{\operatorname*{CR},j}. \label{v1globaleven}%
\end{equation}
From (\ref{divKeven=0}) we conclude that $\mathbf{v}_{q}$ satisfies
(\ref{intcondeven}).

Next, we verify (\ref{qcorecteven}). Let $\mathbf{y}\in\mathcal{C}%
_{\mathcal{T}}\left(  \eta\right)  $ and recall the notation and convention as
in (\ref{tznumbconv}). Let $K_{\ell}^{\mathbf{y}}=\mathfrak{K}\left(
\mathbf{y}\right)  $. Then (\ref{qcorecteven}) follows from%
\begin{align*}
A_{\mathcal{T},\mathbf{y}}\left(  \operatorname*{div}\nolimits_{\mathcal{T}%
}\mathbf{v}_{q}\right)   &  =\left(  -1\right)  ^{\ell}\left(
\operatorname{div}\left.  \mathbf{v}_{q}\right\vert _{K_{\ell}^{\mathbf{y}}%
}\right)  \left(  \mathbf{y}\right)  =\left(  -1\right)  ^{\ell}\left(
\operatorname{div}\left.  \sum_{j=1}^{2}\delta_{K_{\ell}^{\mathbf{y}},j}%
%TCIMACRO{\TeXButton{boldpsi}{\mbox{\boldmath$ \psi$}}}%
%BeginExpansion
\mbox{\boldmath$ \psi$}%
%EndExpansion
_{k,K_{\ell}^{\mathbf{y}}}^{\operatorname*{CR},j}\right\vert _{K_{\ell
}^{\mathbf{y}}}\right)  \left(  \mathbf{y}\right) \\
&  =A_{\mathcal{T},\mathbf{y}}\left(  q\right)  .
\end{align*}
The estimate%
\[
\left\Vert \nabla%
%TCIMACRO{\TeXButton{boldpsi}{\mbox{\boldmath$ \psi$}}}%
%BeginExpansion
\mbox{\boldmath$ \psi$}%
%EndExpansion
_{k,K}^{\operatorname*{CR},j}\right\Vert _{\mathbb{L}^{2}\left(  K\right)
}\leq C\sqrt{\log\left(  k+1\right)  }%
\]
for a constant $C$ which only depends on the shape-regularity of the mesh and
$\alpha_{\Omega}$ follows directly from (\ref{Btildeest}) and the final
estimate (\ref{vq1}) is derived by repeating the arguments as in the proof of
Lemma \ref{LemProject}.%
%TCIMACRO{\TeXButton{End Proof}{\endproof}}%
%BeginExpansion
\endproof
%EndExpansion

This lemma allows us to extend Definition \ref{DefPiCR} to the case of even
$k$ by defining the coefficients $%
%TCIMACRO{\TeXButton{bolddelta}{\mbox{\boldmath$ \delta$}}}%
%BeginExpansion
\mbox{\boldmath$ \delta$}%
%EndExpansion
_{K}$ by (\ref{Defalphazeven}) and the functions $%
%TCIMACRO{\TeXButton{boldpsi}{\mbox{\boldmath$ \psi$}}}%
%BeginExpansion
\mbox{\boldmath$ \psi$}%
%EndExpansion
_{k,K}^{\operatorname*{CR}}$ by (\ref{defboldpsieven}) and set (cf.
(\ref{v1globaleven})) $\Pi_{k}^{\operatorname*{CR}}q:=\mathbf{v}_{q}$. Since
$\left(  I-\Pi_{k}^{\operatorname*{CR}}\right)  q\in M_{\eta,k-1}%
^{\operatorname*{SV}}\left(  \mathcal{T}\right)  $ we may apply the further
steps in the proof of Lemma \ref{LemmaFinOdd} to obtain the inf-sup stability
for even $k$.

\begin{theorem}
Let $\mathcal{T}$ be a conforming triangulation satisfies (\ref{asseven}).
Then, the inf-sup constant for the corresponding Crouzeix-Raviart
discretization satisfies%
\[
c_{\mathcal{T},k}\geq c_{\mathcal{T}}\left(  \log\left(  k+1\right)  \right)
^{-1/2}%
\]
for some constant $c_{\mathcal{T}}$ which depends only on the shape-regularity
of the mesh and $\alpha_{\Omega}$.
\end{theorem}

\section{Conclusion\label{SecConcl}}

In this paper, we have derived lower bounds for the inf-sup constant for
Crouzeix-Raviart elements for the Stokes equation which are explicit with
respect to the polynomial degree $k$ and are independent of the mesh size.

\begin{enumerate}
\item The inf-sup constant can be bounded from below by $c_{\mathcal{T},k}\geq
c_{\mathcal{T}}\left(  \log\left(  k+1\right)  \right)  ^{-1/2}$ if

\begin{enumerate}
\item for odd $k\geq3$,

\begin{enumerate}
\item $\mathcal{T}$ has at least one interior point and

\item for $k\geq5$, $\mathcal{T}$ has no acute critical point,
\end{enumerate}

\item for even $k\geq4$, $\mathcal{T}$ contains more than one triangle,
\end{enumerate}

\item If for odd $k$, condition 1.a.ii. is not satisfied but a step-by-step
construction (\ref{stepbystepextension}) for some $L\geq1$ is possible, then,
$c_{\mathcal{T},k}\geq c_{\mathcal{T}}\left(  \log\left(  k+1\right)  \right)
^{-\left(  L+1\right)  /2}$.
\end{enumerate}

Finally, we compare these findings with some other stable pairs of Stokes
elements on triangulations in the literature. The element $\left(
\mathbf{S}_{k,0}\left(  \mathcal{T}\right)  ,\mathbb{P}_{k-2,0}\left(
\mathcal{T}\right)  \right)  $ has a discrete inf-sup constant which can be
estimated from below by $Ck^{-3}$ (see \cite{Schwab_Suri_hp_infsup},
\cite{Stenberg_Suri_hp}). The discrete inf-sup constant for the Scott-Vogelius
element $\left(  \mathbf{S}_{k,0}\left(  \mathcal{T}\right)  ,M_{0,k-1}%
^{\operatorname*{SV}}\left(  \mathcal{T}\right)  \right)  $ for $k\geq4$ can
be estimated from below by $c\Theta_{\min}k^{-m}$ for some integer $m$
sufficiently large (see \cite{ScottVogelius}, \cite{vogelius1983right}). The
pressure-wired Stokes element $\left(  \mathbf{S}_{k,0}\left(  \mathcal{T}%
\right)  ,M_{\eta,k-1}^{\operatorname*{SV}}\left(  \mathcal{T}\right)
\right)  $ in \ref{Sauter_eta_wired} (again for $k\geq4$) is a mesh-robust
generalization of the Scott-Vogelius element with a lower bound of the inf-sup
constant of the form $c\left(  \Theta_{\min}+\eta\right)  $. In
\cite{Ainsworth_parker_I}, a conforming stable pair $\left(  \mathbf{X}%
_{k}\left(  \mathcal{T}\right)  ,M_{k-1}\left(  \mathcal{T}\right)  \right)  $
of Stokes elements on triangulations is introduced and it is proved that the
discrete inf-sup constant can be estimated from below by $c/\tilde{\Theta
}_{\min}$ for a constant $c$ independent of $h$ and $k$ and $\tilde{\Theta
}_{\min}:=\min_{\mathbf{z}\in\mathcal{V}_{\partial\Omega}\left(
\mathcal{T}\right)  \backslash\mathcal{C}_{\mathcal{T}}}\Theta\left(
\mathbf{z}\right)  $. However, the implementation requires finite elements for
the velocity with $C^{1}$ continuity at the triangle vertices and pressures
which are continuous in the triangle vertices.

\begin{acknowledgement}
Thanks are due to Benedikt Gr\"{a}\ss le, HU Berlin, for fruitful discussions
on Lemma \ref{Lemeta0}.

I am grateful to my colleagues from TU Vienna, Profs. Joachim Sch\"{o}berl and
Markus Melenk. Joachim showed by numerical experiments that the lower bound of
the inf-sup constant in the first arxiv version of the paper, namely
$k^{-1/4}$, might be too pessimistic and \textquotedblleft it should be at
least $k^{-1/6}$\textquotedblright\ and Markus raised the suspicion that the
interpolation argument might be too pessimistic for Legendre polynomials.
\end{acknowledgement}

\appendix

\section{The inverse of the matrix $\mathbf{T}_{\ell}+%
%TCIMACRO{\TeXButton{boldcapdelta}{\mbox{\boldmath$ \Delta$}}}%
%BeginExpansion
\mbox{\boldmath$ \Delta$}%
%EndExpansion
_{\ell}$ in (\ref{DefMlTlDl})\label{AppT_L}}

\begin{lemma}
\label{Tlest}There exists $\eta_{2}>0$ which only depends on the
shape-regularity of the mesh and $\alpha_{\Omega}$ such that for any
$0\leq\eta<\eta_{2}$ the matrix $\mathbf{T}_{\ell}+%
%TCIMACRO{\TeXButton{boldcapdelta}{\mbox{\boldmath$ \Delta$}}}%
%BeginExpansion
\mbox{\boldmath$ \Delta$}%
%EndExpansion
_{\ell}$ in (\ref{DefMlTlDl}) is invertible and there exists a constant $C$
depending only on the shape-regularity of the mesh and $\alpha_{\Omega}$ such
that (cf. Notation \ref{Notation})%
\[
\left\Vert \left(  \mathbf{T}_{\ell}+%
%TCIMACRO{\TeXButton{boldcapdelta}{\mbox{\boldmath$ \Delta$}}}%
%BeginExpansion
\mbox{\boldmath$ \Delta$}%
%EndExpansion
_{\ell}\right)  ^{-1}\right\Vert \leq C.
\]

\end{lemma}

Note that the matrix $\mathbf{T}_{\ell}$ in (\ref{defTl}) is the same as the
matrix $\mathbf{T}_{n,%
%TCIMACRO{\TeXButton{boldalpha}{\mbox{\boldmath$ \alpha$}}}%
%BeginExpansion
\mbox{\boldmath$ \alpha$}%
%EndExpansion
}$ which has been analysed in \cite[(3.36)]{CCSS_CR_1}. In particular, the
formula%
\[
\det\mathbf{T}_{\ell}=\frac{\sin\left(  \sum_{j=1}^{n_{\ell}+1}\alpha
_{\ell,j,1}\right)  }{%
%TCIMACRO{\dprod \nolimits_{j=1}^{n_{\ell}+1}}%
%BeginExpansion
{\displaystyle\prod\nolimits_{j=1}^{n_{\ell}+1}}
%EndExpansion
\sin\alpha_{\ell,j,1}}%
\]
was proved.

Next we show that the sum $\sum_{j=1}^{n_{\ell}+1}\alpha_{\ell,j,1}$ is
bounded away from $0$ and $\pi$. The bound $\sum_{j=1}^{n_{\ell}+1}%
\alpha_{\ell,j,1}\geq\varphi_{\mathcal{T}}$ follows from Remark \ref{Remangle}%
. Since $\mathbf{z}_{\ell,j}$, $1\leq j\leq n_{\ell}$ are $\eta-$critical
points the sum of both angles adjacent to $E_{\ell,j}$ at $\mathbf{z}_{\ell
,j}$ satisfy%
\[
\sin\left(  \alpha_{\ell,j,2}+\alpha_{\ell,j+1,3}\right)  \leq\eta\text{. }%
\]
We write $\alpha_{\ell,j,2}+\alpha_{\ell,j+1,3}=:\pi+\delta_{j}$. From the
proof of Lemma \ref{Lemangle}, in particular from the estimate (\ref{deltaest}%
) we conclude that $\left\vert \delta_{j}\right\vert \leq c_{2}\eta$.

Since all points $\mathbf{z}_{\ell,j}$ are edge-connected to the same point
$\mathbf{z}_{\ell}$, the number $n_{\ell}$ is bounded from above by a constant
$n_{\max}$ which only depends on the shape-regularity of the mesh. Hence,%
\begin{align*}
\sum_{j=1}^{n_{\ell}+1}\alpha_{\ell,j,1}  &  =\sum_{j=1}^{n_{\ell}+1}\left(
\pi-\alpha_{\ell,j,2}-\alpha_{\ell,j,3}\right)  =\left(  n_{\ell}+1\right)
\pi-\alpha_{\ell,1,3}-\alpha_{\ell,n_{\ell}+1,2}-\sum_{j=1}^{n_{\ell}}\left(
\alpha_{\ell,j,2}+\alpha_{\ell,j+1,3}\right) \\
&  =\left(  n_{\ell}+1\right)  \pi-\alpha_{\ell,1,3}-\alpha_{\ell,n_{\ell
}+1,2}-\sum_{j=1}^{n_{\ell}}\left(  \pi+\delta_{j}\right)  =\pi-\alpha
_{\ell,1,3}-\alpha_{\ell,n_{\ell}+1,2}-\sum_{j=1}^{n_{\ell}}\delta_{j}\\
&  \leq\pi-2\varphi_{\mathcal{T}}+n_{\max}c_{2}\eta.
\end{align*}
By adjusting the constant $\eta_{0}$ in Lemma \ref{Lemangle} to $\eta
_{1}:=\min\left\{  \eta_{0},\varphi_{\mathcal{T}}/\left(  n_{\max}%
c_{2}\right)  \right\}  $ it follows that%
\[
\sum_{j=1}^{n_{\ell}+1}\alpha_{\ell,j,1}\leq\pi-\varphi_{\mathcal{T}}\text{.}%
\]
By using the trivial estimate $0<\sin\alpha_{\ell,j,1}\leq1$, we may conclude
that%
\begin{equation}
\det\mathbf{T}_{\ell}=\frac{\sin\left(  \sum_{j=1}^{n_{\ell}+1}\alpha
_{\ell,j,1}\right)  }{%
%TCIMACRO{\dprod \nolimits_{j=1}^{n_{\ell}+1}}%
%BeginExpansion
{\displaystyle\prod\nolimits_{j=1}^{n_{\ell}+1}}
%EndExpansion
\sin\alpha_{\ell,j,1}}\geq\sin\varphi_{\mathcal{T}}>0. \label{detTlest}%
\end{equation}
Note that the entries in the matrix $\mathbf{T}_{\ell}$ (cf. \ref{defTl})
satisfy%
\begin{equation}
\left\vert \left(  \mathbf{T}_{\ell}\right)  _{\mathbf{y},\mathbf{z}%
}\right\vert \leq\frac{1}{\sin^{2}\varphi_{\mathcal{T}}} \label{estelTl}%
\end{equation}
and hence the Frobenius norm $\left\Vert \mathbf{T}_{\ell}\right\Vert
_{\operatorname{F}}$ can be estimated by%
\[
\left\Vert \mathbf{T}_{\ell}\right\Vert _{\operatorname{F}}\leq\frac{3n_{\ell
}}{\sin^{2}\varphi_{\mathcal{T}}}\leq\frac{3n_{\max}}{\sin^{2}\varphi
_{\mathcal{T}}}.
\]
It is well known that $\left\Vert \mathbf{T}_{\ell}\right\Vert \leq\left\Vert
\mathbf{T}_{\ell}\right\Vert _{\operatorname{F}}$ and hence the bound on
$\left\Vert \mathbf{T}_{\ell}\right\Vert $ follows.

We combine (\ref{detTlest}), (\ref{estelTl}), and $n_{\ell}\leq n_{\max}$ to
obtain by Cramer's rule that there exists a constant $C$ which only depends on
the shape-regularity such that%
\[
\left\vert \left(  \mathbf{T}_{\ell}^{-1}\right)  _{\mathbf{y},\mathbf{z}%
}\right\vert \leq C.
\]
By the same arguments as before we conclude that $\left\Vert \mathbf{T}_{\ell
}^{-1}\right\Vert \leq\tilde{C}$ for a constant $\tilde{C}$ which only depends
on the shape-regularity of the mesh.

Next we estimate $\left\Vert
%TCIMACRO{\TeXButton{boldcapdelta}{\mbox{\boldmath$ \Delta$}}}%
%BeginExpansion
\mbox{\boldmath$ \Delta$}%
%EndExpansion
_{\ell}\right\Vert $. Since $%
%TCIMACRO{\TeXButton{boldcapdelta}{\mbox{\boldmath$ \Delta$}}}%
%BeginExpansion
\mbox{\boldmath$ \Delta$}%
%EndExpansion
_{\ell}$ is diagonal it suffices to estimate the diagonal entries%
\[
\left\vert \frac{\sin\left(  \alpha_{\ell,j,2}+\alpha_{\ell,j+1,3}\right)
}{\sin\alpha_{\ell,j,2}\sin\alpha_{\ell,j+1,3}}\right\vert =\left\vert
\frac{\sin\left(  \pi+\delta_{\ell}\right)  }{\sin\alpha_{\ell,j,2}\sin
\alpha_{\ell,j+1,3}}\right\vert =\left\vert \frac{\sin\delta_{\ell}}%
{\sin\alpha_{\ell,j,2}\sin\alpha_{\ell,j+1,3}}\right\vert \leq\frac{c_{2}\eta
}{\sin^{2}\varphi_{\mathcal{T}}}.
\]
We write $\mathbf{T}_{\ell}+%
%TCIMACRO{\TeXButton{boldcapdelta}{\mbox{\boldmath$ \Delta$}}}%
%BeginExpansion
\mbox{\boldmath$ \Delta$}%
%EndExpansion
_{\ell}=\mathbf{T}_{\ell}\left(  \mathbf{I}+\mathbf{T}_{\ell}^{-1}%
%TCIMACRO{\TeXButton{boldcapdelta}{\mbox{\boldmath$ \Delta$}}}%
%BeginExpansion
\mbox{\boldmath$ \Delta$}%
%EndExpansion
_{\ell}\right)  $ and obtain%
\[
\left\Vert \mathbf{T}_{\ell}^{-1}%
%TCIMACRO{\TeXButton{boldcapdelta}{\mbox{\boldmath$ \Delta$}}}%
%BeginExpansion
\mbox{\boldmath$ \Delta$}%
%EndExpansion
_{\ell}\right\Vert \leq\left\Vert \mathbf{T}_{\ell}^{-1}\right\Vert
\left\Vert
%TCIMACRO{\TeXButton{boldcapdelta}{\mbox{\boldmath$ \Delta$}}}%
%BeginExpansion
\mbox{\boldmath$ \Delta$}%
%EndExpansion
_{\ell}\right\Vert \leq\tilde{C}\frac{c_{2}\eta}{\sin^{2}\varphi_{\mathcal{T}%
}}.
\]
Next, we adjust the upper bound $\eta_{1}$ by setting $\eta_{2}:=\min\left\{
\eta_{1},\frac{\sin^{2}\varphi_{\mathcal{T}}}{2\tilde{C}c_{2}}\right\}  $ to
obtain $\left\Vert \mathbf{T}_{\ell}^{-1}%
%TCIMACRO{\TeXButton{boldcapdelta}{\mbox{\boldmath$ \Delta$}}}%
%BeginExpansion
\mbox{\boldmath$ \Delta$}%
%EndExpansion
_{\ell}\right\Vert \leq1/2$ with implies the invertibility of $\mathbf{T}%
_{\ell}+%
%TCIMACRO{\TeXButton{boldcapdelta}{\mbox{\boldmath$ \Delta$}}}%
%BeginExpansion
\mbox{\boldmath$ \Delta$}%
%EndExpansion
_{\ell}$ with bound%
\[
\left\Vert \left(  \mathbf{T}_{\ell}+%
%TCIMACRO{\TeXButton{boldcapdelta}{\mbox{\boldmath$ \Delta$}}}%
%BeginExpansion
\mbox{\boldmath$ \Delta$}%
%EndExpansion
_{\ell}\right)  ^{-1}\right\Vert \leq2\tilde{C}.
\]%
%TCIMACRO{\TeXButton{End Proof}{\endproof}}%
%BeginExpansion
\endproof
%EndExpansion

\section{Estimate of the $H^{1/2}$ norm of traces of non-conforming
Crouzeix-Raviart functions\label{AppH1/2Est}}

In this appendix, we prove the norm estimate (\ref{Btildeest}) for $\tilde
{B}_{k,K}^{\operatorname*{CR}}$ and (\ref{pwgrad}) for $\tilde{B}%
_{k,E}^{\operatorname*{CR}}$. We first introduce some norms and semi-norms on
the unit interval $I:=\left[  -1,1\right]  $ in a formal way:%
\begin{equation}%
\begin{array}
[c]{ll}%
\left\vert u\right\vert _{H^{1/2}\left(  I\right)  }:=\left(
%TCIMACRO{\dint _{-1}^{1}}%
%BeginExpansion
{\displaystyle\int_{-1}^{1}}
%EndExpansion%
%TCIMACRO{\dint _{-1}^{1}}%
%BeginExpansion
{\displaystyle\int_{-1}^{1}}
%EndExpansion
\left\vert \dfrac{u\left(  s\right)  -u\left(  t\right)  }{s-t}\right\vert
^{2}dsdt\right)  ^{1/2}, & \left\Vert u\right\Vert _{H^{1/2}\left(  I\right)
}:=\left(  \left\Vert u\right\Vert _{L^{2}\left(  I\right)  }^{2}+\left\vert
u\right\vert _{H^{1/2}\left(  I\right)  }^{2}\right)  ^{1/2},\\
\left\vert u\right\vert _{H_{(0,}^{1/2}\left(  I\right)  }:=\left(
%TCIMACRO{\dint _{-1}^{1}}%
%BeginExpansion
{\displaystyle\int_{-1}^{1}}
%EndExpansion
\dfrac{\left\vert u\left(  s\right)  \right\vert ^{2}}{1+s}ds\right)
^{1/2}, & \left\Vert u\right\Vert _{H_{(0,}^{1/2}\left(  I\right)  }:=\left(
\left\Vert u\right\Vert _{H^{1/2}\left(  I\right)  }^{2}+\left\vert
u\right\vert _{H_{(0,}^{1/2}\left(  I\right)  }^{2}\right)  ^{1/2},\\
\left\vert u\right\vert _{H_{,0)}^{1/2}\left(  I\right)  }:=\left(
%TCIMACRO{\dint _{-1}^{1}}%
%BeginExpansion
{\displaystyle\int_{-1}^{1}}
%EndExpansion
\dfrac{\left\vert u\left(  s\right)  \right\vert ^{2}}{1-s}ds\right)
^{1/2}, & \left\Vert u\right\Vert _{H_{,0)}^{1/2}\left(  I\right)  }:=\left(
\left\Vert u\right\Vert _{H^{1/2}\left(  I\right)  }^{2}+\left\vert
u\right\vert _{H_{,0)}^{1/2}\left(  I\right)  }^{2}\right)  ^{1/2},\\
\left\vert u\right\vert _{H_{00}^{1/2}\left(  I\right)  }:=\left(  \left\vert
u\right\vert _{H_{(0,}^{1/2}\left(  I\right)  }^{2}+\left\vert u\right\vert
_{H_{,0)}^{1/2}\left(  I\right)  }^{2}\right)  ^{1/2}, & \left\Vert
u\right\Vert _{H_{00}^{1/2}\left(  I\right)  }:=\left(  \left\Vert
u\right\Vert _{H^{1/2}\left(  I\right)  }^{2}+\left\vert u\right\vert
_{H_{00}^{1/2}\left(  I\right)  }^{2}\right)  ^{1/2}.
\end{array}
\label{halfseminorms}%
\end{equation}

\begin{lemma}
\label{LemBKH1/2norm}Let $k\geq4$ be even and $K\in\mathcal{T}$. Let $\left.
\tilde{B}_{k,K}^{\operatorname*{CR}}\right\vert _{\partial K}$ be defined by
(\ref{traceKequal}). Then there exists an absolute constant $C$ depending only
on the shape regularity of $\mathcal{T}$ such that%
\[
\left\Vert \tilde{B}_{k,K}^{\operatorname*{CR}}\right\Vert _{H^{1/2}\left(
\partial K\right)  }\leq C\sqrt{\log\left(  k+1\right)  }.
\]

\end{lemma}%

%TCIMACRO{\TeXButton{Proof}{\proof}}%
%BeginExpansion
\proof
%EndExpansion
We first prove the estimate for the reference element $\widehat{K}$. By
construction (see (\ref{traceKequal})) the function $\tilde{B}_{k,\widehat{K}%
}^{\operatorname*{CR}}$ coincides with $B_{k,\widehat{K}}^{\operatorname*{CR}%
}$ on $\partial\widehat{K}$. Let the vertices of $\widehat{K}$ be numbered
counterclockwise and denoted by $\widehat{\mathbf{z}}_{i}$, $1\leq i\leq3$.
The edge opposite to $\widehat{\mathbf{z}}_{i}$ is $\widehat{E}_{i}=\left[
\widehat{\mathbf{z}}_{i+1},\widehat{\mathbf{z}}_{i-1}\right]  $ (with cyclic
numbering convention $\mathbf{z}_{3+1}:=\mathbf{z}_{1}$ and $\mathbf{z}%
_{1-1}:=\mathbf{z}_{3}$). We choose the pullbacks to $\left[  -1,1\right]  $
by%
\begin{equation}
\phi_{i}\left(  s\right)  :=\widehat{\mathbf{z}}_{i+1}+s\left(  \widehat
{\mathbf{z}}_{i-1}-\widehat{\mathbf{z}}_{i+1}\right)  \quad i=1,2,3
\label{defphii}%
\end{equation}
and observe $\left.  \tilde{B}_{k,\widehat{K}}^{\operatorname*{CR}}\right\vert
_{\widehat{E}_{i}}\circ\phi_{i}=L_{k}$. From \cite{Grisvard85} (see also
\cite[p 1870]{Ainsworth_parker_II}) we deduce that the $H^{1/2}\left(
\partial\widehat{K}\right)  $ norm is equivalent to
\begin{equation}
\left\vert
%TCIMACRO{\TeXButton{l}{\kern-.1em}}%
%BeginExpansion
\kern-.1em%
%EndExpansion
\left\vert
%TCIMACRO{\TeXButton{l}{\kern-.1em}}%
%BeginExpansion
\kern-.1em%
%EndExpansion
\left\vert v\right\vert
%TCIMACRO{\TeXButton{l}{\kern-.1em}}%
%BeginExpansion
\kern-.1em%
%EndExpansion
\right\vert
%TCIMACRO{\TeXButton{l}{\kern-.1em}}%
%BeginExpansion
\kern-.1em%
%EndExpansion
\right\vert _{H^{1/2}\left(  \partial\widehat{K}\right)  }:=\left(  \sum
_{i=1}^{3}\left(  \left\Vert v_{i}\right\Vert _{H^{1/2}\left(  I\right)  }%
^{2}+\left\vert d_{i}\right\vert _{H_{(0,}^{1/2}\left(  I\right)  }%
^{2}\right)  \right)  ^{1/2}, \label{def3norm}%
\end{equation}
where for $v\in H^{1/2}\left(  \partial\widehat{K}\right)  $:%
\begin{equation}
v_{i}:=\left.  v\right\vert _{E_{i}}\circ\phi_{i},\quad d_{i}\left(  s\right)
:=v_{i-1}\left(  s\right)  -v_{i+1}\left(  -s\right)  . \label{defdi}%
\end{equation}
In our application (and even $k$) we have $v_{i-1}\left(  s\right)
=L_{k}\left(  s\right)  =L_{k}\left(  -s\right)  =v_{i+1}\left(  -s\right)  $
so that $d_{i}=0$. We use $\left\Vert L_{k}\right\Vert _{L^{2}\left(
I\right)  }=\sqrt{2/\left(  2k+1\right)  }$ to get%
\begin{equation}
\left\vert
%TCIMACRO{\TeXButton{l}{\kern-.1em}}%
%BeginExpansion
\kern-.1em%
%EndExpansion
\left\vert
%TCIMACRO{\TeXButton{l}{\kern-.1em}}%
%BeginExpansion
\kern-.1em%
%EndExpansion
\left\vert L_{k}\right\vert
%TCIMACRO{\TeXButton{l}{\kern-.1em}}%
%BeginExpansion
\kern-.1em%
%EndExpansion
\right\vert
%TCIMACRO{\TeXButton{l}{\kern-.1em}}%
%BeginExpansion
\kern-.1em%
%EndExpansion
\right\vert _{H^{1/2}\left(  \partial\widehat{K}\right)  }\leq\left(  \frac
{6}{2k+1}+3\left\vert L_{k}\right\vert _{H^{1/2}\left(  I\right)  }%
^{2}\right)  ^{1/2}. \label{Lkhalfway}%
\end{equation}
This integral can be evaluated analytically for $v=L_{k}$ (see Lem.
\ref{Lemhalbnorm}) and we obtain%
\begin{equation}
\left\vert L_{k}\right\vert _{H^{1/2}\left(  I\right)  }=2\left(  \sum
_{\ell=1}^{k}\frac{1}{\ell}\right)  ^{1/2}\leq\sqrt{C\log\left(  k+1\right)
}\qquad\forall k=1,2,\ldots\label{LkHalfEst}%
\end{equation}
for a generic constant $C>0$. This leads to the final estimate on the
reference element:%
\begin{align*}
\left\Vert \tilde{B}_{k,\widehat{K}}^{\operatorname*{CR}}\right\Vert
_{H^{1/2}\left(  \partial\widehat{K}\right)  }  &  \leq C\left\vert
%TCIMACRO{\TeXButton{l}{\kern-.1em}}%
%BeginExpansion
\kern-.1em%
%EndExpansion
\left\vert
%TCIMACRO{\TeXButton{l}{\kern-.1em}}%
%BeginExpansion
\kern-.1em%
%EndExpansion
\left\vert \tilde{B}_{k,\widehat{K}}^{\operatorname*{CR}}\right\vert
%TCIMACRO{\TeXButton{l}{\kern-.1em}}%
%BeginExpansion
\kern-.1em%
%EndExpansion
\right\vert
%TCIMACRO{\TeXButton{l}{\kern-.1em}}%
%BeginExpansion
\kern-.1em%
%EndExpansion
\right\vert _{H^{1/2}\left(  \partial\widehat{K}\right)  }\overset
{\text{(\ref{Lkhalfway}), (\ref{LkHalfEst})}}{\leq}C\left(  \frac{1}%
{2k+1}+\log\left(  k+1\right)  \right)  ^{1/2}\\
&  \leq C\sqrt{\log\left(  k+1\right)  }.
\end{align*}
For a triangle $K\in\mathcal{T}$, let $\phi_{K}:\widehat{K}\rightarrow K$ be
an affine pullback and set $\tilde{B}_{k,\widehat{K}}^{\operatorname*{CR}%
}=\tilde{B}_{k,K}^{\operatorname*{CR}}\circ\phi_{K}$. Then, the transformation
rule for integrals yields%
\begin{align}
\left\Vert \tilde{B}_{k,K}^{\operatorname*{CR}}\right\Vert _{H^{1/2}\left(
\partial K\right)  }  &  =\left(  \left\Vert \tilde{B}_{k,K}%
^{\operatorname*{CR}}\right\Vert _{L^{2}\left(  \partial K\right)  }%
^{2}+\left\vert \tilde{B}_{k,K}^{\operatorname*{CR}}\right\vert _{H^{1/2}%
\left(  \partial K\right)  }^{2}\right)  ^{1/2}\nonumber\\
&  =C\left(  h_{K}\left\Vert \tilde{B}_{k,\widehat{K}}^{\operatorname*{CR}%
}\right\Vert _{L^{2}\left(  \partial\widehat{K}\right)  }^{2}+\left\vert
\tilde{B}_{k,\widehat{K}}^{\operatorname*{CR}}\right\vert _{H^{1/2}\left(
\partial\widehat{K}\right)  }^{2}\right)  ^{1/2}, \label{H1/2normTransform}%
\end{align}
where $C$ only depends on the shape regularity of the mesh. The leads to the
claim.%
%TCIMACRO{\TeXButton{End Proof}{\endproof}}%
%BeginExpansion
\endproof
%EndExpansion

The case of odd $k\geq5$ is considered in the following lemma.

\begin{lemma}
\label{LemBEH1/2norm}Let $k\geq5$ be odd. For $E\in\mathcal{E}_{\Omega}\left(
\mathcal{T}\right)  $, let the function $\tilde{B}_{k,E}^{\operatorname*{CR}}$
be as in the proof of Lemma \ref{Lemtrianglemean}. Then there exists a generic
constant $C$ depending only on the shape regularity of $\mathcal{T}$ such that
for any $K\in\mathcal{T}_{E}$, it holds%
\[
\left\Vert \tilde{B}_{k,E}^{\operatorname*{CR}}\right\Vert _{H^{1/2}\left(
\partial K\right)  }\leq C\sqrt{\log\left(  k+1\right)  }.
\]

\end{lemma}%

%TCIMACRO{\TeXButton{Proof}{\proof}}%
%BeginExpansion
\proof
%EndExpansion
Let $E\in\mathcal{E}_{\Omega}\left(  \mathcal{T}\right)  $ and $K\in
\mathcal{T}_{E}$. Similarly as in (\ref{H1/2normTransform}) we have%
\[
\left\Vert \tilde{B}_{k,E}^{\operatorname*{CR}}\right\Vert _{H^{1/2}\left(
\partial K\right)  }\leq C\left\Vert \tilde{B}_{k,\widehat{E}}%
^{\operatorname*{CR}}\right\Vert _{H^{1/2}\left(  \partial\widehat{K}\right)
}%
\]
with $\tilde{B}_{k,\widehat{E}}^{\operatorname*{CR}}=\tilde{B}_{k,E}%
^{\operatorname*{CR}}\circ\phi_{K}$ and affine pullback $\phi_{K}:\widehat
{K}\rightarrow K$. Number the vertices in $\widehat{K}$ counterclockwise
$\widehat{\mathbf{z}}_{i}$, $1\leq i\leq3$, such that $\widehat{\mathbf{z}%
}_{3}$ is opposite to $\widehat{E}:=\phi_{K}^{-1}\left(  E\right)  $. The edge
in $\partial\widehat{K}$ opposite to $\widehat{\mathbf{z}}_{i}$ is denoted by
$\widehat{E}_{i}$ and this implies $\widehat{E}=\widehat{E}_{3}$. The edgewise
pullbacks $\phi_{i}$ are defined as in (\ref{defphii}). We employ the
equivalence of the $H^{1/2}\left(  \partial\widehat{K}\right)  $ norm with
$\left\vert
%TCIMACRO{\TeXButton{l}{\kern-.1em}}%
%BeginExpansion
\kern-.1em%
%EndExpansion
\left\vert
%TCIMACRO{\TeXButton{l}{\kern-.1em}}%
%BeginExpansion
\kern-.1em%
%EndExpansion
\left\vert \cdot\right\vert
%TCIMACRO{\TeXButton{l}{\kern-.1em}}%
%BeginExpansion
\kern-.1em%
%EndExpansion
\right\vert
%TCIMACRO{\TeXButton{l}{\kern-.1em}}%
%BeginExpansion
\kern-.1em%
%EndExpansion
\right\vert _{H^{1/2}\left(  \partial\widehat{K}\right)  }$ (see
(\ref{def3norm})) and obtain%
\[
\left\Vert \tilde{B}_{k,\widehat{E}}^{\operatorname*{CR}}\right\Vert
_{H^{1/2}\left(  \partial\widehat{K}\right)  }\leq C\left(  \sum_{i=1}%
^{3}\left(  \left\Vert v_{i}\right\Vert _{H^{1/2}\left(  I\right)  }%
^{2}+\left\vert d_{i}\right\vert _{H_{(0,}^{1/2}\left(  I\right)  }%
^{2}\right)  \right)  ^{1/2}%
\]
with%
\[
v_{i}:=\left.  \tilde{B}_{k,\widehat{E}}^{\operatorname*{CR}}\right\vert
_{\widehat{E}_{i}}\circ\phi_{i}\quad\text{and\quad}d_{i}\left(  s\right)
:=v_{i-1}\left(  s\right)  -v_{i+1}\left(  -s\right)  \text{.}%
\]
Note that%
\begin{equation}
v_{1}\left(  s\right)  =L_{k}\left(  -s\right)  ,\quad v_{2}\left(  s\right)
=L_{k}\left(  s\right)  ,\quad v_{3}=L_{k-1}-\tilde{w}_{3} \label{defvineu}%
\end{equation}
with $\tilde{w}_{3}:=L_{k-1}^{\prime}\left(  -1\right)  \tilde{\psi}_{k}%
^{-}-L_{k-1}^{\prime}\left(  1\right)  \tilde{\psi}_{k}^{+}$ and $\tilde{\psi
}_{k}^{\pm}$ as in (\ref{defphiktilde}). The antisymmetry of the Legendre
polynomial for odd $k$ implies $d_{3}=0$ and%
\begin{equation}
\left\Vert \tilde{B}_{k,\widehat{E}}^{\operatorname*{CR}}\right\Vert
_{H^{1/2}\left(  \partial\widehat{K}\right)  }\leq C\left(  \left\Vert
L_{k}\right\Vert _{H^{1/2}\left(  I\right)  }^{2}+\left\Vert L_{k-1}%
\right\Vert _{H^{1/2}\left(  I\right)  }^{2}+\left\Vert \tilde{w}%
_{3}\right\Vert _{H^{1/2}\left(  I\right)  }^{2}+\sum_{i=1}^{2}\left\vert
d_{i}\right\vert _{H_{(0,}^{1/2}\left(  I\right)  }^{2}\right)  ^{1/2}.
\label{3rdterm}%
\end{equation}
The estimates $\left\Vert L_{j}\right\Vert _{H^{1/2}\left(  I\right)  }%
^{2}\leq C\log\left(  k+1\right)  $ for $j\in\left\{  k-1,k\right\}  $ follow
from (\ref{LkHalfEst}). For the last term, we employ%
\begin{align*}
\left\vert d_{i}\right\vert _{H_{(0,}^{1/2}\left(  I\right)  }  &
\overset{\text{(\ref{defdi}), (\ref{defvineu})}}{\leq}\left\vert L_{k}%
+L_{k-1}\right\vert _{H_{(0,}^{1/2}\left(  I\right)  }+\left\vert \tilde
{w}_{3}\right\vert _{H_{(0,}^{1/2}\left(  I\right)  }\\
&  \overset{\text{(\ref{halfseminorms})}}{\leq}\left\vert L_{k}+L_{k-1}%
\right\vert _{H_{(0,}^{1/2}\left(  I\right)  }+\left\vert \tilde{w}%
_{3}\right\vert _{H_{00}^{1/2}\left(  I\right)  }.
\end{align*}
The last term can be estimated by taking into account $\left\vert
L_{k}^{\prime}\left(  \pm1\right)  \right\vert =\binom{k+1}{2}$ (obtained,
e.g., by evaluating and differentiating \cite[18.5.8 for the choice
$\alpha=\beta=0$ and 18.5.10 for the choice $\lambda=1/2.$]{NIST:DLMF} at
$\pm1$):%
\[
\left\Vert \tilde{w}_{3}\right\Vert _{H_{00}^{1/2}\left(  I\right)  }\leq
C\left(  k+1\right)  ^{2}\left(  \left\Vert \tilde{\psi}_{k}^{-}\right\Vert
_{H_{00}^{1/2}\left(  I\right)  }+\left\Vert \tilde{\psi}_{k}^{+}\right\Vert
_{H_{00}^{1/2}\left(  I\right)  }\right)  \overset{\text{\cite[(A.5)]%
{Ainsworth_parker_I}}}{\leq}C.
\]
For the third term in (\ref{3rdterm}), we simply employ $\left\Vert \tilde
{w}_{3}\right\Vert _{H^{1/2}\left(  I\right)  }\leq\left\Vert \tilde{w}%
_{3}\right\Vert _{H_{00}^{1/2}\left(  I\right)  }$ so that%
\[
\left\Vert \tilde{B}_{k,\widehat{E}}^{\operatorname*{CR}}\right\Vert
_{H^{1/2}\left(  \partial\widehat{K}\right)  }\leq C\left(  \sqrt{\log\left(
k+1\right)  }+\left\vert L_{k}+L_{k-1}\right\vert _{H_{(0,}^{1/2}\left(
I\right)  }\right)  .
\]
The last integral can be evaluated analytically: By using the recurrence
relation (cf. \cite[18.9.1]{NIST:DLMF})
\[
L_{k}\left(  x\right)  =\frac{2k-1}{k}xL_{k-1}\left(  x\right)  -\frac{k-1}%
{k}L_{k-2}\left(  x\right)
\]
we obtain%
\begin{equation}
L_{k}\left(  x\right)  +L_{k-1}\left(  x\right)  =\frac{2k-1}{k}\left(
x+1\right)  L_{k-1}\left(  x\right)  -\frac{k-1}{k}\left(  L_{k-1}\left(
x\right)  +L_{k-2}\left(  x\right)  \right)  . \label{recLeg}%
\end{equation}
This and the orthogonality properties of the Legendre polynomials lead to%
\begin{align}
I_{k}  &  :=\left\vert L_{k}+L_{k-1}\right\vert _{H_{(0,}^{1/2}\left(
I\right)  }^{2}=\int_{-1}^{1}\frac{\left(  \frac{2k-1}{k}\left(  x+1\right)
L_{k-1}\left(  x\right)  -\frac{k-1}{k}\left(  L_{k-1}\left(  x\right)
+L_{k-2}\left(  x\right)  \right)  \right)  ^{2}}{x+1}dx\nonumber\\
&  =\left(  \frac{2k-1}{k}\right)  ^{2}\int_{-1}^{1}\left(  x+1\right)
L_{k-1}^{2}\left(  x\right)  dx-\frac{4\left(  k-1\right)  }{k^{2}}+\left(
\frac{k-1}{k}\right)  ^{2}I_{k-1}. \label{int2ndline}%
\end{align}
The integral in (\ref{int2ndline}) will be evaluated in
(\ref{firstLegendreIntegral}). We obtain%
\[
\left(  \frac{2k-1}{k}\right)  ^{2}\int_{-1}^{1}\left(  x+1\right)
L_{k-1}^{2}\left(  x\right)  dx=\frac{2\left(  2k-1\right)  }{k^{2}}%
\]
and the explicit recursion formula%
\[
I_{k}=\frac{2}{k^{2}}+\left(  \frac{k-1}{k}\right)  ^{2}I_{k-1}%
\]
with starting value%
\[
I_{1}=\int_{-1}^{1}\frac{\left(  L_{1}\left(  x\right)  +L_{0}\left(
x\right)  \right)  ^{2}}{x+1}dx=\int_{-1}^{1}\frac{\left(  1+x\right)  ^{2}%
}{x+1}dx=2.
\]
It is easy to verify that $I_{k}=2/k$ solves the recursion and hence,%
\[
\left\vert L_{k}+L_{k-1}\right\vert _{H_{(0,}^{1/2}\left(  I\right)  }%
=\sqrt{2/k}.
\]
From this, the assertion follows.%
%TCIMACRO{\TeXButton{End Proof}{\endproof}}%
%BeginExpansion
\endproof
%EndExpansion

\section{Analytic evaluation of some integrals involving Legendre
polynomials\label{LegPoly}}

In the proof of Lemma \ref{LemProject} some integrals over Legendre
polynomials appear and we present here their explicit evaluation.

\begin{lemma}
\label{LemExInt}For $k\geq0$, it holds%
\begin{equation}
\int_{-1}^{1}L_{k}^{2}\left(  t\right)  =\int_{-1}^{1}\left(  t+1\right)
L_{k}^{2}\left(  t\right)  dt=\frac{2}{2k+1} \label{firstLegendreIntegral}%
\end{equation}

and%
\begin{equation}
\int_{-1}^{1}\left(  L_{k}^{\prime}\left(  t\right)  \right)  ^{2}dt=\int
_{-1}^{1}\left(  t+1\right)  \left(  L_{k}^{\prime}\left(  t\right)  \right)
^{2}dt=k\left(  k+1\right)  . \label{secondLegendreIntegral}%
\end{equation}

\end{lemma}

%

%TCIMACRO{\TeXButton{Proof}{\proof}}%
%BeginExpansion
\proof
%EndExpansion
The relation $\int_{-1}^{1}L_{k}^{2}\left(  t\right)  dt=2/\left(
2k+1\right)  $ follows from \cite[7.221(1)]{gradstein}.

For $k=0$, these relations follows from $L_{0}\left(  t\right)  =1$. Let
$k\geq1$. The recurrence relation in \cite[18.9.1, Table 18.9.1]{NIST:DLMF} imply%

\begin{equation}
\left(  t+1\right)  L_{k}\left(  t\right)  =\frac{k+1}{2k+1}L_{k+1}\left(
t\right)  +L_{k}\left(  t\right)  +\frac{k}{2k+1}L_{k-1}\left(  t\right)  .
\label{recLeg2}%
\end{equation}
Substituting $\left(  t+1\right)  L_{k}\left(  t\right)  $ under the integral
in (\ref{firstLegendreIntegral}) by this and taking into account the
orthogonality relations of the Legendre polynomials leads to%
\[
\int_{-1}^{1}\left(  t+1\right)  L_{k}^{2}\left(  t\right)  dt=\int_{-1}%
^{1}L_{k}^{2}\left(  t\right)  dt=\frac{2}{2k+1}.
\]
For the second integral (\ref{secondLegendreIntegral}) we employ integration
by parts:%
\begin{equation}
\int_{-1}^{1}\left(  t+1\right)  \left(  L_{k}^{\prime}\left(  t\right)
\right)  ^{2}dt=\left.  \left(  t+1\right)  L_{k}^{\prime}\left(  t\right)
L_{k}\left(  t\right)  \right\vert _{-1}^{1}-\int_{-1}^{1}g_{k}\left(
t\right)  L_{k}\left(  t\right)  dt\quad\text{for }g\left(  t\right)
:=\left(  \left(  t+1\right)  L_{k}^{\prime}\left(  t\right)  \right)
^{\prime}. \label{Lkprimeest}%
\end{equation}
Since $g\in\mathbb{P}_{k-1}$ the orthogonality properties of Legendre
polynomials imply that the integral in the right-hand side of
(\ref{Lkprimeest}) is zero. By using $L_{k}\left(  \pm1\right)  =\left(
\pm1\right)  ^{k}$ (cf. \cite[Table 18.6.1]{NIST:DLMF}) and $L_{k}^{\prime
}\left(  \pm1\right)  =\left(  \pm1\right)  ^{k+1}\binom{k+1}{2}$ (cf.
\cite[combine 18.9.15 with Table 18.6.1]{NIST:DLMF}) we get%
\[
\int_{-1}^{1}\left(  t+1\right)  \left(  L_{k}^{\prime}\left(  t\right)
\right)  ^{2}dt=2L_{k}^{\prime}\left(  1\right)  L_{k}\left(  1\right)
=k\left(  k+1\right)  .
\]
Finally%
\[
\int_{-1}^{1}\left(  L_{k}^{\prime}\left(  t\right)  \right)  ^{2}dt=\left.
L_{k}^{\prime}\left(  t\right)  L_{k}\left(  t\right)  \right\vert _{t=-1}%
^{1}-\int_{-1}^{1}L_{k}^{\prime\prime}\left(  t\right)  L_{k}\left(  t\right)
dt.
\]
The last integral is zero due the orthogonality of the Legendre polynomials.
The endpoint values of $L_{k}$ and $L_{k}^{\prime}$ lead to the assertion.%
%TCIMACRO{\TeXButton{End Proof}{\endproof}}%
%BeginExpansion
\endproof
%EndExpansion

In the final part of this section, we will compute the value of $\left\vert
L_{k}\right\vert _{H^{1/2}\left(  I\right)  }$ explicitly. We set%

\[
I_{k}\left(  s\right)  :=\left\vert L_{k}\right\vert _{H^{s}\left(  I\right)
}^{2}=\int_{-1}^{1}\int_{-1}^{1}\frac{\left(  L_{k}\left(  x\right)
-L_{k}\left(  y\right)  \right)  ^{2}}{\left\vert x-y\right\vert ^{1+2s}%
}dydx.
\]

\begin{lemma}
\label{Lemhalbnorm}It holds%
\[
I_{k}\left(  1/2\right)  =4\sum_{\ell=1}^{k}\frac{1}{\ell}.
\]

\end{lemma}%

%TCIMACRO{\TeXButton{Proof}{\proof}}%
%BeginExpansion
\proof
%EndExpansion
We write%
\[
\left(  L_{k}\left(  y\right)  -L_{k}\left(  x\right)  \right)  ^{2}%
=\sum_{\ell=2}^{2k}\frac{\kappa_{\ell}\left(  x\right)  }{\ell!}\left(
y-x\right)  ^{\ell}%
\]
with%
\begin{align*}
\kappa_{\ell}\left(  x\right)   &  :=\left.  \frac{d^{\ell}}{dy^{\ell}}\left(
\left(  L_{k}\left(  y\right)  -L_{k}\left(  x\right)  \right)  ^{2}\right)
\right\vert _{y\leftarrow x}\\
&  =\sum_{r=0}^{\ell}\binom{\ell}{r}\left.  \left(  \frac{d^{r}}{dy^{r}%
}\left(  L_{k}\left(  y\right)  -L_{k}\left(  x\right)  \right)  \right)
\frac{d^{\ell-r}}{dy^{\ell-r}}\left(  L_{k}\left(  y\right)  -L_{k}\left(
x\right)  \right)  \right\vert _{y\leftarrow x}\\
&  =\sum_{r=1}^{\ell-1}\binom{\ell}{r}\left.  \left(  \frac{d^{r}}{dy^{r}%
}\left(  L_{k}\left(  y\right)  -L_{k}\left(  x\right)  \right)  \right)
\frac{d^{\ell-r}}{dy^{\ell-r}}\left(  L_{k}\left(  y\right)  -L_{k}\left(
x\right)  \right)  \right\vert _{y\leftarrow x}\\
&  =\sum_{r=1}^{\ell-1}\binom{\ell}{r}L_{k}^{\left(  r\right)  }\left(
x\right)  L_{k}^{\left(  \ell-r\right)  }\left(  x\right)  =\left(  L_{k}%
^{2}\right)  ^{\left(  \ell\right)  }\left(  x\right)  -2L_{k}\left(
x\right)  L_{k}^{\left(  \ell\right)  }\left(  x\right)  .
\end{align*}
Hence,%
\[
I_{k}\left(  s\right)  :=I_{k}^{\operatorname*{I}}\left(  s\right)
-I_{k}^{\operatorname*{II}}\left(  s\right)
\]
with%
\begin{align*}
I_{k}^{\operatorname*{I}}\left(  s\right)   &  =\sum_{\ell=2}^{2k}\frac
{1}{\ell!}\int_{-1}^{1}\int_{-1}^{1}\frac{\left(  L_{k}^{2}\right)  ^{\left(
\ell\right)  }\left(  x\right)  \left(  y-x\right)  ^{\ell}}{\left\vert
x-y\right\vert ^{1+2s}}dydx,\\
I_{k}^{\operatorname*{II}}\left(  s\right)   &  =2\sum_{\ell=2}^{2k}\frac
{1}{\ell!}\int_{-1}^{1}\int_{-1}^{1}\frac{L_{k}\left(  x\right)
L_{k}^{\left(  \ell\right)  }\left(  x\right)  \left(  y-x\right)  ^{\ell}%
}{\left\vert x-y\right\vert ^{1+2s}}dydx.
\end{align*}
We perform the integration with respect to $y$ explicitly and get%
\begin{align*}
I_{k}^{\operatorname*{I}}\left(  s\right)   &  =\sum_{\ell=2}^{2k}\frac
{1}{\ell!}\int_{-1}^{1}\left(  L_{k}^{2}\right)  ^{\left(  \ell\right)
}\left(  x\right)  w_{\ell,s}\left(  x\right)  dx,\\
I_{k}^{\operatorname*{II}}\left(  s\right)   &  =2\sum_{\ell=2}^{2k}\frac
{1}{\ell!}\int_{-1}^{1}L_{k}\left(  x\right)  L_{k}^{\left(  \ell\right)
}\left(  x\right)  w_{\ell,s}\left(  x\right)  dx
\end{align*}
with%
\[
w_{\ell,s}\left(  x\right)  :=\frac{\left(  -1\right)  ^{\ell}\left(
1+x\right)  ^{\ell-2s}+\left(  1-x\right)  ^{\ell-2s}}{\ell-2s}.
\]
From now on we restrict to the case $s=1/2$. Since $w_{\ell,1/2}L_{k}^{\left(
\ell\right)  }$ is a polynomial of maximal degree $k-1$ the second integral
vanishes: $I_{k}^{\operatorname*{II}}\left(  1/2\right)  =0$. We apply
recursively integration by parts and use the orthogonality of the Legendre
polynomials to obtain%
\[
I_{k}\left(  1/2\right)  =I_{k}^{\operatorname*{I}}\left(  1/2\right)
=\sum_{\ell=2}^{2k}\sum_{m=0}^{2k-\ell}\frac{\left(  -1\right)  ^{m}2^{m+2}%
}{\left(  m+\ell\right)  \left(  m+\ell-1\right)  }\frac{1}{\left(
m+1\right)  !}\left(  L_{k}^{2}\right)  ^{\left(  m+1\right)  }\left(
1\right)  .
\]
We interchange the ordering of the summation, introduce the new variable
$t=\ell+m-2$, and obtain%
\begin{align*}
I_{k}\left(  1/2\right)   &  =\sum_{m=0}^{2k-2}\sum_{\ell=2}^{2k-m}%
\frac{\left(  -1\right)  ^{m}2^{m+2}}{\left(  m+\ell\right)  \left(
m+\ell-1\right)  }\frac{1}{\left(  m+1\right)  !}\left(  L_{k}^{2}\right)
^{\left(  m+1\right)  }\left(  1\right) \\
&  =\sum_{m=0}^{2k-2}\frac{\left(  -1\right)  ^{m}2^{m+2}}{\left(  m+1\right)
!}\left(  L_{k}^{2}\right)  ^{\left(  m+1\right)  }\left(  1\right)
\sum_{t=m}^{2k-2}\frac{1}{\left(  t+2\right)  \left(  t+1\right)  }.
\end{align*}
By using a telescoping sum argument, it is easy to verify that the inner sum
equals%
\[
\sum_{t=m}^{2k-2}\frac{1}{\left(  t+2\right)  \left(  t+1\right)  }=\sum
_{t=m}^{2k-2}\left(  \frac{1}{t+1}-\frac{1}{t+2}\right)  =\frac{1}{m+1}%
-\frac{1}{2k}.
\]
Hence,%
\begin{equation}
I_{k}\left(  1/2\right)  =\frac{1}{k}I_{k}^{\operatorname*{III}}\left(
1/2\right)  -2I_{k}^{\operatorname*{IV}}\left(  1/2\right)  , \label{Iksplit2}%
\end{equation}
for%
\begin{align*}
I_{k}^{\operatorname*{III}}\left(  1/2\right)   &  :=-\sum_{m=0}^{2k-2}%
\frac{\left(  -1\right)  ^{m}2^{m+1}}{\left(  m+1\right)  !}\left(  L_{k}%
^{2}\right)  ^{\left(  m+1\right)  }\left(  1\right)  ,\\
I_{k}^{\operatorname*{IV}}\left(  1/2\right)   &  :=\sum_{m=0}^{2k-2}%
\frac{\left(  -1\right)  ^{m+1}2^{m+1}}{\left(  m+1\right)  !\left(
m+1\right)  }\left(  L_{k}^{2}\right)  ^{\left(  m+1\right)  }\left(
1\right)  .
\end{align*}
For $I_{k}^{\operatorname*{III}}$ and $k\geq1$ we get%
\begin{align}
I_{k}^{\operatorname*{III}}\left(  1/2\right)   &  =\sum_{m=1}^{2k-1}%
\frac{\left(  -1\right)  ^{m}2^{m}}{m!}\left(  L_{k}^{2}\right)  ^{\left(
m\right)  }\left(  1\right) \nonumber\\
&  =-L_{k}^{2}\left(  1\right)  -\frac{2^{2k}}{\left(  2k\right)  !}\left(
L_{k}^{2}\right)  ^{\left(  2k\right)  }\left(  1\right)  +\sum_{m=0}%
^{2k}\frac{\left(  -1\right)  ^{m}2^{m}}{m!}\left(  L_{k}^{2}\right)
^{\left(  m\right)  }\left(  1\right)  . \label{IkIII2}%
\end{align}
We use the endpoint formula $L_{k}^{\left(  m\right)  }\left(  \pm1\right)
=0$ for $m>k$ and for $m\in\left\{  0,1,\ldots,k\right\}  :$%
\[
L_{k}^{\left(  m\right)  }\left(  \pm1\right)  \overset{\text{\cite[22.5.37,
22.4.2]{Abramowitz}}}{=}\left(  \pm1\right)  ^{k+m}\left(  2m-1\right)
!!\binom{k+m}{k-m}=\left(  \pm1\right)  ^{k+m}\frac{1}{\left(  2m\right)
!!}\frac{\left(  k+m\right)  !}{\left(  k-m\right)  !}.
\]
This leads to $L_{k}^{2}\left(  1\right)  =1$ and the Leibniz rule for
differentiation yields%
\begin{align*}
\frac{2^{2k}}{\left(  2k\right)  !}\left(  L_{k}^{2}\right)  ^{\left(
2k\right)  }\left(  1\right)   &  =\frac{2^{2k}}{\left(  2k\right)  !}%
\sum_{\ell=0}^{2k}\binom{2k}{\ell}L_{k}^{\left(  \ell\right)  }\left(
1\right)  L_{k}^{\left(  2k-\ell\right)  }\left(  1\right)  =\frac{2^{2k}%
}{\left(  2k\right)  !}\binom{2k}{k}\left(  L_{k}^{\left(  k\right)  }\left(
1\right)  \right)  ^{2}\\
&  =\frac{2^{2k}}{\left(  2k\right)  !}\binom{2k}{k}\left(  \frac{\left(
2k\right)  !}{\left(  2k\right)  !!}\right)  ^{2}=\frac{2^{2k}}{\left(
2k\right)  !}\frac{\left(  2k\right)  !}{\left(  k!\right)  ^{2}}\frac{\left(
2k\right)  !^{2}}{2^{2k}\left(  k\right)  !^{2}}=\frac{\left(  2k\right)
!^{2}}{\left(  k\right)  !^{4}}.
\end{align*}
The last sum in (\ref{IkIII2}) is the Taylor expansion of $L_{k}^{2}$ about
$x=1$, evaluated at $x=-1$, i.e.,%
\[
\sum_{m=0}^{2k}\frac{\left(  -1\right)  ^{m}2^{m}}{m!}\left(  L_{k}%
^{2}\right)  ^{\left(  m\right)  }\left(  1\right)  =\left(  L_{k}^{2}\right)
\left(  -1\right)  =1.
\]
Hence,%
\begin{equation}
I_{k}^{\operatorname*{III}}\left(  1/2\right)  =-\frac{\left(  2k\right)
!^{2}}{\left(  k\right)  !^{4}}. \label{IkIIIfinal}%
\end{equation}
For the quantity $I_{k}^{\operatorname*{IV}}$ we obtain by similar arguments%
\begin{align*}
I_{k}^{\operatorname*{IV}}\left(  1/2\right)   &  =\sum_{m=1}^{2k-1}%
\frac{\left(  -1\right)  ^{m}2^{m}}{m!m}\left(  L_{k}^{2}\right)  ^{\left(
m\right)  }\left(  1\right) \\
&  =-\frac{2^{2k-1}}{\left(  2k\right)  !k}\left(  L_{k}^{2}\right)  ^{\left(
2k\right)  }\left(  1\right)  -\int_{-1}^{1}\frac{1}{s-1}\sum_{m=1}^{2k}%
\frac{1}{m!}\left(  L_{k}^{2}\right)  ^{\left(  m\right)  }\left(  1\right)
\left(  s-1\right)  ^{m}ds\\
&  =-\frac{1}{2k}\frac{\left(  2k\right)  !^{2}}{\left(  k\right)  !^{4}}%
-\int_{-1}^{1}\frac{L_{k}^{2}\left(  s\right)  -L_{k}^{2}\left(  1\right)
}{s-1}ds.
\end{align*}
The combination of this with (\ref{Iksplit2}), (\ref{IkIIIfinal}) leads to%
\[
I_{k}\left(  1/2\right)  =2\int_{-1}^{1}\frac{L_{k}^{2}\left(  s\right)
-1}{s-1}ds=2\int_{-1}^{1}\frac{L_{k}\left(  s\right)  -1}{s-1}\left(
L_{k}\left(  s\right)  +1\right)  ds.
\]
Since $\frac{L_{k}\left(  s\right)  -1}{s-1}$ is a polynomial of maximal
degree $k-1$ the orthogonality of $L_{k}$ leads to%
\begin{equation}
I_{k}\left(  1/2\right)  =2\int_{-1}^{1}\frac{L_{k}\left(  s\right)  -1}%
{s-1}ds. \label{Iklinint}%
\end{equation}
We employ the recursion formula in \cite[Table 18.9.1]{NIST:DLMF} for $k\geq2$%
\[
L_{k}\left(  s\right)  =\frac{2k-1}{k}sL_{k-1}\left(  s\right)  -\frac{k-1}%
{k}L_{k-2}\left(  s\right)
\]
so that
\begin{align}
I_{k}\left(  1/2\right)   &  =2\int_{-1}^{1}\frac{\frac{2k-1}{k}\left(
sL_{k-1}\left(  s\right)  -1\right)  -\frac{k-1}{k}\left(  L_{k-2}\left(
s\right)  -1\right)  }{s-1}ds\nonumber\\
&  =\frac{2k-1}{k}2\int_{-1}^{1}\left(  L_{k-1}\left(  s\right)
+\frac{L_{k-1}\left(  s\right)  -1}{s-1}\right)  ds-\frac{k-1}{k}%
I_{k-2}\left(  1/2\right) \nonumber\\
&  =\frac{2k-1}{k}I_{k-1}\left(  1/2\right)  -\frac{k-1}{k}I_{k-2}\left(
1/2\right)  . \label{recfinal}%
\end{align}
From (\ref{Iklinint}) and $L_{0}\left(  s\right)  =1$, $L_{1}\left(  s\right)
=s$, $L_{2}\left(  s\right)  =\left(  3s^{2}-1\right)  /2$ we get%
\begin{equation}%
\begin{array}
[c]{l}%
I_{0}\left(  1/2\right)  =0,\\
I_{1}\left(  1/2\right)  =2\int_{-1}^{1}\frac{L_{1}\left(  s\right)  -1}%
{s-1}ds=4,\\
I_{2}\left(  1/2\right)  =2\int_{-1}^{1}\frac{L_{2}\left(  s\right)  -1}%
{s-1}ds=6.
\end{array}
\label{initvalrec}%
\end{equation}
Now, it is easy to verify that $I_{k}=4\sum_{\ell=1}^{k}\frac{1}{\ell}$
satisfies the recursion (\ref{recfinal}) and the initial value conditions
(\ref{initvalrec}).%
%TCIMACRO{\TeXButton{End Proof}{\endproof}}%
%BeginExpansion
\endproof
%EndExpansion

\section{Norm equivalence for Crouzeix-Raviart spaces\label{SecNormEq}}

It is well known that for $V:=H_{0}^{1}\left(  \Omega\right)
+\operatorname*{CR}_{k,0}\left(  \mathcal{T}\right)  $, the norms $\left(
\left\Vert \nabla_{\mathcal{T}}u\right\Vert _{L^{2}\left(  \Omega\right)
}^{2}+\left\Vert u\right\Vert _{L^{2}\left(  \Omega\right)  }^{2}\right)
^{1/2}$ and $\left\Vert \nabla_{\mathcal{T}}u\right\Vert _{L^{2}\left(
\Omega\right)  }$ are equivalent. In this section, we state estimates for the
constants in these equivalencies -- the proof is a repetition of the
well-known arguments for the case $k=1$ (see, e.g., \cite[Lem. 36.6]%
{ErnGuermondII}).

\begin{theorem}
\label{TheoDiscFried}There exists a constant $C>0$ depending only on the
shape-regularity of the mesh and the domain $\Omega$ such that%
\[
\left\Vert u\right\Vert _{H^{1}\left(  \mathcal{T}\right)  }\leq\left(
\left\Vert \nabla_{\mathcal{T}}u\right\Vert _{L^{2}\left(  \Omega\right)
}^{2}+\left\Vert u\right\Vert _{L^{2}\left(  \Omega\right)  }^{2}\right)
^{1/2}\leq C\left\Vert u\right\Vert _{H^{1}\left(  \mathcal{T}\right)  }%
\quad\forall u\in V.
\]
In particular $C$ is independent of the polynomial degree $k\geq1$ and the
mesh size $h_{\mathcal{T}}$.
\end{theorem}

%

%TCIMACRO{\TeXButton{Proof}{\proof}}%
%BeginExpansion
\proof
%EndExpansion
We prove this result only under the regularity assumption that the Poisson
problem:%
\[
\text{find }\phi\in H_{0}^{1}\left(  \Omega\right)  \quad\text{s.t.\quad
}\left(  \nabla\phi,\nabla v\right)  _{\mathbf{L}^{2}\left(  \Omega\right)
}=\left(  f,v\right)  _{L^{2}\left(  \Omega\right)  }\quad\forall v\in
H_{0}^{1}\left(  \Omega\right)
\]
is $H^{2}$ regular. For less regularity we refer to \cite[Lem. 36.6]%
{ErnGuermondII}. For $u\in V$, we have%
\begin{equation}
\left\Vert u\right\Vert _{L^{2}\left(  \Omega\right)  }=\sup_{v\in
L^{2}\left(  \Omega\right)  \backslash\left\{  0\right\}  }\frac{\left(
u,v\right)  _{L^{2}\left(  \Omega\right)  }}{\left\Vert v\right\Vert
_{L^{2}\left(  \Omega\right)  }}. \label{defuL2}%
\end{equation}
For $v\in L^{2}\left(  \Omega\right)  $, there exists some $\mathbf{w}%
\in\mathbf{H}^{1}\left(  \Omega\right)  $ such that $\operatorname*{div}%
\mathbf{w}=v$ and $\left\Vert \mathbf{w}\right\Vert _{\mathbf{H}^{1}\left(
\Omega\right)  }\leq C_{\Omega}\left\Vert v\right\Vert _{L^{2}\left(
\Omega\right)  }$ for a constant $C_{\Omega}$ which only depends on $\Omega$.
Hence,%
\[
\left(  u,v\right)  _{L^{2}\left(  \Omega\right)  }=\left(
u,\operatorname*{div}\mathbf{w}\right)  _{L^{2}\left(  \Omega\right)
}=-\left(  \nabla_{\mathcal{T}}u,\mathbf{w}\right)  _{L^{2}\left(
\Omega\right)  }+\sum_{K\in\mathcal{T}}\int_{\partial K}\left\langle
\mathbf{w},\mathbf{n}_{K}\right\rangle u,
\]
where $\mathbf{n}_{K}$ is the unit normal vector pointing to the exterior of
$K$. Next, we rewrite the sum over the triangle boundaries as a sum over the
edges. For $E\in\mathcal{E}_{\Omega}\left(  \mathcal{T}\right)  $ we fix the
direction of a unit vector $\mathbf{n}_{E}$ which is orthogonal to $E$ and for
$E\in\mathcal{E}_{\partial\Omega}\left(  \mathcal{T}\right)  $ let
$\mathbf{n}_{E}$ denote the unit vector, orthogonal to $E$, pointing to the
exterior of $\Omega$. Then%
\[
\left(  u,v\right)  _{L^{2}\left(  \Omega\right)  }=\left(
u,\operatorname*{div}\mathbf{w}\right)  _{L^{2}\left(  \Omega\right)
}=-\left(  \nabla_{\mathcal{T}}u,\mathbf{w}\right)  _{L^{2}\left(
\Omega\right)  }-\sum_{E\in\mathcal{E}_{\Omega}\left(  \mathcal{T}\right)
}\int_{E}\left\langle \mathbf{w},\mathbf{n}_{E}\right\rangle \left[  u\right]
_{E}+\sum_{E\in\mathcal{E}_{\partial\Omega}\left(  \mathcal{T}\right)  }%
\int_{E}\left\langle \mathbf{w},\mathbf{n}_{E}\right\rangle u.
\]
Let $K_{E}\in\mathcal{T}_{E}$ be fixed and let $\mathbf{q}_{E}\in\left(
\mathbb{P}_{0}\left(  K_{E}\right)  \right)  ^{2}$ be the function with
constant value $\frac{1}{\left\vert E\right\vert }\int_{E}\mathbf{w}$. The
orthogonality conditions of the Crouzeix-Raviart elements across edges (see
(\ref{PCR0RB})) imply%
\begin{align}
\left(  u,v\right)  _{L^{2}\left(  \Omega\right)  }  &  =-\left(
\nabla_{\mathcal{T}}u,\mathbf{w}\right)  _{L^{2}\left(  \Omega\right)  }%
-\sum_{E\in\mathcal{E}_{\Omega}\left(  \mathcal{T}\right)  }\int
_{E}\left\langle \mathbf{w}-\mathbf{q}_{E},\mathbf{n}_{E}\right\rangle \left[
u\right]  _{E}+\sum_{E\in\mathcal{E}_{\partial\Omega}\left(  \mathcal{T}%
\right)  }\int_{E}\left\langle \mathbf{w}-\mathbf{q}_{E},\mathbf{n}%
_{E}\right\rangle u\label{L2estbroken}\\
&  \leq\left\Vert \nabla_{\mathcal{T}}u\right\Vert _{\mathbf{L}^{2}\left(
\Omega\right)  }\left\Vert \mathbf{w}\right\Vert _{\mathbf{L}^{2}\left(
\Omega\right)  }\nonumber\\
&  +\sum_{E\in\mathcal{E}_{\Omega}\left(  \mathcal{T}\right)  }\left\Vert
\left[  u\right]  _{E}\right\Vert _{L^{2}\left(  E\right)  }\left\Vert
\mathbf{w}-\mathbf{q}_{E}\right\Vert _{\mathbf{L}^{2}\left(  E\right)  }%
+\sum_{E\in\mathcal{E}_{\partial\Omega}\left(  \mathcal{T}\right)  }\left\Vert
u\right\Vert _{L^{2}\left(  E\right)  }\left\Vert \mathbf{w}-\mathbf{q}%
_{E}\right\Vert _{\mathbf{L}^{2}\left(  E\right)  }.\nonumber
\end{align}
We employ first a weighted trace inequality (see, e.g., \cite[Lem.
12.15]{ErnGuermondI}) and then a Poincar\'{e}-Steklov estimate (see, e.g.,
\cite[(12.17) for $p=2$ and $s=1$.]{ErnGuermondII}) to get for $h_{E}%
:=\left\vert E\right\vert $%
\begin{equation}
\left\Vert \mathbf{w}-\mathbf{q}_{E}\right\Vert _{\mathbf{L}^{2}\left(
E\right)  }\leq C\left(  h_{E}^{-1/2}\left\Vert \mathbf{w}-\mathbf{q}%
_{E}\right\Vert _{\mathbf{L}^{2}\left(  K_{E}\right)  }+h_{E}^{1/2}\left\Vert
\nabla\left(  \mathbf{w}-\mathbf{q}_{E}\right)  \right\Vert _{\mathbb{L}%
^{2}\left(  K_{E}\right)  }\right)  \leq Ch_{K_{E}}^{1/2}\left\Vert
\nabla\mathbf{w}\right\Vert _{\mathbb{L}^{2}\left(  K_{E}\right)  },
\label{wbroken}%
\end{equation}
where $C$ only depends on the shape-regularity of the mesh.

Next we estimate the jump of $u$ across $E$. For $E\in\mathcal{E}_{\Omega
}\left(  \mathcal{T}\right)  $, we define $u_{E}\in\mathbb{P}_{0}\left(
\mathcal{T}_{E}\right)  $ as the function with constant value $\frac
{1}{\left\vert E\right\vert }\int_{E}\left.  u\right\vert _{K}$ on
$K\in\mathcal{T}_{E}$ and observe $\left[  u_{E}\right]  _{E}=0$. Hence,%
\begin{align*}
\left\Vert \left[  u\right]  _{E}\right\Vert _{L^{2}\left(  E\right)  }  &
=\left\Vert \left[  u-u_{E}\right]  _{E}\right\Vert _{L^{2}\left(  E\right)
}\leq\sum_{K\in\mathcal{T}_{E}}\left\Vert \left.  \left(  u-u_{E}\right)
\right\vert _{K}\right\Vert _{L^{2}\left(  E\right)  }\\
&  \leq\sum_{K\in\mathcal{T}_{E}}\left(  h_{E}^{-1/2}\left\Vert u-u_{E}%
\right\Vert _{L^{2}\left(  K\right)  }+h_{E}^{1/2}\left\Vert \nabla
u\right\Vert _{\mathbf{L}^{2}\left(  K\right)  }\right)  \leq C\sum
_{K\in\mathcal{T}_{E}}h_{K}^{1/2}\left\Vert \nabla u\right\Vert _{\mathbf{L}%
^{2}\left(  K\right)  },
\end{align*}
for a constant $C$ which only depends on the shape-regularity of the mesh. For
$E\in\mathcal{E}_{\partial\Omega}\left(  \mathcal{T}\right)  $ the estimate
$\left\Vert u\right\Vert _{L^{2}\left(  E\right)  }\leq h_{K}^{1/2}\left\Vert
\nabla u\right\Vert _{\mathbf{L}^{2}\left(  K\right)  }$ for $K\in
\mathcal{T}_{E}$ follows in a similar fashion. The combination of
(\ref{L2estbroken}) with (\ref{wbroken}) and the two trace estimates for $u$
leads to%
\begin{align*}
\left(  u,v\right)  _{L^{2}\left(  \Omega\right)  }  &  \leq\left\Vert
\nabla_{\mathcal{T}}u\right\Vert _{\mathbf{L}^{2}\left(  \Omega\right)
}\left\Vert \mathbf{w}\right\Vert _{\mathbf{L}^{2}\left(  \Omega\right)
}+C\sum_{E\in\mathcal{E}\left(  \mathcal{T}\right)  }\left\Vert \nabla
_{\mathcal{T}}u\right\Vert _{\mathbf{L}^{2}\left(  \omega_{E}\right)
}\left\Vert \nabla\mathbf{w}\right\Vert _{\mathbf{L}^{2}\left(  \omega
_{E}\right)  }\\
&  \leq C\left\Vert \nabla_{\mathcal{T}}u\right\Vert _{\mathbf{L}^{2}\left(
\Omega\right)  }\left\Vert \mathbf{w}\right\Vert _{\mathbf{H}^{1}\left(
\Omega\right)  }\leq CC_{\Omega}\left\Vert \nabla_{\mathcal{T}}u\right\Vert
_{\mathbf{L}^{2}\left(  \Omega\right)  }\left\Vert q\right\Vert _{L^{2}\left(
\Omega\right)  }.
\end{align*}
Using this estimate in (\ref{defuL2}) finishes the proof.%
%TCIMACRO{\TeXButton{End Proof}{\endproof}}%
%BeginExpansion
\endproof
%EndExpansion

\bibliographystyle{abbrv}
\bibliography{nlailu}

\end{document}